\newif\iflinenumbers
\linenumbersfalse

\documentclass[a4paper, preprint, 3p]{elsarticle}

\usepackage{natbib}
\usepackage{lmodern}
\usepackage[utf8]{inputenc}
\usepackage[english]{babel}

\usepackage{microtype}  % Improved typography

\usepackage{bm}
\usepackage{algorithm}
\usepackage[noEnd=false]{algpseudocodex}
\usepackage{graphicx}
\usepackage{amsmath}
\usepackage{wrapfig}
\usepackage{caption}
\usepackage{subcaption}
\usepackage{multirow}
\usepackage{booktabs}
\usepackage{mathtools}
\usepackage{csquotes}
\usepackage{xspace}
\usepackage{transparent}
\usepackage{import}
\usepackage{derivative}
\usepackage{sansmath}
\usepackage{amssymb}
\usepackage{amsthm}
\usepackage{subfig}
\usepackage{thmtools}
\usepackage[normalem]{ulem}
\usepackage[export]{adjustbox}

\usepackage{lineno}
\newcommand*\patchAmsMathEnvironmentForLineno[1]{%
  \expandafter\let\csname old#1\expandafter\endcsname\csname #1\endcsname
  \expandafter\let\csname oldend#1\expandafter\endcsname\csname end#1\endcsname
  \renewenvironment{#1}%
     {\linenomath\csname old#1\endcsname}%
     {\csname oldend#1\endcsname\endlinenomath}}%
\newcommand*\patchBothAmsMathEnvironmentsForLineno[1]{%
  \patchAmsMathEnvironmentForLineno{#1}%
  \patchAmsMathEnvironmentForLineno{#1*}}%
\AtBeginDocument{%
\patchBothAmsMathEnvironmentsForLineno{equation}%
\patchBothAmsMathEnvironmentsForLineno{align}%
\patchBothAmsMathEnvironmentsForLineno{flalign}%
\patchBothAmsMathEnvironmentsForLineno{alignat}%
\patchBothAmsMathEnvironmentsForLineno{gather}%
\patchBothAmsMathEnvironmentsForLineno{multline}%
}
\iflinenumbers
\linenumbers %TODO
\fi

\usepackage{todonotes}
\setuptodonotes{caption={}}
\algrenewcommand\textproc{}% Used for small function names

\graphicspath{{figures/}}

\usepackage{tikz}
\usepackage{tikz-3dplot}
\usepackage{pgfplots}
\usetikzlibrary{3d, calc, decorations.pathreplacing, angles, quotes, arrows.meta, positioning}

\date{\today}

\begin{document}

\begin{frontmatter}

\title{Robust and efficient pre-processing techniques for particle-based methods including dynamic boundary generation}

\author[1]{Niklas S. Neher\corref{cor1}}
\ead{neher@hlrs.de}
\cortext[cor1]{Corresponding author}

\author[2]{Erik Faulhaber}

\author[3]{Sven Berger}

\author[5]{Christian Weißenfels}

\author[2]{Gregor J. Gassner}

\author[4]{Michael~Schlottke-Lakemper}

\address[1]{High-Performance Computing Center Stuttgart, University of Stuttgart, Germany}
\address[2]{Department of Mathematics and Computer Science, University of Cologne, Germany}
\address[3]{Institute of Surface Science, Helmholtz--Zentrum hereon, Germany}
\address[4]{High-Performance Scientific Computing, Institute of Mathematics, University of
Augsburg, Germany}
\address[5]{Institute of Materials Resource Management, University of
Augsburg, Germany}

\begin{abstract}
  \iflinenumbers
  \begin{linenumbers}
    % !TeX root = main.tex
Obtaining high-quality particle distributions for stable and accurate particle-based simulations poses significant challenges, especially for complex geometries.
We introduce a preprocessing technique for 2D and 3D geometries, optimized for smoothed particle hydrodynamics (SPH) and other particle-based methods.
Our pipeline begins with the generation of a resolution-adaptive point cloud near the geometry's surface employing a face-based neighborhood search.
This point cloud forms the basis for a signed distance field,
enabling efficient, localized computations near surface regions.
To create an initial particle configuration, we apply a hierarchical winding number method for fast and accurate inside-outside segmentation.
Particle positions are then relaxed using an SPH-inspired scheme, which also serves to pack boundary particles.
This ensures full kernel support and promotes isotropic distributions while preserving the geometry interface.
By leveraging the meshless nature of particle-based methods,
our approach does not require connectivity information and is thus straightforward to integrate into existing particle-based frameworks.
It is robust to imperfect input geometries and memory-efficient without compromising performance.
Moreover, our experiments demonstrate that with increasingly higher resolution, the
resulting particle distribution converges to the exact geometry.

  \end{linenumbers}
  \else
    % !TeX root = main.tex
Obtaining high-quality particle distributions for stable and accurate particle-based simulations poses significant challenges, especially for complex geometries.
We introduce a preprocessing technique for 2D and 3D geometries, optimized for smoothed particle hydrodynamics (SPH) and other particle-based methods.
Our pipeline begins with the generation of a resolution-adaptive point cloud near the geometry's surface employing a face-based neighborhood search.
This point cloud forms the basis for a signed distance field,
enabling efficient, localized computations near surface regions.
To create an initial particle configuration, we apply a hierarchical winding number method for fast and accurate inside-outside segmentation.
Particle positions are then relaxed using an SPH-inspired scheme, which also serves to pack boundary particles.
This ensures full kernel support and promotes isotropic distributions while preserving the geometry interface.
By leveraging the meshless nature of particle-based methods,
our approach does not require connectivity information and is thus straightforward to integrate into existing particle-based frameworks.
It is robust to imperfect input geometries and memory-efficient without compromising performance.
Moreover, our experiments demonstrate that with increasingly higher resolution, the
resulting particle distribution converges to the exact geometry.

  \fi
\end{abstract}

\end{frontmatter}

% !TeX root = ../main.tex
\section{Introduction}
\label{sec:introduction}
Unlike mesh-based methods, particle-based methods represent the continuum using discrete,
unconnected points or particles that carry physical properties and state variables.
In many such methods, e.g., smoothed particle hydrodynamics (SPH),
interactions between particles are mediated by a weighting function or kernel that allows field quantities to be estimated through local summations.
For classical kernel-based approaches such as SPH, the accuracy of the approximation depends strongly on the quality of the particle distribution~\cite{Adami2013, Sun2017, Diehl2015}.
In methods like SPH, a uniform and isotropic distribution—meaning that the particles are arranged without any preferred direction either locally or globally—is essential for ensuring that the kernel interpolation accurately reproduces the continuum field~\cite{Price2012}.
While alternative formulations based on least-squares approximation, such as~\cite{Dilts1999}, can in principle achieve second-order accuracy independently of the particle distribution,
this typically comes at the cost of higher computational and algorithmic complexity.

To obtain a suitable first particle distribution, a common approach is to perform an initial relaxation step with high numerical damping to improve uniformity before the actual simulation, as described in Monaghan et al.~\cite{Monaghan_1994}.
However, this approach depends on the initial configuration being of sufficiently high quality to ensure that the relaxation process remains stable.
While such relaxation procedures can improve accuracy in fluid dynamics, they are not directly applicable in all contexts.
For instance, in structural mechanics the initial particle configuration serves as the reference state for all measurements of deformation and stress~\cite{O_Connor_2021,Sun2019}.
Thus, any changes during relaxation would invalidate subsequent calculations.

To address this issue, the particle distribution is typically optimized during a preprocessing stage prior to the simulation, thereby eliminating the need for a relaxation step.
This prompts the question of which objectives the preprocessing must satisfy.
According to Diehl et al.~\cite{Diehl2015}, an optimal particle configuration method should achieve isotropic particle distributions at both local and global scales.
Furthermore, high interpolation accuracy can only be guaranteed by a uniform lattice configuration.
At the same time, the method should be sufficiently versatile to reproduce any spatial configuration.
This includes geometries represented as numerical datasets, such as triangulated surface meshes or polygonal representations.

To meet these requirements, the preprocessing typically consists of two steps~\cite{Zhu2021,Negi2021,Colagrossi2012}.
First, an initial particle configuration representing the geometry is generated using inside-outside segmentation techniques.
To perform the segmentation, a uniform particle grid can be generated analytically, for example using a Cartesian or hexagonal lattice.
While this yields a perfectly isotropic distribution after segmentation, such configurations are typically not body-fitted.
As a result, non-trivial geometries are not accurately resolved near the surface, particularly in complex regions, leading to discretization errors.
To overcome this, a second step is performed, further adjusting the particle positions.
In the following, we briefly review established approaches for generating the initial particle configuration,
before turning to the second step, commonly referred to as particle packing.

A detailed overview of point cloud generation methods for the initial particle configuration is given by Suchde et al.~\cite{Suchde2022}.
The point cloud is commonly generated by directly placing particles on a lattice structure \cite{Dominguez2011, Zhu2021, Yu2023a}
or by generating particles on a volume element mesh~\cite{Heimbs2011, Akinci2013CouplingES}.
In the former, a particle grid is created that spans the entire bounding box of the geometry,
where subsequent optimization of the bounding box can reduce the number of particles.
Next, the particles are segmented into those inside and outside the geometry.
The segmentation algorithm must be efficient and robust, as generating the initial configuration for large problems (large number of particles) or complex geometries (large number of faces) can be computationally expensive.
One approach is to perform the segmentation via a level-set method~\cite{Fu2019,Zhu2021,Yu2023a,Ji2020,Ji2021},
which typically requires the construction of a background mesh to store the signed distance values at each cell.
However, because the level-set function transitions smoothly to zero at the interface, the distinction between interior and exterior points becomes less pronounced near the surface.
Consequently, in regions with complex geometrical features, resolving steep gradients of the signed distance function demands a highly refined mesh.
To alleviate the computational and memory costs associated with such fine discretizations, the narrow band method~\cite{Gomez2005} is often employed,
whereby the mesh is higher resolved in a limited region around the interface.
Alternatively, background mesh--free approaches such as point-in-polygon or point-in-polyhedron tests can be used.

For the particle packing, there are various approaches to generate an isotropic and body-fitted particle distribution.
Colagrossi et al.~\cite{Colagrossi2012} utilize the inherent properties of SPH for packing.
In their approach, particles are evolved using the momentum equations with a damping term, effectively ``packing'' them.
To confine the inner particles, the geometry is enclosed by analytically generated fixed solid particles.
The interface between the particles to be packed and the fixed particles forms the surface profile of the geometry.
This means that the accuracy of the geometry representation depends on the fixed particles and is therefore suitable only for simple geometries.
Negi et al.~\cite{Negi2021} take a different approach. Instead of segmenting a geometry, a rectangular domain enclosing a geometry is sampled with particles.
In this domain, particles are assigned to boundary particles if they are near the geometry surface represented by boundary nodes.
The boundary particles conform to the surface of the geometry during packing by projecting them onto the surface.
This method requires a set of points that discretizes the surface of the geometry.
The boundary nodes need to be shifted by half a particle spacing,
e.g., if the actual boundary surface is exactly between solid and fluid particles.
This becomes challenging when the geometry has sharp features, as the nodes tend to self-intersect,
and requires an initialization step for the input geometry, as the geometry surface is sampled analytically using nodes.
One of the first methods that works with interpolation of tabulated data and that does not require analytical solutions was proposed by Diehl et al.~\cite{Diehl2015}
and is based on a weighted Voronoi tessellation algorithm.
Fu et al.~\cite{Fu2019} proposed a simpler method that does not require an expensive and complex Voronoi algorithm,
in which they use the inter-particle pressure force for particle relaxation and represent the surface of the geometry by a zero--level-set function,
allowing particles outside the geometry to be fixed accordingly.
Zhu et al.~\cite{Zhu2021} employ a similar method, where particle relaxation is driven by the transport-velocity formulation (TVF)~\cite{Adami2013}.
Since full kernel support must be maintained to achieve a force equilibrium, even though there are no particles outside the geometry, boundary conditions are a major challenge.
To remedy this, Ji et al.~\cite{Ji2021} apply a correction method for a set of fixed boundary particles which then mimics full kernel support.
For this, the boundary particles must be arranged accordingly on the geometry surface, which is not trivial.
In general, generating ghost particles for the boundary conditions can be quite challenging for complex geometries.
Thus, Zhu et al.~\cite{Zhu2021} do not use any ghost particles but only constrain the particles to the surface by updating their positions according to the bounding method by Fu et al.~\cite{Fu2019}.
However, full kernel support is not ensured, which leads to clumping effects on the surface, especially in areas with high surface curvature.
Yu et al.~\cite{Yu2023a} addressed this issue with a static confinement boundary condition where they compute the missing kernel support outside the geometry and add it to the pressure force calculation.

The focus of our work lies on developing a preprocessing method capable of handling arbitrarily complex geometries without the need for manual parameter tuning.
We particularly target biomechanical applications, in which a complex fluid domain must be populated with an isotropic particle distribution.
Furthermore, the method should enable the automated generation of an accurate boundary representation for any preprocessed geometry.
Equally important, the approach must be purely particle-based to allow for seamless integration into existing SPH frameworks, while maintaining robustness and computational efficiency.

In this paper, the initial particle configuration is generated using the hierarchical winding number approach introduced by Jacobson et al.~\cite{Jacobson2013}.
To the best of our knowledge, no other particle-based method or framework has been proposed that leverages the winding number for inside-outside segmentation in this context.
This approach eliminates the need for a background mesh and resolves the geometry exactly,
as the winding number provides a jump discontinuity across the geometry surface.
Additionally, this method is efficient and robust, as will be demonstrated later.

For particle packing, we combine and extend several of the approaches mentioned above.
To constrain particles to the surface, we apply a method analogous to those described in Fu et al.~\cite{Fu2019}, Zhu et al.~\cite{Zhu2021}, and Yu et al.~\cite{Yu2023a}.
A face-based neighborhood search (NHS) is developed to efficiently generate a resolution-dependent and memory-efficient signed distance field (SDF) near the surface.
This SDF is stored in a fixed point cloud localized around the surface.
During the packing process, Shepard interpolation~\cite{Shepard1968} is employed to project the SDF information onto nearby particles, ensuring that they remain constrained to the geometry's surface.
  Beyond the memory efficiency gained by storing signed distance values locally rather than on a global grid, this approach is motivated by two main benefits:
  First, the use of Shepard interpolation eliminates the requirement for global mesh connectivity.
  Information is interpolated locally, which allows the level-set points to be distributed arbitrarily.
  Consequently, we avoid the constraint of a structured grid and the need to sample the entire computational domain.
  This also facilitates integration into existing particle-based frameworks, since standard particle-based NHS algorithms can be employed directly.
  Second, our interpolation strategy does not rely on Gaussian filtering of structured mesh nodes to smooth the level-set function,
  as done for example in \cite{Marrone2010,Yu2023a}.
  Instead, smoothing is inherently achieved through the SPH interpolation itself.
Additionally, we introduce a boundary particle generation method that supports a dynamic packing of the boundary particles.
The approach natively guarantees full kernel support, with the boundary particles distributed uniformly.
Moreover, the boundary particles together with the inner particles form an isotropic distribution while preserving the geometry interface.
These features not only improve the accuracy of geometric representation but also facilitates rapid convergence of the packing process,
leading to a faster stabilization toward the steady state, as can be demonstrated.

To summarize, our approach provides two key advancements.
First, we introduce an automated method for creating boundary particles that are dynamically packed.
Second, we develop a pointcloud-based technique for the level-set method.
Combined with the hierarchical winding number approach, this method offers several advantages, as outlined above.
A detailed description of the entire preprocessing pipeline is provided in the following section.

The proposed techniques are entirely particle-based and can be seamlessly integrated into existing SPH frameworks.
Our implementation utilizes the open-source software package TrixiParticles.jl \cite{Neher2025}, where all algorithms described here are implemented.
We validate the quality of the particle configuration in 2D and 3D based on the particle number density, demonstrating that increased resolution yields convergence toward the exact geometry.
Finally, a performance analysis for each processing step confirms the computational efficiency of the technique.
The numerical setups for all results shown in this paper, including the necessary post-processing code, can be found in our reproducibility repository~\cite{Neher2025reproducibility}.
The paper is organized as follows.
Sec~\ref{sec:overview} gives an overview of the entire preprocessing pipeline.
In Sec.~\ref{sec:sph}, a brief overview of the SPH method and the time integration scheme is included for context.
Sec.~\ref{sec:signed_distance_point_cloud} details the generation of the point cloud representing the SDF.
Sec.~\ref{sec:initial_particle_configuration} describes the inside-outside segmentation method for creating the initial particle configuration.
Sec.~\ref{sec:packing} discusses the packing algorithm to achieve a body-fitted particle distribution.
In Sec.~\ref{sec:results} validation cases are shown for 2D and 3D, and a simulation study is conducted.
Sec.~\ref{sec:performance} provides a performance analysis and in the final section, concluding remarks are given.

\section{Overview of the Preprocessing Pipeline}
\label{sec:overview}
\begin{figure}[!h]
  \centering
  \begin{subfigure}[b]{0.24\textwidth}
    \centering
    \includegraphics[width=0.95\textwidth]{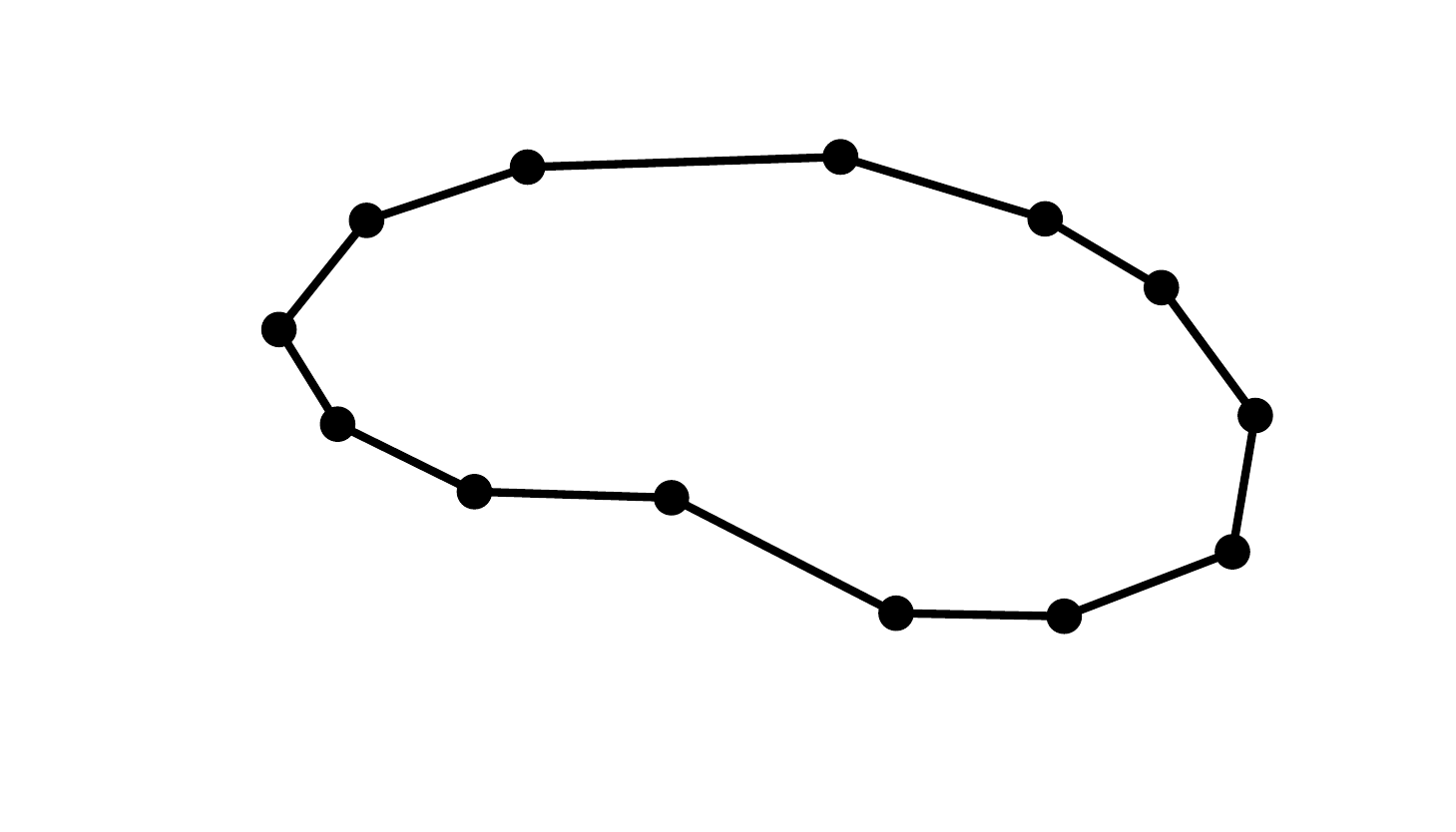}
    \caption{Geometry representation.}
    \label{fig:potato_geometry}
  \end{subfigure}
  \hfill
  \begin{subfigure}[b]{0.28\textwidth}
    \centering
    \includegraphics[width=0.95\textwidth]{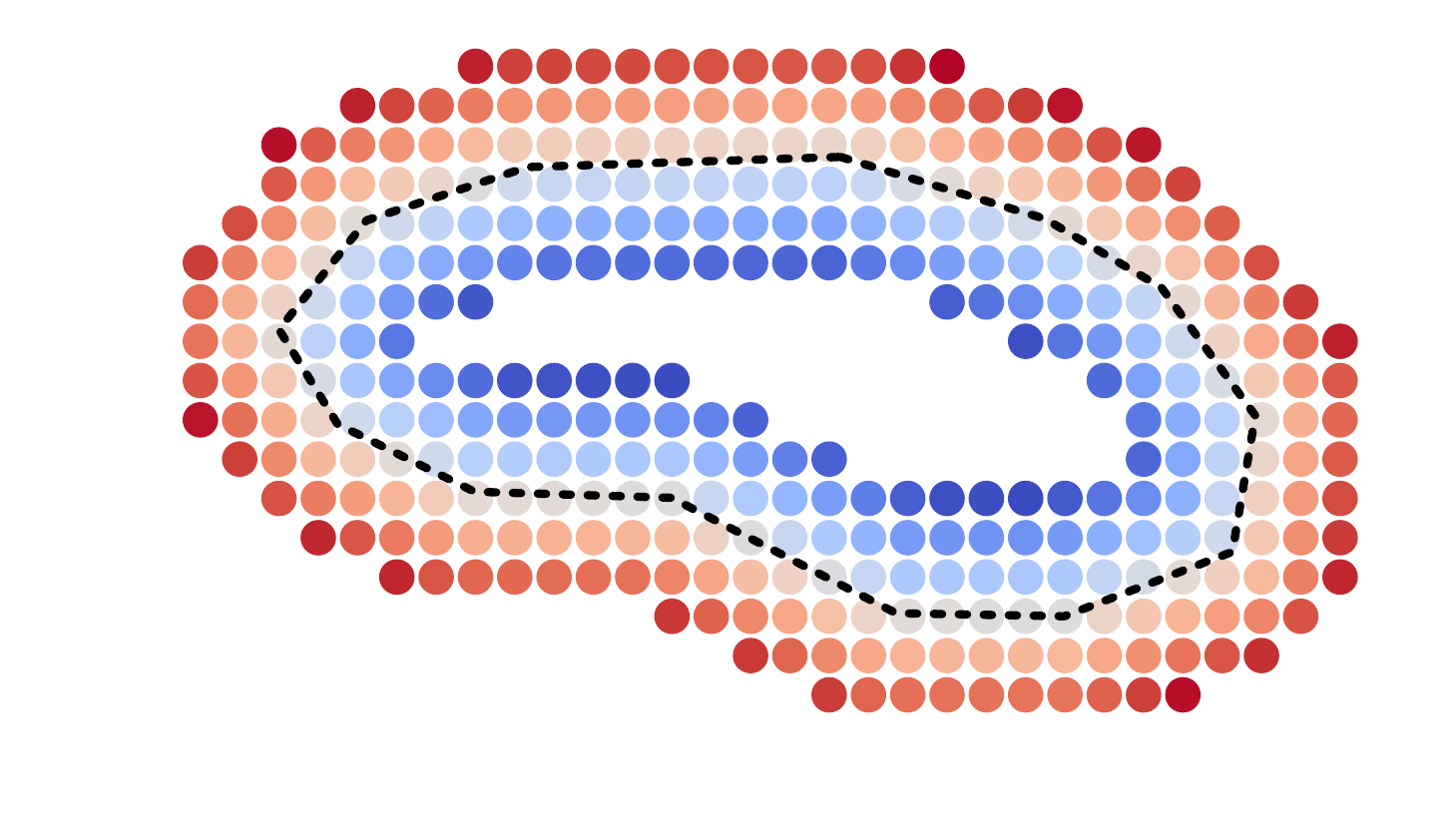}
    \caption{Signed distances to the surface.}
    \label{fig:potato_sdf}
  \end{subfigure}
  \begin{subfigure}[b]{0.22\textwidth}
    \centering
    \includegraphics[width=0.95\textwidth]{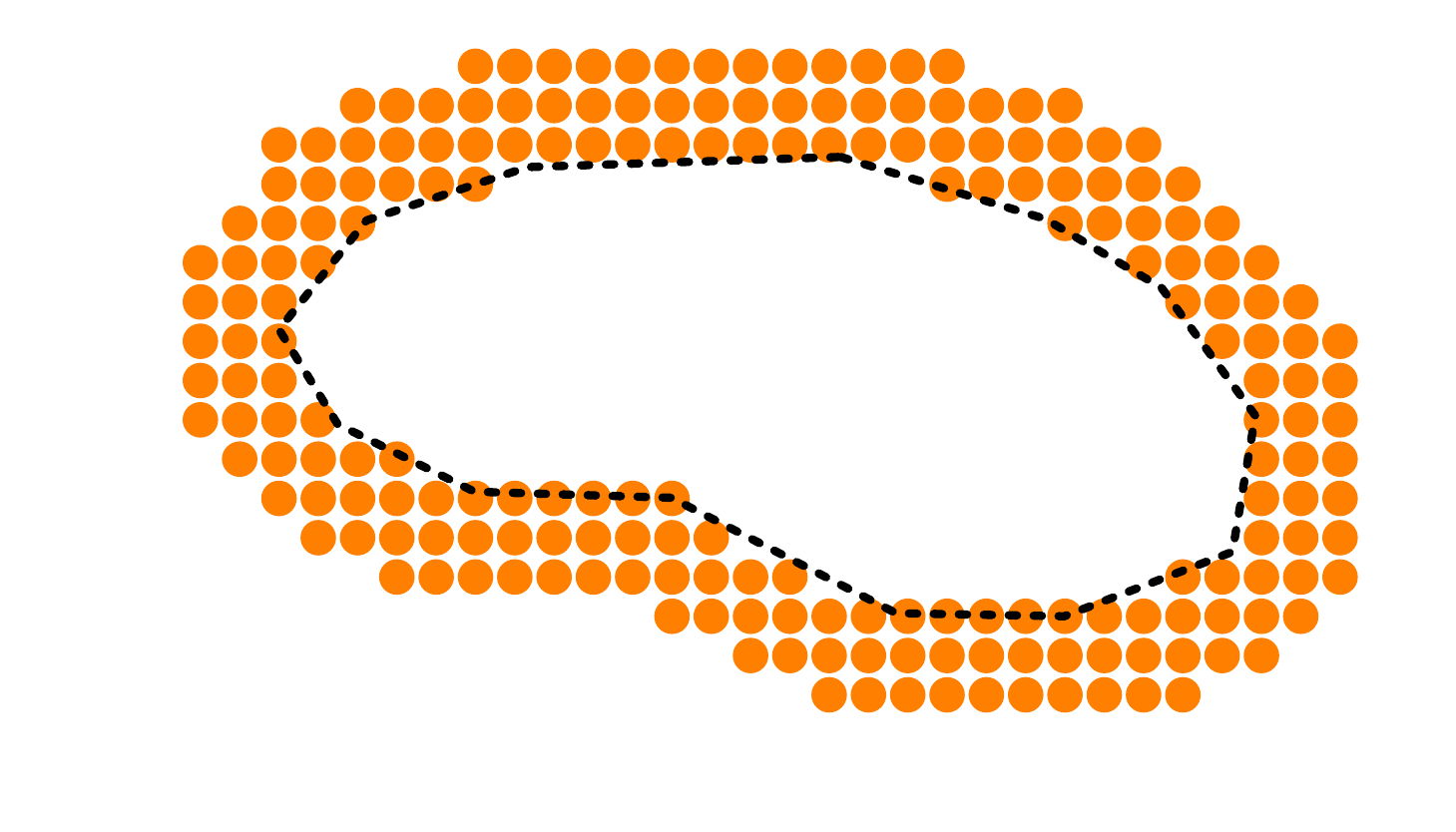}
    \caption{Boundary particles.}
    \label{fig:initial_point_grid_boundary}
  \end{subfigure}
  \hfill
  \begin{subfigure}[b]{0.22\textwidth}
    \centering
    \includegraphics[width=0.95\textwidth]{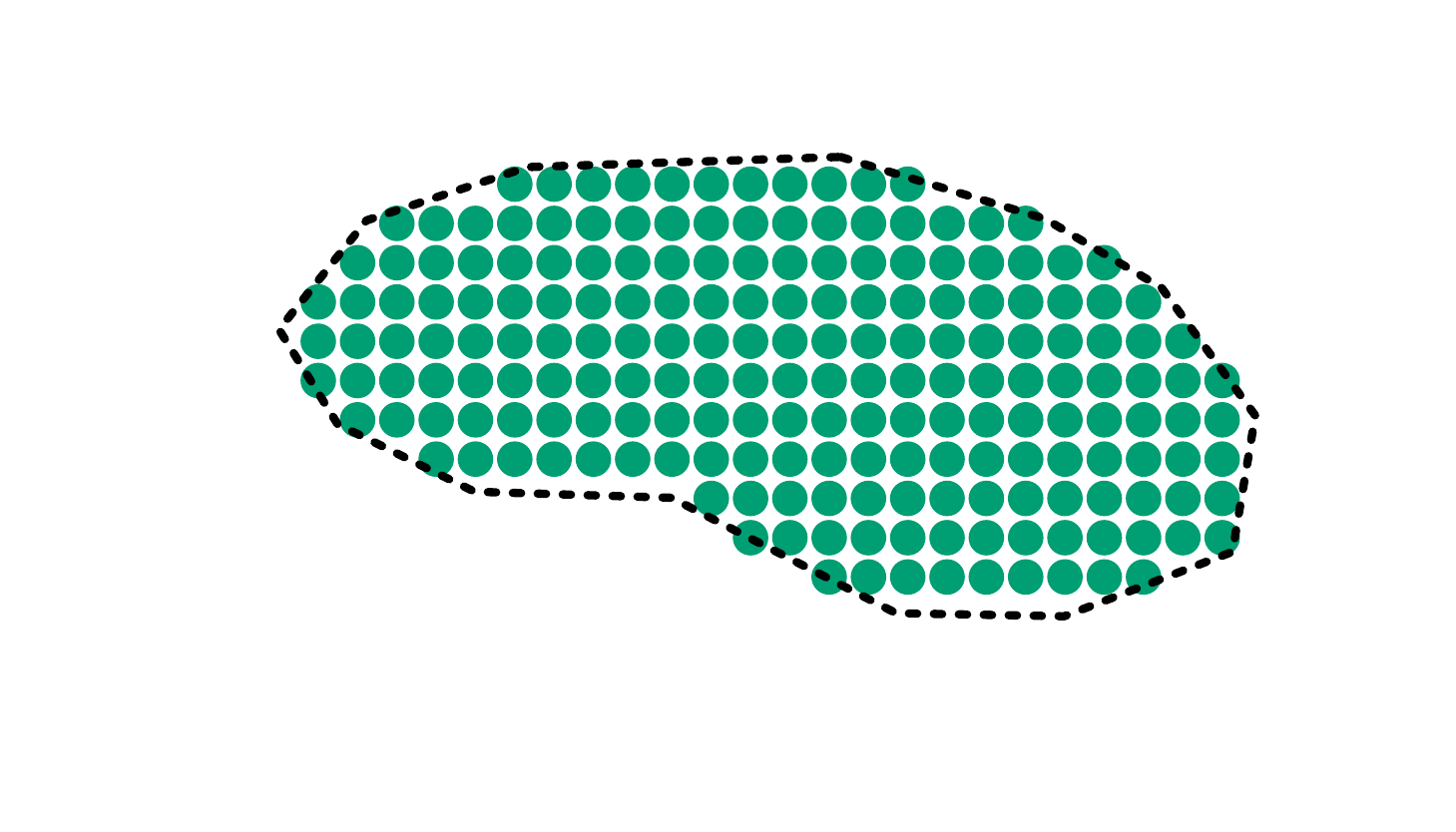}
    \caption{Interior particles.}
    \label{fig:initial_point_grid}
  \end{subfigure}
  \caption{(a)~visualization of the polygonal traversal of an arbitrary geometry surface.
  (b)~point grid surrounding the geometry's surface, where the points store the signed distances to the surface, visualized as a color map.
  (c)~and (d)~show the initial configuration of the boundary and interior particles, respectively.
  }
  \label{fig:exemplary_geometry}
\end{figure}

\begin{figure}[!h]
  \centering
  \begin{subfigure}[b]{0.45\textwidth}
    \centering
    \includegraphics[width=0.55\textwidth]{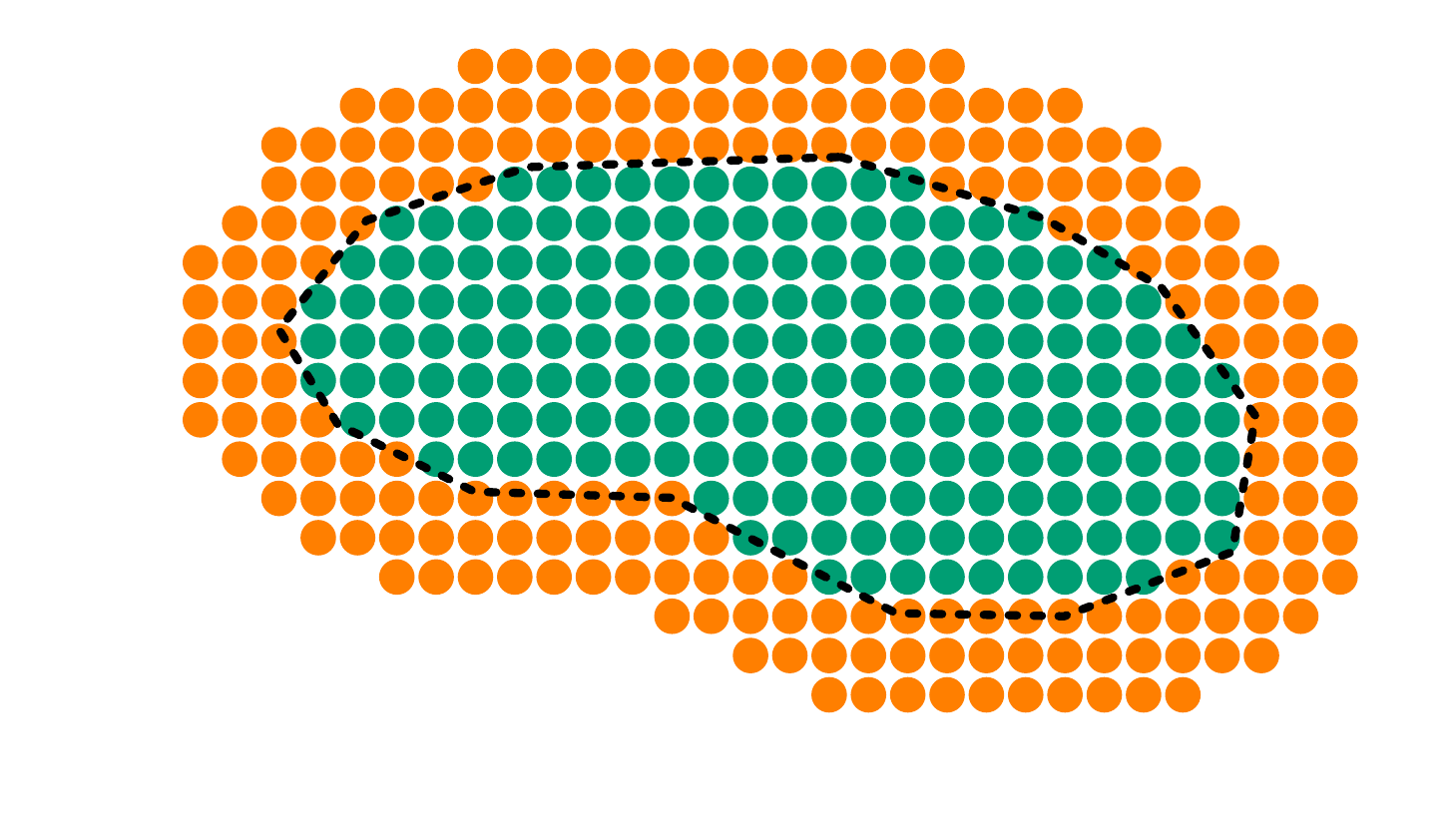}
    \caption{Initial configuration of interior (green) and boundary (orange) particles.}
    \label{fig:initial_point_grid_both}
  \end{subfigure}
  \hfill
  \begin{subfigure}[b]{0.45\textwidth}
    \centering
    \includegraphics[width=0.55\textwidth]{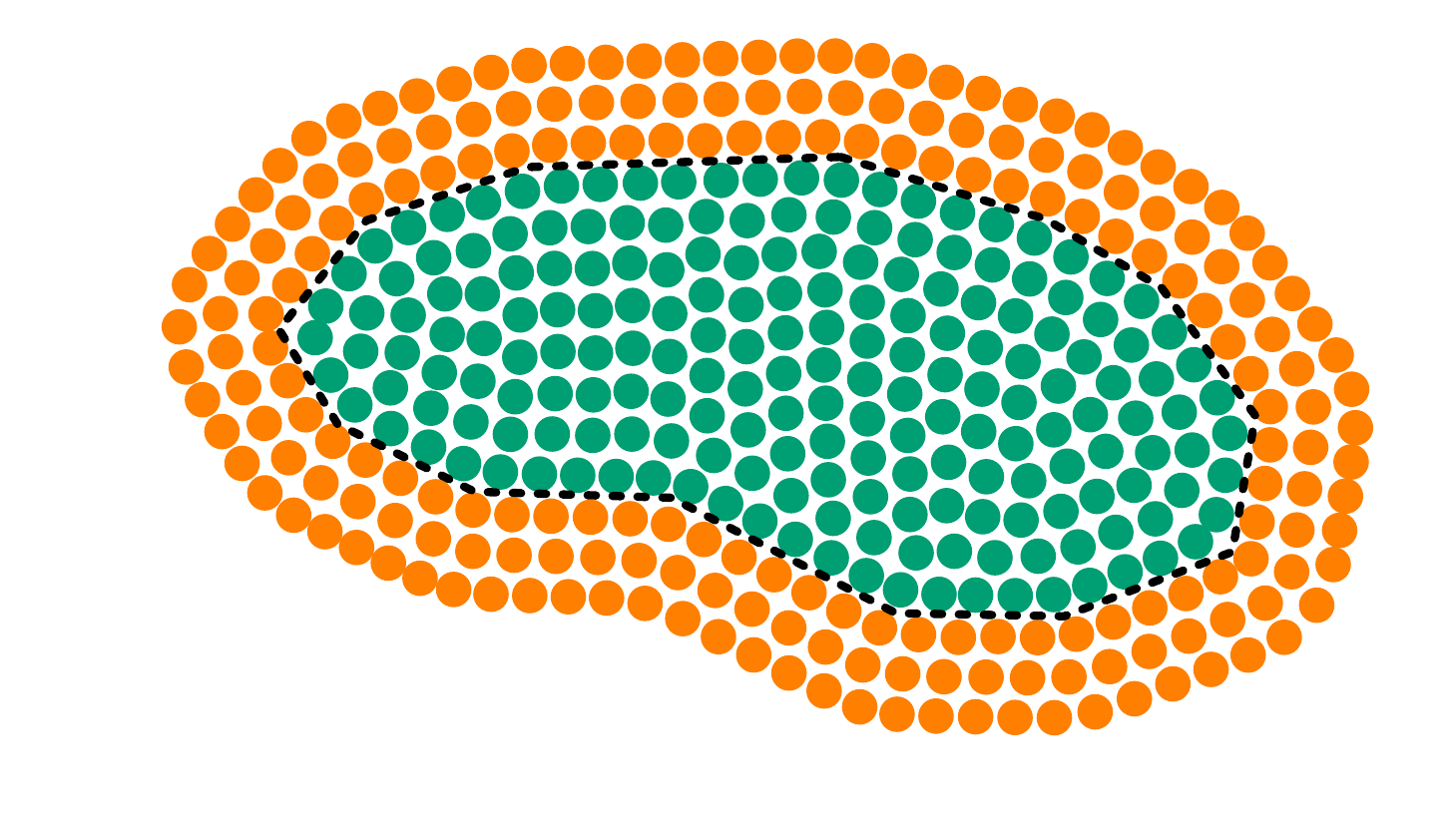}
    \caption{Final configuration of interior (green) and boundary (orange) particles after packing.}
    \label{fig:final_particle_distribution_boundary}
  \end{subfigure}
  \caption{(a)~shows the grid structure of the initial configuration, which does not accurately represent the geometry surface.
           (b)~displays the final configuration of the particles, which are isotropically distributed while accurately preserving the geometry surface.}
  \label{fig:final_distribution_arbitrary_geometry}
\end{figure}

The preprocessing pipeline consists of the following steps:
\begin{itemize}
  \item[(1)] Load geometry (Fig.~\ref{fig:potato_geometry}).
  \item[(2)] Compute signed distance field (Fig.~\ref{fig:potato_sdf}).
  \item[(3)] Generate boundary particles (Fig.~\ref{fig:initial_point_grid_boundary}).
  \item[(4)] Create initial sampling of interior particles with inside-outside segmentation (Fig.~\ref{fig:initial_point_grid}).
  \item[(5)] Pack interior and boundary particles (Fig.~\ref{fig:final_particle_distribution_boundary}).
\end{itemize}
In the first step, the geometry is loaded from a numerical representation of the object's surface.
The input data can either be a 3D triangulated surface mesh (e.g., in STL\footnote{A file format that represents a raw, unstructured triangulated surface.}) or a 2D polygonal traversal of the geometry (see Fig.~\ref{fig:potato_geometry} for an example).
Thus, the following methods apply to geometries, where the surface is discretely represented by oriented edges (in 2D) or triangles (in 3D), collectively referred to as faces, unless noted otherwise.
The second step generates the signed distance field (SDF), which is necessary for surface detection in the final packing step.
As illustrated in Fig.~\ref{fig:potato_sdf}, the SDF is computed only within a narrow band around the geometry's surface by using a face-based neighborhood search (NHS).
This NHS is used exclusively for the efficient generation of the SDF.
The points storing the SDF are fixed in space and should not be confused with the particles to be packed; they solely serve to represent the level-set information.
In the third step, the initial configuration of the boundary particles is generated (orange particles in Fig.~\ref{fig:initial_point_grid_boundary}).
To create the boundary particles, the positions of the SDF points located outside the geometry and within a predefined boundary thickness are copied.
In the fourth step, the initial configuration of the interior particles (green particles in Fig.~\ref{fig:initial_point_grid}) is generated.
Here, the hierarchical winding number approach \cite{Jacobson2013} is used to efficiently perform an inside-outside segmentation of particles at their initial, lattice-based positions,
even when far away from the surface.
After steps (1) to (4), the initial configuration of both interior and boundary particles has been obtained, as illustrated in Fig.~\ref{fig:initial_point_grid_both}.
However, the interface of the geometry surface is not yet well resolved.
Therefore, in the final step (5), a packing algorithm is applied, utilizing the SDF to simultaneously optimize the positions of both interior and boundary particles.
This yields an isotropic particle distribution while accurately preserving the geometry surface, as illustrated in Fig.~\ref{fig:final_particle_distribution_boundary}.
In the subsequent sections, each step of the pipeline is explained in detail.

% !TeX root = ../main.tex

\section{Smoothed Particle Hydrodynamics}
\label{sec:sph}
Before describing the individual steps of the preprocessing, smoothed particle hydrodynamics (SPH) is introduced as the foundation of the particle-based approach employed in this work.
The following section provides a brief overview of the SPH methodology and the time integration scheme to establish the necessary background.

\subsection{Spatial discretization}
\label{sec:spatial_discretization}
SPH is a versatile, mesh-free simulation method based on a Lagrangian formulation widely employed in fields ranging from astrophysics~\cite{Monaghan1977}
to fluid dynamics~\cite{Monaghan1994} and solid mechanics~\cite{Gray2001}.
Its particle-based approach naturally accommodates large deformations, free surfaces, and complex geometries.
The SPH method continuously approximates any spatial function $f$ by a convolution integral with a kernel function $W$,
\begin{align}
    f(\mathbf{r}) \approx (f * W)(\mathbf{r}) &=  \int_{\Omega} f(\mathbf{r}')
    W(\mathbf{r}-\mathbf{r}', h) \mathrm{d}\mathbf{r}' \\
    &= \int_{\Omega} \frac{f(\mathbf{r}')}{\rho(\mathbf{r}')}
    W(\mathbf{r}-\mathbf{r}', h)\underbrace{\rho(\mathbf{r}')
        \mathrm{d}\mathbf{r}'}_{\mathrm{d}m'},
    \label{eq:convolution_integral_sph}
\end{align}
where $\mathbf{r}$ and $\mathbf{r}'$ are coordinates in the domain $\Omega$ and $\rho$ denotes the local density.
In this work, the quintic spline kernel \cite{Schoenberg1946} is used, which is given by
\begin{equation}
    W(\mathbf{r}, h) = \frac{1}{h^d} w(\Vert\mathbf{r}\Vert /h),
\end{equation}
with
\begin{equation}
w(q) = \sigma \begin{cases}
    (3 - q)^5 - 6(2 - q)^5 + 15(1 - q)^5    & \text{if } 0 \leq q < 1, \\
    (3 - q)^5 - 6(2 - q)^5                  & \text{if } 1 \leq q < 2, \\
    (3 - q)^5                               & \text{if } 2 \leq q < 3, \\
    0                                       & \text{if } q \geq 3.
\end{cases}
\end{equation}
Here, $d$ is the number of dimensions, $h$ denotes the kernel's smoothing length, and $\sigma$ is a normalization constant given by
$\sigma = \left[\frac{1}{120}, \frac{7}{478 \pi}, \frac{1}{120\pi}\right]$ in $[1, 2, 3]$ dimensions.
We use $\Vert \cdot \Vert$ to denote the Euclidean norm throughout this paper.
For SPH discretizations, the continuous integral in Eq.~\eqref{eq:convolution_integral_sph} is replaced by a sum over discrete interpolation points, the particles,
\begin{equation}
    (f * W)(\mathbf{r}_i) \approx \sum_{j \in \mathcal{N}} f_j
    \frac{m_j}{\rho(\mathbf{r}_j)}  W(\mathbf{r}_i - \mathbf{r}_j, h),
    \label{eq:discretized_sph_formulation}
\end{equation}
where $\mathbf{r}_i$ is the position of the $i$-th particle, $\mathbf{r}_j$ and $m_j$ are the position and mass of the $j$-th neighboring particle
and $\mathcal{N}$ is the set of neighboring particles within the kernel's compact support.
In the following, we abbreviate $W_{ij} \coloneqq W(\mathbf{r}_i-\mathbf{r}_j,  h)$.
By setting $f_i = \rho_i$, the density $\rho_i = \rho(\mathbf{r}_i)$ at particle position $i$ is computed by the SPH density estimate as
\begin{equation}
    \rho_i = \sum_{j \in \mathcal{N}} m_j W_{ij}.
    \label{eq:summation_density}
\end{equation}
Gradients can be calculated by differentiating Eq.~\eqref{eq:discretized_sph_formulation} with respect to the particle position,
\begin{equation}
    \nabla f \approx \sum_{j \in \mathcal{N}} f_j \frac{m_j}{\rho_j} \nabla W_{ij},
    \label{eq:gradient_sph}
\end{equation}
where $\nabla  W_{ij}$ is the gradient of the kernel function, for which an analytical expression is known \cite{Schoenberg1946}.

In general, the formulation in Eq.~\eqref{eq:gradient_sph} leads to poor gradient estimates.
To improve these estimates, alternative formulations have been derived by Price et al.~\cite{Price2012} and Monaghan et al.~\cite{Monaghan2005}.
However, this poor gradient estimation can still be exploited for a physical particle packing as described in the next section.

\subsection{Governing equations}
To rearrange the particles into a stationary distribution, we employ the TVF of Adami et al.~\cite{Adami2013}.
In the following, we strictly follow their notation and formulation, and give a brief summary for completeness.
The advection velocity $\mathbf{\tilde{v}}_i$ of particle $i$ is used to evolve the position of the particle~$\mathbf{r}_i$ from one time step to the next by
\begin{equation}
\frac{\mathrm{d} \mathbf{r}_i}{\mathrm{d}t} = \mathbf{\tilde{v}}_i .
\end{equation}
It is obtained in each new time step by computing
\begin{equation}
\mathbf{\tilde{v}}_i (t + \Delta t) = \mathbf{v}_i (t) + \Delta t \left(\frac{\tilde{\mathrm{d}} \mathbf{v}_i}{\mathrm{d}t} + \mathbf{a}_p \right),
\end{equation}
where $\Delta t$ is the time step size and $\mathbf{v}_i$ is the momentum velocity.
$\tilde{\mathrm{d}} \mathbf{v}_i / \mathrm{d}t$ is zero for the packing algorithm and only shown here for completeness.
Likewise, the momentum velocity~$\mathbf{v}_i$ is set to zero at the beginning of each time step, as we want to achieve a fully stationary state.
The discretized form of the remaining term is
\begin{equation}
    \mathbf{a}_p = -\frac{1}{\rho_i} \nabla p_{\text{b}} \approx  -\frac{2p_{\text{b}}}{\rho_i} \sum_{j \in \mathcal{N}} \frac{m_j}{\rho_j}  \nabla W_{ij},
    \label{eq:pressure_force}
\end{equation}
where $p_{\text{b}}$ is a constant background pressure.
Note that although in the continuous case $\nabla p_{\text{b}} =  \bm{0}$, the discretization is not 0th-order consistent for non-uniform particle distribution,
which means that there is a non-vanishing contribution when particles are disordered.
This implies that $p_{\text{b}}$ drives a particle's trajectory, thereby promoting uniform pressure distributions.
For a uniform particle distribution with a constant background pressure and uniform density and mass,
the acceleration vanishes due to the properties of the kernel function
\begin{equation}
    \sum_j \nabla_i W_{ij} = \bm{0},
    \label{eq:zero_gradient}
\end{equation}
which only holds under the assumption of a full kernel support (see Sec.~\ref{sec:dynamic_boundary_packing}).
The background pressure can be chosen arbitrarily large when the time step criterion is adjusted.
This implies that $p_{\text{b}}$ scales the simulation time, but the final state remains unchanged after the same number of time steps.
In other words, increasing $p_b$ reduces the time step size, but simultaneously accelerates the relaxation process,
such that after the same number of time steps the system reaches the same physical state, even though the simulation time differs.
In the following, we choose $p_{\text{b}} = 1$ for simplicity.

\subsection{Time integration}
\label{sec:time_integration}
We use an explicit Runge-Kutta method (\texttt{RDPK3SpFSAL35}~\cite{Ranocha2021} if not mentioned otherwise) provided by the OrdinaryDiffEq.jl interface~\cite{rackauckas2017differentialequations} with automatic error-based adaptive time stepping.
The advection velocity $\mathbf{\tilde{v}}$ and the positions of the particles are updated across the stages.
After each full time step, the momentum and advection velocities are reset to zero.

% !TeX root = ../main.tex
\section{Point Cloud with Signed Distance Field}
\label{sec:signed_distance_point_cloud}
This section focuses on the second step of the overall pipeline, as described in Sec.~\ref{sec:overview}. We first provide a brief overview of the level-set function,
followed by a description of the SDF generation.

The SDF is defined by the level-set function $\phi(\mathbf{x})$ and is used to obtain the information about a surface $ \partial \Omega$
enclosing an object $\Omega$ in 2D or 3D space.
\begin{equation}
    \{\mathbf{X} \in \mathbb{R}^d \mid \phi(\mathbf{X}) = 0\} = \partial \Omega
\end{equation}
represents the object surface.
The sign of $\phi$ determines if the point is inside or outside of the object:
\begin{equation}
    \phi(\mathbf{x}) = \begin{cases}
        -\text{dist}(\mathbf{x}, \partial \Omega) & \text{if } \mathbf{x} \in \Omega, \\
        \text{dist}(\mathbf{x}, \partial \Omega) & \text{if } \mathbf{x} \notin \Omega,
    \end{cases}
    \label{eq:level_set_definition}
\end{equation}
where the function $\text{dist}(\mathbf{x}, \partial \Omega)$ gives the distance from a point $\mathbf{x}$ to the surface $\partial \Omega$.
Also, the normal direction  $\mathbf{n}$ to the surface $\partial \Omega$ at an arbitrary point $\mathbf{x}$ to the surface $ \partial \Omega$ is given by
\begin{equation}
    \mathbf{n} = \frac{\nabla \phi}{\Vert\nabla \phi\Vert}.
    \label{eq:level_set_normal}
\end{equation}

In practice, for a discrete case, $\phi_i = \phi(\mathbf{x}_i)$ is evaluated at an arbitrary position $\mathbf{x}_i$
by traversing each face of the triangle mesh to find the minimum distance to the surface.
This means that the SDF is discretized by the positions at which $\phi$ is evaluated.

In this work, we do not use the SDF to perform an inside-outside segmentation
(which would require resolving a domain enclosing the entire geometry; see Sec.~\ref{sec:segmentation}).
Instead, the SDF is employed solely for surface detection (see Sec.~\ref{sec:bounding_method}).
Accordingly, it is sufficient to evaluate $\phi$ only in the vicinity of the surface.
This localized approach is both memory-efficient and computationally advantageous,
since it enables the use of a face-based NHS for generating the SDF.
Before introducing this NHS, we first explain why it is sufficient to compute the SDF only in the vicinity of the surface.

\subsection{Signed distance field interpolation}
\label{sec:signed_distance_field_interpolation}
Our goal is to interpolate the distance and normal direction to the surface for particles located near the geometry.
To this end, we generate a point cloud close to the surface.
The points of the cloud should not be confused with the particles to be packed.
These points are fixed in space and store the information of~$\phi$ and $\mathbf{n}$ for their respective positions.
In principle, the positions of these points can be chosen arbitrarily; in this work, however, we sample them on a regular grid.
\begin{figure}[h!]
    \centering
    \includegraphics[width=0.4\textwidth]{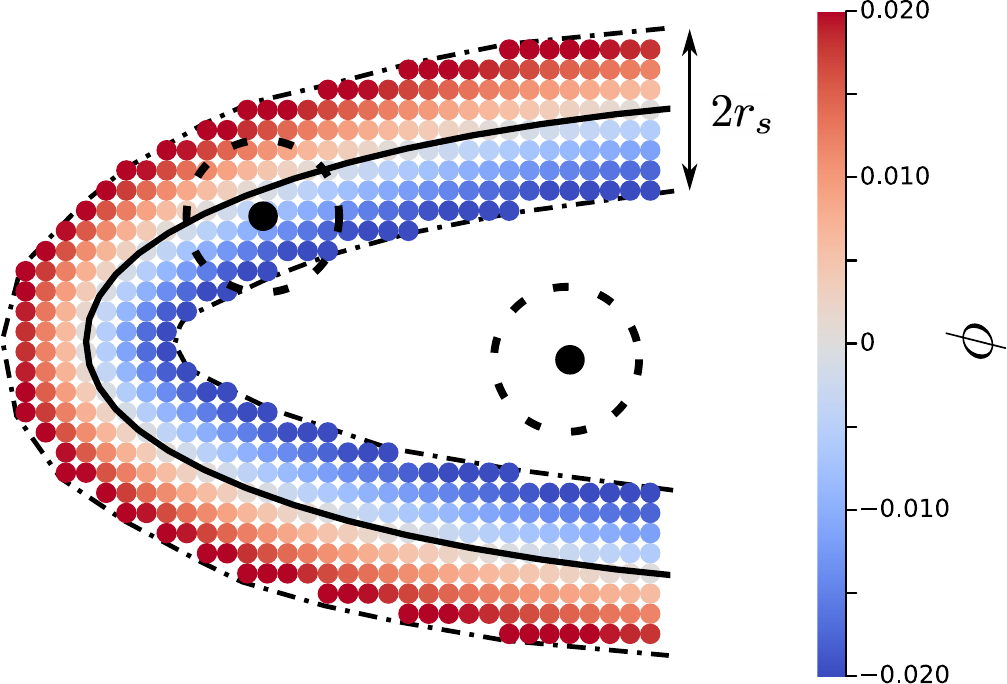}
    \caption{Interpolation of the SDF onto the particles (black dots), where the dashed circles represent the compact support of the kernel.
            The SDF is shown with colored points. The solid black line represents the geometry.
            The dashed-dotted lines represent the thickness of the convex hull of the point cloud.}%
    \label{fig:sdf_interpolation}
\end{figure}%
We interpolate the information of the SDF to the particle positions using Shepard interpolation~\cite{Shepard1968}:
\begin{equation}
    \phi_i =  \frac{\sum_{j \in \mathcal{N}_p} \phi_j W_{ij}}{\sum_{j \in \mathcal{N}_p} W_{ij}},
    \label{eq:interpolate phi}
\end{equation}
\begin{equation}
    \mathbf{n}_i = \frac{\sum_{j \in \mathcal{N}_p} \mathbf{n}_j W_{ij}}{\sum_{j \in \mathcal{N}_p} W_{ij}},
    \label{eq:interpolate n}
\end{equation}
where $\phi_i$ is the signed distance at the position of the $i$th particle, $\mathbf{n}_i$ is the normal direction to the particle's closest surface region, and $\mathcal{N}_p$ is the set of neighboring points in the point cloud representing the SDF.
Note that the interpolation, being based on the smoothing kernel with a compact support $r_s$, only requires $\phi$ values within a search radius of $r_s$.

Fig.~\ref{fig:sdf_interpolation} illustrates the interpolation, where two particles are depicted (black dots) with the corresponding~$r_s$ (dashed line).
The SDF is visualized by colored points representing the distance values, with the solid black line indicating the geometry surface.
Particles located near the surface interact with the field, while those deeper inside the geometry do not.
Consequently, it is sufficient to construct a point cloud along the surface with a thickness of~$2r_s$ (outlined by the dashed-dotted lines),
providing the necessary information for interpolation.
Note that the normal vectors are directly derived from the surface (see Eq.~\ref{eq:interpolate n}) and not by numerical differentiation of the gradient field,
as is usually the case with the level-set method (see Eq.~\eqref{eq:level_set_normal}).

\subsection{Face-based neighborhood search}
\label{sec:face_based_nhs}
The following describes the initialization of the face-based NHS, which is only required for generating the SDF and is not used afterward.
The procedure for the face-based NHS is also presented in Alg.~\ref{alg:face-nhs}.
A performance analysis of the introduced NHS is presented in Sec.~\ref{sec:performance}.
We use a NHS method which is commonly used in particle-based methods, as detailed in \cite{Ihmsen2010} and \cite{Chalela2021}.
However, the difference is that, unlike with a particle NHS, cells are not populated with particle indices, but with face indices.
The challenge here is that a face is not a single point for which the cell coordinates can be easily calculated,
but a triangle spanned by three vertices (in 3D) or an edge defined by two vertices (in 2D).

We divide the domain enclosing the entire geometry into a regular grid, where the position of the cells of the grid is given in 2D or 3D by
\begin{equation}
    \left( \left\lfloor \frac{x_1}{r_s} \right\rfloor, \left\lfloor \frac{x_2}{r_s} \right\rfloor \right) \quad \text{or} \quad
    \left( \left\lfloor \frac{x_1}{r_s} \right\rfloor, \left\lfloor \frac{x_2}{r_s} \right\rfloor, \left\lfloor \frac{x_3}{r_s} \right\rfloor \right),
    \label{eq:cell_coords}
\end{equation}
where $x$, $y$, $z$ are the space coordinates and $r_s$ is the search radius.
A potentially infinite domain is represented by cells of size $r_s$ in a hash table, where only filled cells are stored.
The search radius $r_s$ corresponds to the maximum required signed distance, e.g., the chosen boundary thickness, and should not be smaller than the compact support of the kernel.
In the following, for simplicity, the initialization of the NHS is described using a 2D example with edges instead of triangles.
\begin{figure}[!h]
    \centering
    \begin{subfigure}[b]{0.3\textwidth}
      \centering
      \includegraphics[width=0.75\textwidth]{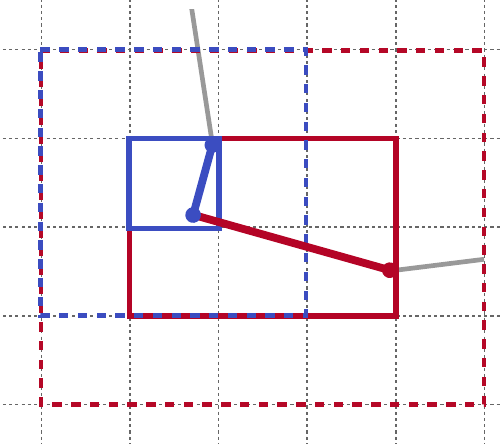}
      \caption{$r_s = 6 \Delta x$ }
      \label{fig:nhs_example_a}
    \end{subfigure}
    \hfill
    \begin{subfigure}[b]{0.3\textwidth}
      \centering
      \includegraphics[width=0.75\textwidth]{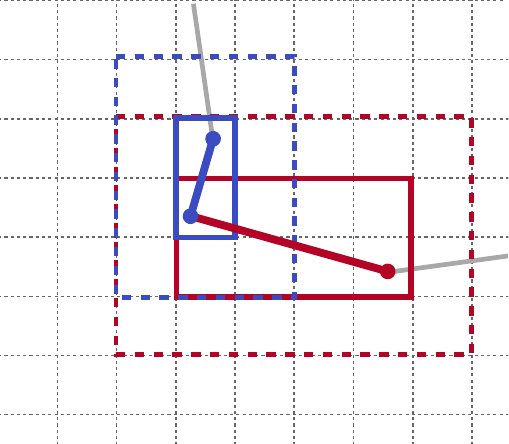}
      \caption{$r_s = 4 \Delta x$ }
      \label{fig:nhs_example_b}
    \end{subfigure}
    \hfill
    \begin{subfigure}[b]{0.3\textwidth}
      \centering
      \includegraphics[width=0.75\textwidth]{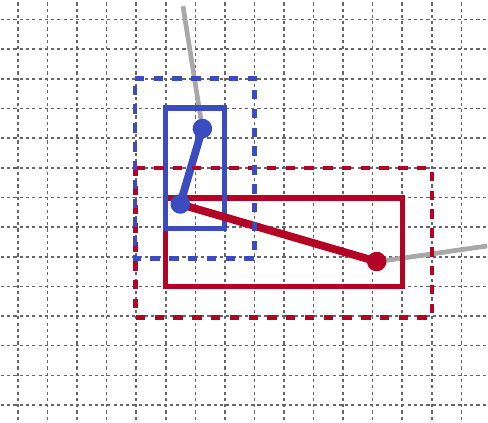}
      \caption{$r_s = 2\Delta x$}
      \label{fig:nhs_example_c}
    \end{subfigure}
    \caption{            Two exemplary edges (blue and red) in cell grids, each with a different search radius $r_s$, where $\Delta x$ is the particle spacing.
    The grey edges are part of the polygon traversal but not considered in this example.
    The bounding boxes of the edges are framed with the corresponding color.
    The dashed boxes represent the direct neighbor cells of the bounding boxes.}%
    \label{fig:nhs_example}
  \end{figure}%
In Fig.~\ref{fig:nhs_example}, two exemplary edges, blue and red, are shown, which are assigned to three different cell sizes.
First, as described in Alg.~\ref{alg:face-nhs} line~\ref{alg:line:bbox}, the bounding box is determined for each edge.
In the figure, it is framed with the corresponding color.
To do this, the cell coordinates for each vertex are computed using Eq.~\eqref{eq:cell_coords} to determine the minimum and maximum cell coordinates.
This allows the Cartesian indices of the cells in the bounding box to be determined using the minimum and maximum cell.
All the cells within the bounding box are filled with the associated edge index (Alg.~\ref{alg:face-nhs}, line~\ref{alg:line:add_face}).
Note that the direct neighboring cells must also be considered to handle edge cases.
For example, in Fig.~\ref{fig:nhs_example_a}, the vertex of the red edge is very close to the bounding box boundary, which is why the neighboring cell must also be considered.
Thus, the edge indices are also propagated to the direct neighboring cells (Alg.~\ref{alg:face-nhs}, line~\ref{alg:line:neighboring_ids}).
Also note that the bounding boxes overlap.
Thus, it is important to ensure that the edge indices in the cells are unique to avoid multiple entries of the same edge.
The performance of the NHS depends on the search radius $r_s$, which depends on the thickness of the SDF.
As illustrated in Fig.~\ref{fig:nhs_example}, the larger $r_s$ is chosen, the more face indices might be stored per cell, which in turn means more neighboring faces in the cell neighbors.
Conversely, a very small $r_s$ requires a larger number of cells, which in turn increases the memory requirement.
This is more detailed in Sec.~\ref{sec:performance_sdf}.
\begin{algorithm}[!t]
    \caption{Face-based neighborhood search}
    \begin{algorithmic}[1]
        \Function{initialize\_NHS}{geometry, search\_radius}
        \State cell\_hashtable $\gets$ \Call{create\_hashtable}{{}}
        \For{face in geometry}
            \State cells $\gets$ \Call{bounding\_box}{face, search\_radius} \label{alg:line:bbox}
            \Comment{returns cell IDs in the bounding box of the face}
            \For{cell in cells}
                 \State \Call{append}{cell\_hashtable[cell], face} \label{alg:line:add_face}
                 \Comment{add face ID to cell}
            \EndFor
        \EndFor

        \State cell\_grid $\gets$ \Call{bounding\_box}{geometry, search\_radius}
        \Comment{returns cell IDs in the bounding box of the geometry}
        \For{cell in cell\_grid}
        \State face\_ids $\gets$ \Call{faces\_in\_neighboring\_cells}{cell}\label{alg:line:neighboring_ids}

            \State cell\_hashtable[cell] $\gets$ \Call{unique}{face\_ids}
            \Comment{make sure all face IDs are unique}

        \EndFor
        \State \textbf{return} cell\_hashtable
        \EndFunction
    \end{algorithmic}
    \label{alg:face-nhs}
\end{algorithm}

\subsection{Generating the signed distance field}
\label{sec:determine_signed_distances}
\begin{figure}[!t]
    \center
    \includegraphics[width=0.7\textwidth]{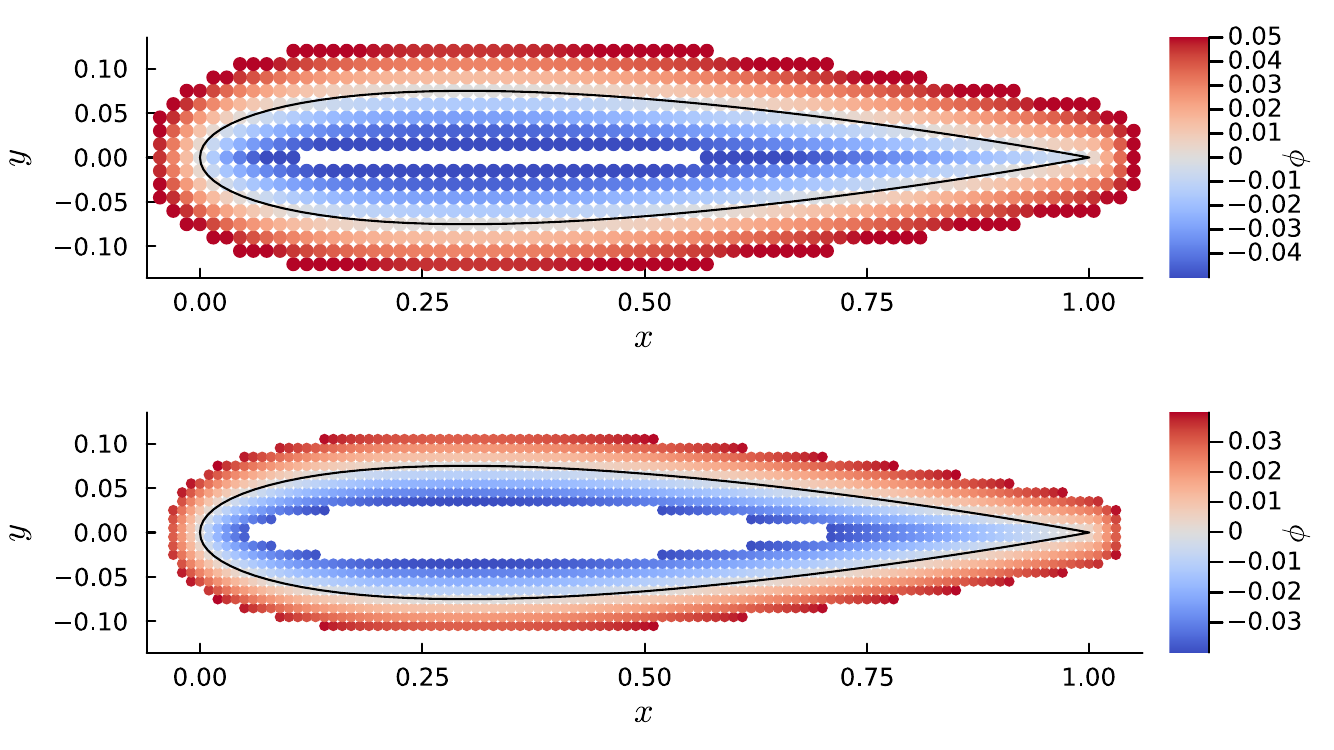}
    \caption{Signed distance point cloud for different resolutions with a maximum signed distance of~$4\Delta x$ requiring a search radius of~$r_s = 4\Delta x$.
    The geometry surface is represented by the black solid line. Top: $\Delta x = 0.015$, bottom: $\Delta x = 0.01$.}%
    \label{fig:airfoil_sdf}
\end{figure}

\begin{algorithm}[!t]
    \caption{Signed distance field}
    \begin{algorithmic}[1]
        \Function{calculate\_signed\_distances}{geometry, points}
        \State NHS $\gets$ \Call{initialize\_NHS}{geometry, search\_radius}
        \Comment{see Alg.~\ref{alg:face-nhs}}

        \State \Call{remove\_points\_in\_empty\_cells}{points, NHS} \label{alg:line:remove_points}
        \For{point in points}
            \State {$\phi_\text{min} \gets$  $\infty$}%revA
            \For{face in neighbor faces of point}
                \State $\phi \gets$ \Call{calculate\_signed\_distance\_to\_face}{point, face}
                \If{{$\phi < \phi_\text{min}$}}%revA
                    \State {$\phi_\text{min} \gets \phi$}%revA
                    \State {$n \gets$ \Call{normal\_vector\_to\_surface}{face}}%revA
                    \State {distances[point] $\gets$ \Call{store}{$\phi$, $n$}}%revA
                \EndIf
            \EndFor
        \EndFor
        \For{point in points}
            \If{$\left| \text{distances[point]} \right|$  > search\_radius}
                \State  \Call{remove\_point}{distances, point}
            \EndIf
        \EndFor
        \EndFunction
    \end{algorithmic}
    \label{alg:sdf}
\end{algorithm}%
To generate the point cloud, we create an initial point grid in the bounding box of the geometry.
Optimizations creating this grid are possible, but not necessary,
as only points in cells filled with face indices are considered.
The signed distances are determined by Alg.~\ref{alg:sdf}.
First, all points in cells without face indices are removed (Alg.~\ref{alg:sdf}, line~\ref{alg:line:remove_points}).
Then, the remaining points are iterated over, and the signed distance is calculated for each neighboring face.
This means that the distance to each neighboring face is determined.
A detailed description of the distance to triangle calculation can be found in~\cite{Baerentzen2002}.
The minimum distance and its associated normal vector are then assigned to the respective point position.
Finally, the positions that do not lie within the specified maximum distance are deleted.
Since only a point cloud along the surface needs to be generated,
the resolution of the SDF can be chosen according to the resolution of the initial particle configuration,
significantly reducing memory and computation compared to full-volume discretizations.
To demonstrate the surface-based SDF generation process, the resulting point cloud for the NACA0015 airfoil for different resolutions is shown in Fig.~\ref{fig:airfoil_sdf}.

% !TeX root = ../main.tex
\section{Initial Particle Configuration}
\label{sec:initial_particle_configuration}
In this section, we describe steps (3) and (4) from Sec.~\ref{sec:overview}.
Boundary particles are generated using the precomputed SDF, while interior particles are created through an inside-outside segmentation process.

\subsection{Boundary particle generation}
\label{sec:boundary_particle_generation}
When generating the point cloud with the SDF along the surface (see Sec.~\ref{sec:signed_distance_point_cloud}), points are also sampled outside the surface.
For each point within the boundary hull, i.e., outside the geometry but within a given distance, the boundary thickness,
the positions of these points are used to sample boundary particles.

\subsection{Inside-outside segmentation}
\label{sec:segmentation}
The objective of Step (4) from Sec.~\ref{sec:overview} is to classify query points as either inside or outside the geometry.
Query points identified as inside are assigned to the interior particle set, while those outside are discarded.
While the initial sampling of the query points can, in principle, be performed arbitrarily,
we sample the bounding box of the geometry with a uniform point grid, yielding $n$ query points.
To determine whether a query point lies inside the geometry, we compute its winding number (as described in Sec.~\ref{sec:winding_number}).
In a naive implementation, this involves evaluating the winding number for every face $m$ of the geometry,
resulting in a computational complexity of $\mathcal{O}(n\cdot m)$.
To address this, we employ the hierarchical winding number algorithm proposed by Jacobson et al.~\cite{Jacobson2013}, which reduces the computational complexity to $\mathcal{O}(n\sqrt{m})$.
This algorithm leverages a hierarchical structure to efficiently partition the geometry and evaluate the winding number.

In the following, we provide a brief introduction to the concept of the winding number,
starting with its definition in the 2D case~\cite{Hormann2001} and subsequently extending it to 3D~\cite{Jacobson2013}.
This is followed by a description of the hierarchical winding number method.
For a more detailed explanation, we refer to the section ``Hierarchical Winding'' in the TrixiParticles.jl documentation~\cite{trixiparticles_documentation}.

\subsubsection{Winding number}
\label{sec:winding_number}
The winding number $w(\mathbf{p}, \mathcal{C})$ of a point $\mathbf{p}$
with respect to a closed Lipschitz curve $\mathcal{C}$ in $\mathbb{R}^2$ is the number of revolutions made around $\mathbf{p}$ when traversing $\mathcal{C}$ \cite{Hormann2001}.
That is, $w$ represents the total number of turns around $\mathcal{C}$.
If $\mathcal{C}$ is traversed in a clockwise direction around the point, $w$ is negative, otherwise it is positive.
The point $\mathbf{p}$ lies outside the curve $\mathcal{C}$ if $w=0$.
For a continuous curve, the winding number is defined by
\begin{equation}
  w(\mathbf{p}) = \frac{1}{2\pi} \oint_{\mathcal{C}} \mathrm{d}\Theta.
\end{equation}
This can be generalized to three dimensions \cite{Jacobson2013} by
\begin{equation}
  w(\mathbf{p}) = \frac{1}{4\pi} \Phi(\mathbf{p}),
\end{equation}
where $\Phi(\mathbf{p}) = \iint_{\partial\Omega} \sin(\phi) \,\mathrm{d}\Theta \mathrm{d}\phi$ is the signed surface area of the projection of the geometry surface $\partial\Omega$ onto the unit sphere centered at~$\mathbf{p}$.
In three dimensions, $w$ counts the signed total number of times the surface~$\partial\Omega$ wraps around the
point~$\mathbf{p}$.

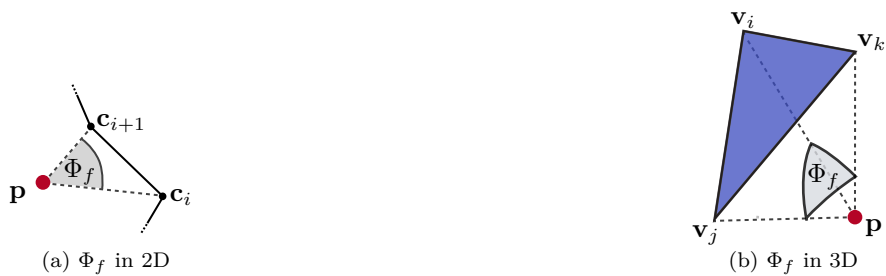
\begin{figure}[!h]
  \centering
  \begin{subfigure}[b]{0.45\textwidth}
    \centering
    \begingroup%
    \makeatletter%
    \providecommand\color[2][]{%
      \errmessage{(Inkscape) Color is used for the text in Inkscape, but the package 'color.sty' is not loaded}%
      \renewcommand\color[2][]{}%
    }%
    \providecommand\transparent[1]{%
      \errmessage{(Inkscape) Transparency is used (non-zero) for the text in Inkscape, but the package 'transparent.sty' is not loaded}%
      \renewcommand\transparent[1]{}%
    }%
    \providecommand\rotatebox[2]{#2}%
    \newcommand*\fsize{\dimexpr\f@size pt\relax}%
    \newcommand*\lineheight[1]{\fontsize{\fsize}{#1\fsize}\selectfont}%
    \ifx\svgwidth\undefined%
      \setlength{\unitlength}{72.95798294bp}%
      \ifx\svgscale\undefined%
        \relax%
      \else%
        \setlength{\unitlength}{\unitlength * \real{\svgscale}}%
      \fi%
    \else%
      \setlength{\unitlength}{\svgwidth}%
    \fi%
    \global\let\svgwidth\undefined%
    \global\let\svgscale\undefined%
    \makeatother%
    \begin{picture}(1,0.75648263)%
      \lineheight{1}%
      \setlength\tabcolsep{0pt}%
      \put(0,0){\includegraphics[width=\unitlength,page=1]{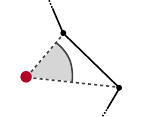}}%
      \put(-0.00171893,0.17923303){\makebox(0,0)[lt]{\lineheight{0}\smash{\begin{tabular}[t]{l}$\mathbf{p}$\end{tabular}}}}%
      \put(0.4454137,0.55190892){\makebox(0,0)[lt]{\lineheight{0}\smash{\begin{tabular}[t]{l}$\mathbf{c}_{i+1}$\end{tabular}}}}%
      % \put(0.22,0.45){\makebox(0,0)[lt]{\lineheight{0}\smash{\begin{tabular}[t]{l}$\mathbf{a}$\end{tabular}}}}%
      % \put(0.4,0.1){\makebox(0,0)[lt]{\lineheight{0}\smash{\begin{tabular}[t]{l}$\mathbf{b}$\end{tabular}}}}%
      \put(0.81983772,0.16834023){\makebox(0,0)[lt]{\lineheight{0}\smash{\begin{tabular}[t]{l}$\mathbf{c}_i$\end{tabular}}}}%
      \put(0.2766603,0.2911645){\makebox(0,0)[lt]{\lineheight{0}\smash{\begin{tabular}[t]{l}$\Phi_f$\end{tabular}}}}%
    \end{picture}%
  \endgroup%

    \caption{$\Phi_f$ in 2D}
    \label{fig:angle_2d}
  \end{subfigure}
  \hfill
  \begin{subfigure}[b]{0.45\textwidth}
    \centering
    \begingroup%
    \makeatletter%
    \providecommand\color[2][]{%
      \errmessage{(Inkscape) Color is used for the text in Inkscape, but the package 'color.sty' is not loaded}%
      \renewcommand\color[2][]{}%
    }%
    \providecommand\transparent[1]{%
      \errmessage{(Inkscape) Transparency is used (non-zero) for the text in Inkscape, but the package 'transparent.sty' is not loaded}%
      \renewcommand\transparent[1]{}%
    }%
    \providecommand\rotatebox[2]{#2}%
    \newcommand*\fsize{\dimexpr\f@size pt\relax}%
    \newcommand*\lineheight[1]{\fontsize{\fsize}{#1\fsize}\selectfont}%
    \ifx\svgwidth\undefined%
      \setlength{\unitlength}{76.30753164bp}%
      \ifx\svgscale\undefined%
        \relax%
      \else%
        \setlength{\unitlength}{\unitlength * \real{\svgscale}}%
      \fi%
    \else%
      \setlength{\unitlength}{\svgwidth}%
    \fi%
    \global\let\svgwidth\undefined%
    \global\let\svgscale\undefined%
    \makeatother%
    \begin{picture}(1,1.05847052)%
      \lineheight{1}%
      \setlength\tabcolsep{0pt}%
      \put(0,0){\includegraphics[width=\unitlength,page=1]{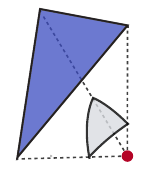}}%
      \put(0.80513996,0.92395654){\makebox(0,0)[lt]{\lineheight{0}\smash{\begin{tabular}[t]{l}$\mathbf{v}_k$\end{tabular}}}}%
      \put(-0.00080604,-0.02){\makebox(0,0)[lt]{\lineheight{0}\smash{\begin{tabular}[t]{l}$\mathbf{v}_j$\end{tabular}}}}%
      \put(0.18569144,1.03){\makebox(0,0)[lt]{\lineheight{0}\smash{\begin{tabular}[t]{l}$\mathbf{v}_i$\end{tabular}}}}%
      \put(0.84445764,0.00865521){\makebox(0,0)[lt]{\lineheight{0}\smash{\begin{tabular}[t]{l}$\mathbf{p}$\end{tabular}}}}%
      \put(0.56,0.25){\makebox(0,0)[lt]{\lineheight{0}\smash{\begin{tabular}[t]{l}$\Phi_f$\end{tabular}}}}%
    \end{picture}%
  \endgroup%

  \caption{$\Phi_f$ in 3D}
    \label{fig:angle_3d}
  \end{subfigure}
  \caption{Signed angles $\Phi_f$ for the 2D and 3D case.
  The point $\mathbf{p}$ is the query point, $\mathbf{c}_i$ and $\mathbf{c}_{i+1}$ are two consecutive vertices representing an edge in a polygon,
  and $\mathbf{v}_i$, $\mathbf{v}_j$ and $\mathbf{v}_k$ are the vertices of an oriented triangle representing a face in a triangle mesh.}
  \label{fig:solid_angles}
\end{figure}

For a discrete geometry, the winding number can be calculated exactly by summing the signed angles between the point $\mathbf{p}$ and the faces of the geometry,
as shown in~\cite{Jacobson2013}:
\begin{equation}
  w(\mathbf{p}) = \frac{1}{\nu \pi} \sum_{f=1}^m \Phi_f,
  \label{eq:sum_w}
\end{equation}
where $\nu = 2$ in 2D and $\nu = 4$ in 3D.
An example of the signed angle $\Phi_f$ for the 2D and 3D cases is depicted in Fig.~\ref{fig:angle_2d} and Fig.~\ref{fig:angle_3d}, respectively.
Note that the result is invariant to the ordering of the faces, as each face contributes independently to $w$.
For the 2D case, $\Phi_f$ is defined as the signed angle between $\mathbf{a} = \mathbf{c}_i - \mathbf{p}$ and
$\mathbf{b} =\mathbf{c}_{i+1} - \mathbf{p}$, where $\mathbf{c}_i$ and $\mathbf{c}_{i+1}$
are two consecutive vertices on the curve,
\begin{equation}
  \Phi_f(\mathbf{p}) = \tan^{-1}\left(\frac{\det(\left[\mathbf{a} \mathbf{b} \right])}{\mathbf{a} \cdot \mathbf{b}}\right).
\end{equation}
Here, $\left[\mathbf{a} \mathbf{b} \right]$ denotes a matrix with $\mathbf{a}$ and $\mathbf{b}$ as column vectors.
In the 3D case, we consider the solid angle of an oriented triangle $\{\mathbf{v}_i, \mathbf{v}_j, \mathbf{v}_k \}$ with respect to $\mathbf{p}$,
which is defined as
\begin{equation}
  \Phi_f(\mathbf{p}) = 2 \tan^{-1}\left(\frac{\det(\left[\mathbf{a} \mathbf{b} \mathbf{c} \right])}{
    \Vert\mathbf{a}\Vert \Vert\mathbf{b}\Vert \Vert\mathbf{c}\Vert + (\mathbf{a} \cdot \mathbf{b}) \Vert\mathbf{c}\Vert  +
    (\mathbf{b} \cdot \mathbf{c}) \Vert\mathbf{a} \Vert  + (\mathbf{c} \cdot \mathbf{a})\Vert\mathbf{b}\Vert }\right),
\end{equation}
where $\mathbf{a}, \mathbf{b}, \mathbf{c}$ are the vectors from $\mathbf{p}$ to the vertices $\mathbf{v}_i$, $\mathbf{v}_j$ and $\mathbf{v}_k$.

\subsubsection{Hierarchical winding}
\label{sec:hierarchical_winding}

The summation property in Eq.~\eqref{eq:sum_w} can be exploited to compute the winding number efficiently.
Let $\mathcal{S}$ be an open surface and let $\bar{\mathcal{S}}$ be an arbitrary closing surface satisfying
\begin{equation}
\partial \bar{\mathcal{S}} = \partial \mathcal{S}.
\end{equation}
Then $\partial\Omega = \bar{\mathcal{S}} \cup \mathcal{S}$ forms a closed oriented surface.
For any query point $\mathbf{p}$ outside of $\partial\Omega$, the following holds:
\begin{equation}
w_{\mathcal{S}}(\mathbf{p}) + w_{\bar{\mathcal{S}}}(\mathbf{p}) = w_{\partial\Omega}(\mathbf{p}) = 0.
\end{equation}
As long as $\mathbf{p}$ is outside of $\partial\Omega$, this implies
\begin{equation}
w_{\mathcal{S}}(\mathbf{p}) = - w_{\bar{\mathcal{S}}}(\mathbf{p}),
\end{equation}
regardless of how $\bar{\mathcal{S}}$ is constructed.
This property can be used in the discrete case to efficiently compute the winding number of a query point
by partitioning the polygon or mesh into a part consisting of a small number of faces and a large remaining part.
For the small part, the winding number is computed directly, and for the large part, a small closing is constructed and its winding number is computed.
The construction of the closing surface is described in more detail in \cite{Jacobson2013,trixiparticles_documentation}.
To efficiently partition the geometry into small parts,
a hierarchy of bounding boxes is constructed by starting with the whole domain and recursively splitting it into two equally sized boxes until the number of faces per box is sufficiently small.
The resulting hierarchy is a binary tree.
The algorithm (\cite{Jacobson2013}, Algorithm~2, p.~5) traverses this binary tree recursively until the leaf containing the query point is found.
The recursion stops with the following criteria:
\begin{itemize}
  \item If the bounding box \(T\) is a leaf, then \(\mathcal{S} \cap T\), the part of \(\mathcal{S}\)
  that lies inside \(T\), is the ``small part'' mentioned above, thus the winding number is evaluated as \(w(\mathbf{p}, \mathcal{S} \cap T)\).
  \item If \(\mathbf{p}\) is outside \(T\), then \(\mathcal{S} \cap T\) is the ``large part'', thus the winding number is evaluated as
  \(-w(\mathbf{p}, \bar{\mathcal{S}} \cap T)\), where \(\bar{\mathcal{S}} \cap T\) is the closing surface of \(\mathcal{S} \cap T\).
\end{itemize}
\begin{figure}[!h]
  \centering
  \begin{subfigure}[b]{0.48\textwidth}
    \centering
    \includegraphics[width=0.5\textwidth]{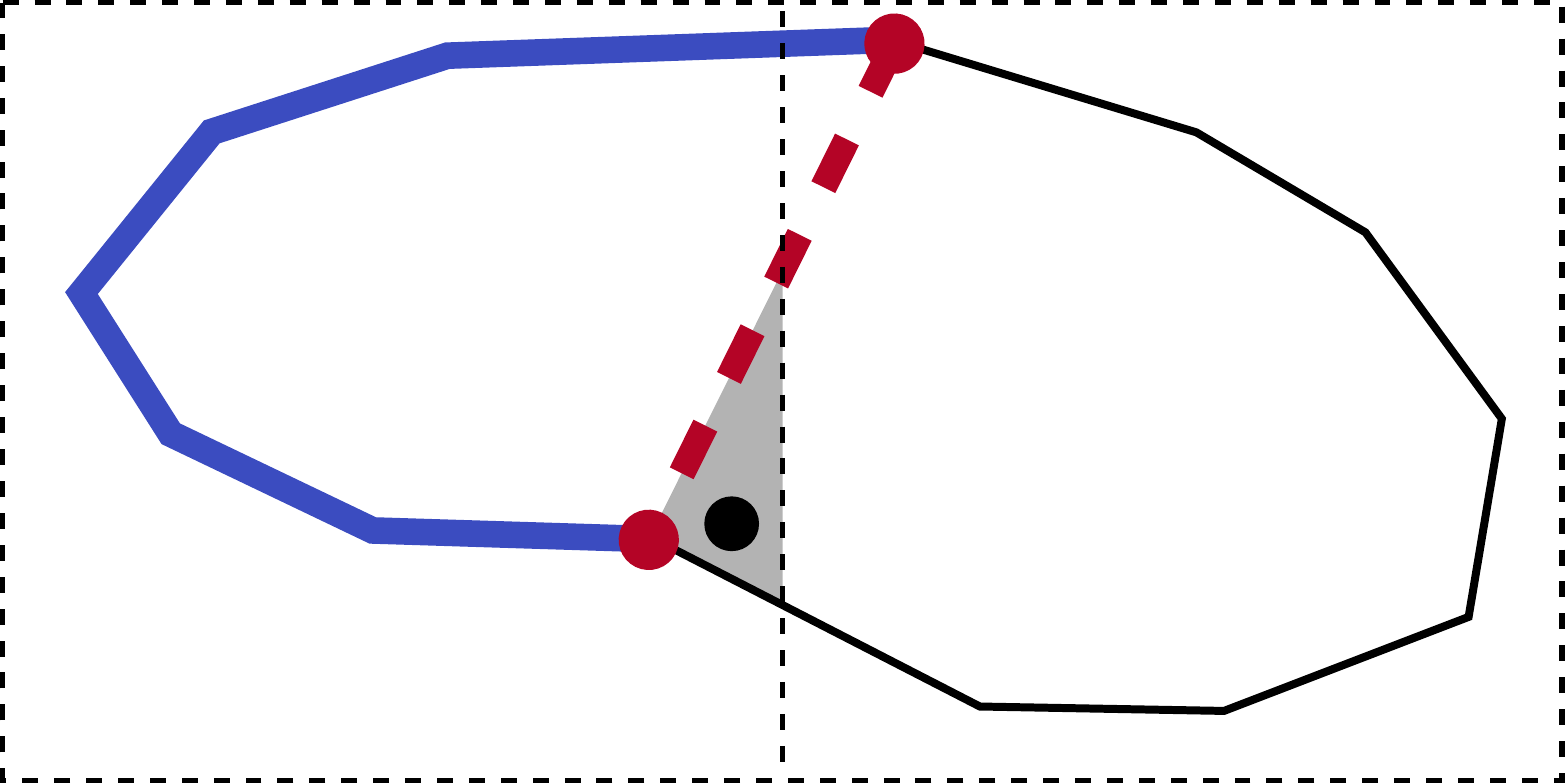}
    \caption{Incorrect evaluation.}
    \label{fig:hierarchical_winding_evaluation_a}
  \end{subfigure}
  \hfill
  \begin{subfigure}[b]{0.48\textwidth}
    \centering
    \includegraphics[width=0.5\textwidth]{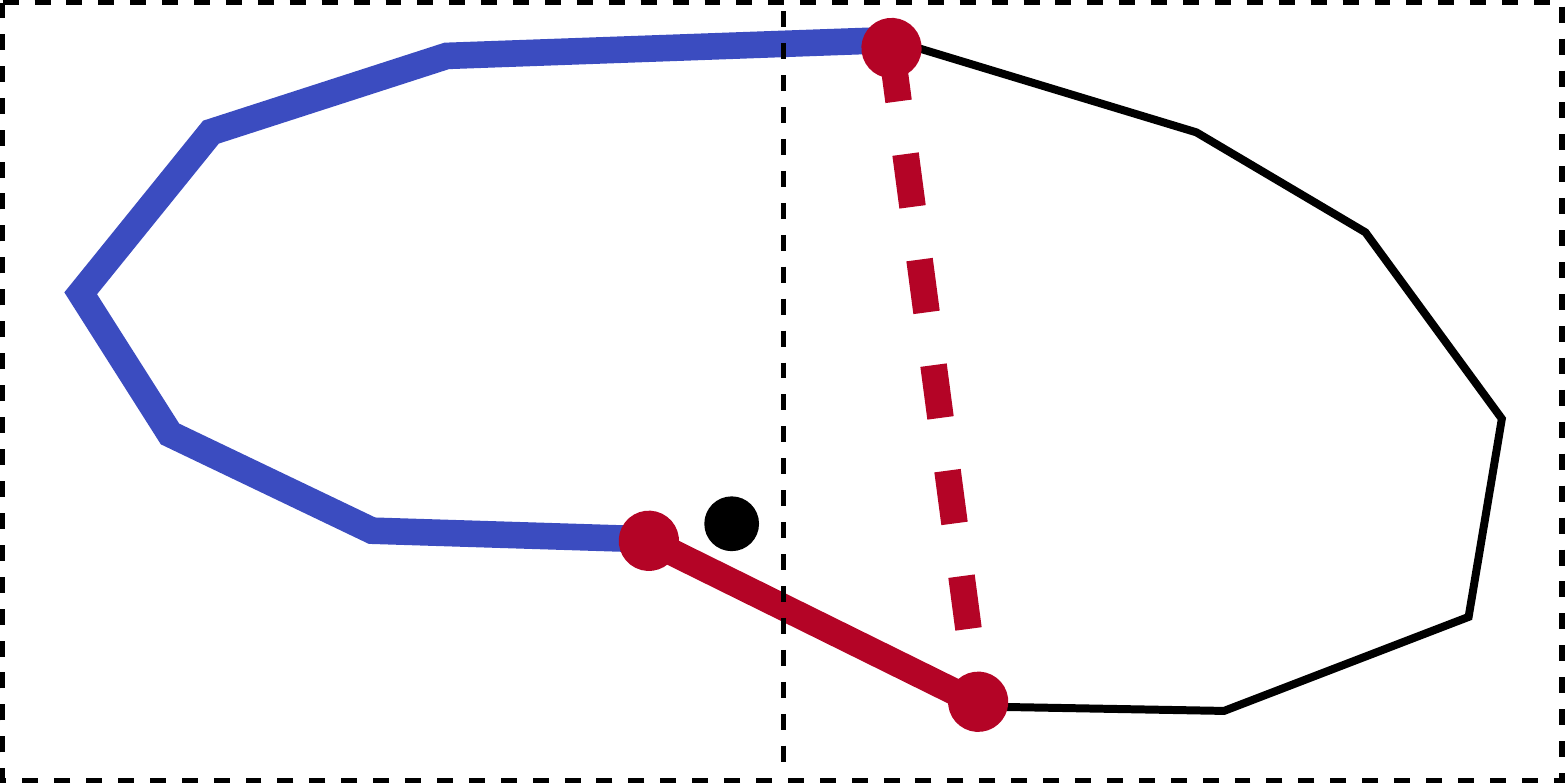}
    \caption{Correct evaluation.}
    \label{fig:hierarchical_winding_evaluation_b}
  \end{subfigure}
  \caption{Evaluation of a query point (black dot) in an arbitrary geometry. The dashed boxes represent the bounding boxes of the hierarchy.
  The blue line represents the faces inside the left bounding box.
  The dashed red line represents the closing determined for the right bounding box and the solid red line represents the intersecting face.
  Fig~(a) shows an incorrect evaluation of the query point in the grey area.
  Fig~(b) shows the correct evaluation due to a different closing.}
  \label{fig:hierarchical_winding_evaluation}
\end{figure}
For discrete geometries, the algorithm described above can lead to an incorrect evaluation of the query points
as can be seen in Fig.~\ref{fig:hierarchical_winding_evaluation_a}.
The query point (black dot) lies in the grey area, which is inside the geometry but is incorrectly classified as outside due to the calculated closing (red dashed line) of the right bounding box.
To resolve this issue, as clarified through personal communication with Jacobson~\citetext{Correspondence with the author of Jacobson et al.~\cite{Jacobson2013}, May 2024},
resizing the bounding boxes ensures a correct evaluation as the exterior vertices forming the closing surface remain within the bounding box.

We avoid resizing by following a different approach and calculating the closure of the bounding box as follows.
We exclude intersecting faces, represented as solid red lines in Fig.~\ref{fig:hierarchical_winding_evaluation_b}, in the calculation of the exterior vertices forming the closing surface.
This way, all exterior vertices are inside the bounding box (right bounding box in Fig.~\ref{fig:hierarchical_winding_evaluation_b}),
and so will be the closing surface (dashed red line).
The intersecting faces are later added with flipped orientation, so that the closing is actually a closing of the exterior plus intersecting faces.

\subsubsection{Robustness}
\label{sec:winding_robustness}
For a completely watertight geometry, the winding number exhibits a jump discontinuity across the surface, i.e., it takes on the values $0$ or $\pm 1$.
In practical applications, however, geometries are often corrupted or not completely watertight.
While other approaches, e.g., ray casting, often fail in such cases unless special provisions are made, the winding number can robustly accommodate holes, non-manifold attachments, and duplicated faces~\cite{Jacobson2013}.
This is achieved by introducing a relaxation factor $\epsilon_w$, which allows a point to be considered inside if its winding number falls within a specified range.
\begin{figure}[!h]
  \centering
  \begin{subfigure}[b]{0.3\textwidth}
    \centering
    \includegraphics[width=0.75\textwidth]{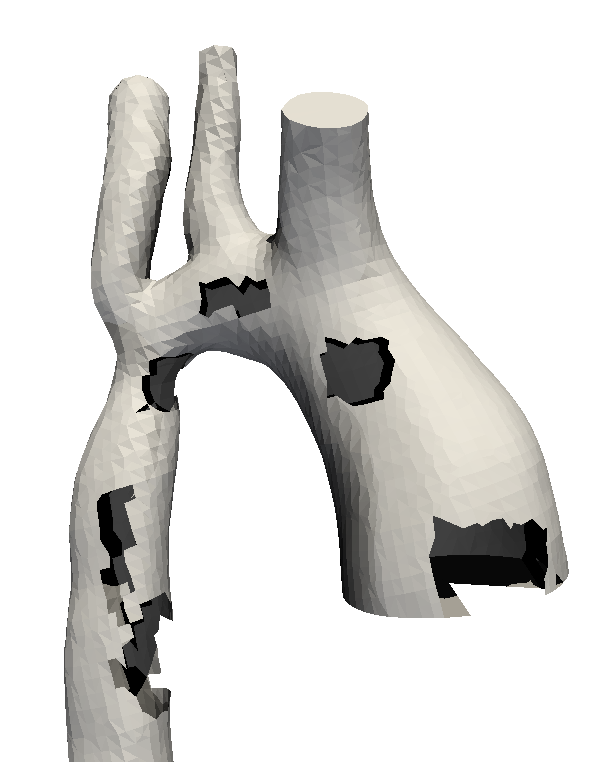}
    \caption{Corrupted geometry}
    \label{fig:corrupted_aorta_a}
  \end{subfigure}
  \hfill
  \begin{subfigure}[b]{0.3\textwidth}
    \centering
    \includegraphics[width=0.75\textwidth]{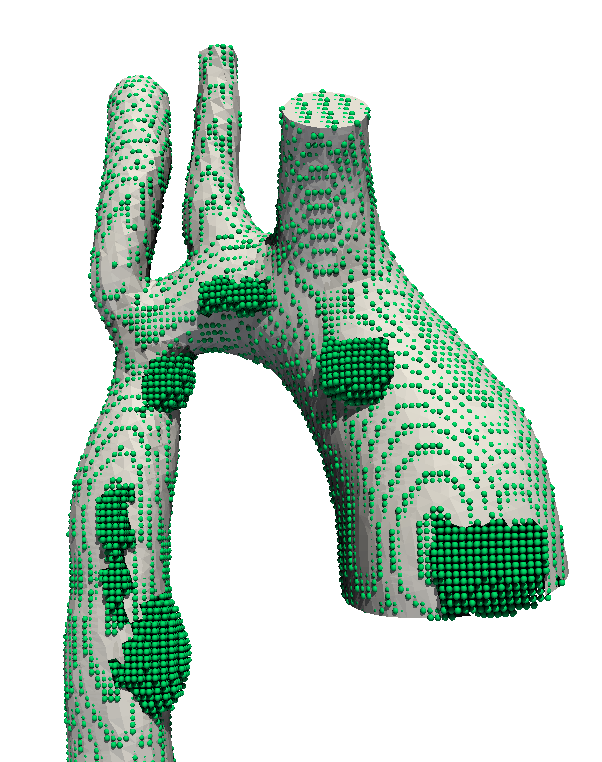}
    \caption{$\epsilon_w = 0.2$}
    \label{fig:corrupted_aorta_b}
  \end{subfigure}
  \hfill
  \begin{subfigure}[b]{0.3\textwidth}
    \centering
    \includegraphics[width=0.75\textwidth]{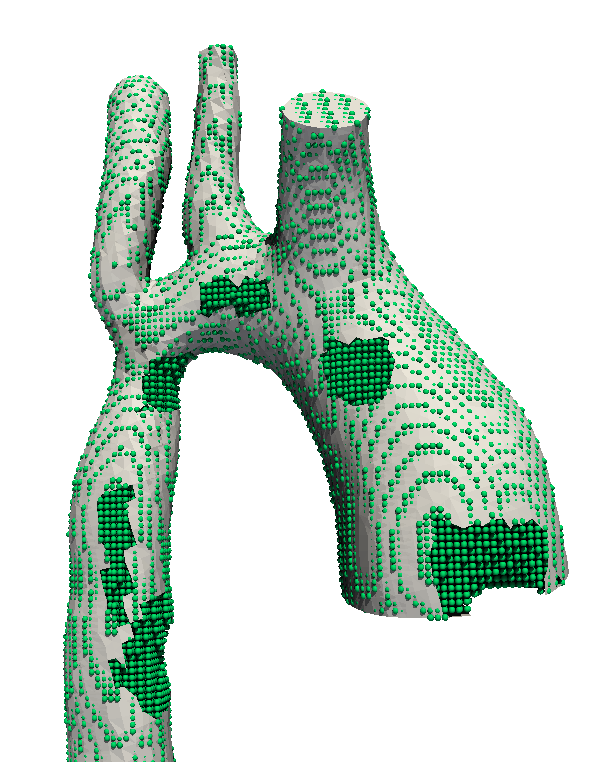}
    \caption{$\epsilon_w = 0.4$}
    \label{fig:corrupted_aorta_c}
  \end{subfigure}
  \caption{Corrupted geometry and the effect of the winding number relaxation factor $\epsilon_w$ on the inside-outside segmentation.
  The green particles represent query points that are considered inside of the geometry.}
  \label{fig:corrupted_aorta}
\end{figure}

Fig.~\ref{fig:corrupted_aorta_a} shows a corrupted geometry with holes.
The hierarchical winding number algorithm was applied to this geometry.
The green particles in Fig.~\ref{fig:corrupted_aorta_b} and Fig.~\ref{fig:corrupted_aorta_c}  satisfy the condition
\begin{equation}
  \mathtt{in\_poly} = \left| w(\mathbf{p}) \right| \geq \epsilon_w.
\end{equation}
It can be seen that for a relaxed condition, the geometry can still be represented.
A low $\epsilon_w$ tends to produce convex results (center), whereas a high $\epsilon_w$ tends to produce concave results (right).

% !TeX root = ../main.tex
\section{Particle Packing}
\label{sec:packing}
\begin{figure}[!h]
    \centering
    \begin{subfigure}[b]{0.48\textwidth}
      \centering
      \includegraphics[width=0.75\textwidth]{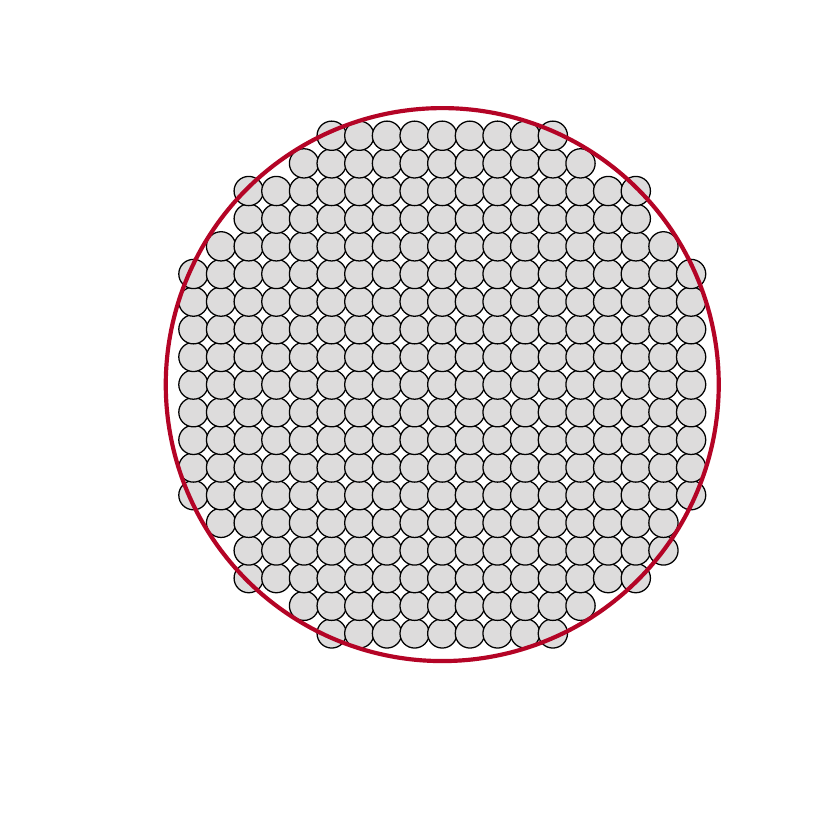}
      \caption{Initial configuration.}
      \label{fig:circle_sampled_a}
    \end{subfigure}
    \hfill
    \begin{subfigure}[b]{0.48\textwidth}
      \centering
      \includegraphics[width=0.75\textwidth]{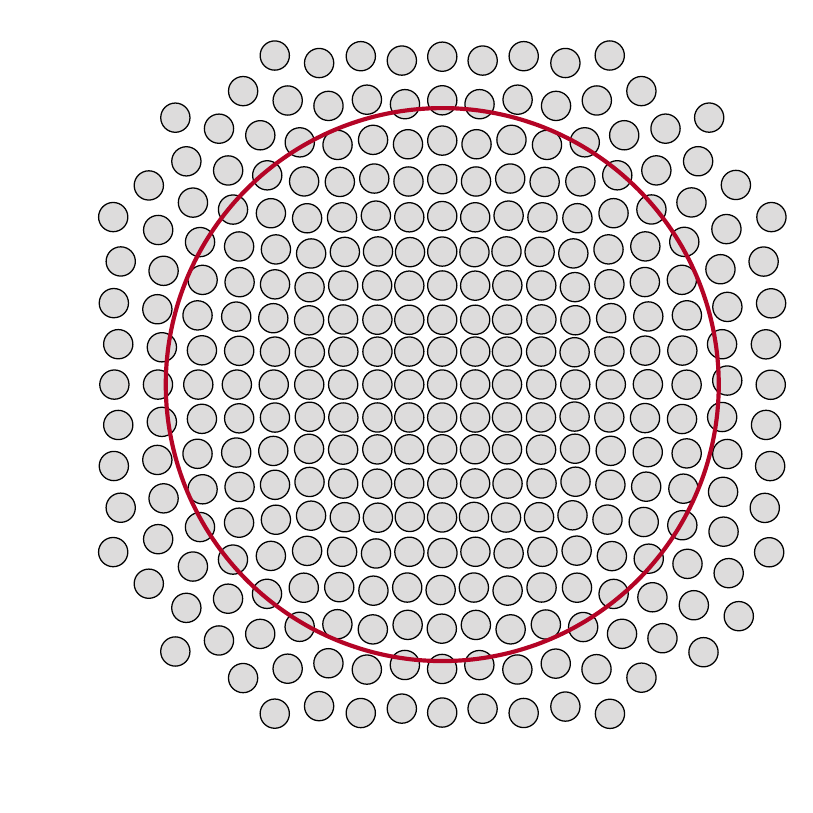}
      \caption{Final configuration.}
      \label{fig:circle_sampled_b}
    \end{subfigure}
    \caption{(a) Initial particle configuration obtained after segmentation.
    (b) Particle configuration after applying the packing algorithm without enforcing boundary constraints.}
    \label{fig:circle_sampled}
  \end{figure}

In this section, we focus on the final step (5), the packing of particles.
The initial particle configuration forms an isotropic grid distribution within the interior, as illustrated in Fig.~\ref{fig:circle_sampled_a}.
However, along the boundaries, particles must be aligned with the surface while the particle distribution should remain isotropic throughout the entire domain.
To achieve this, we employ a packing method developed by Zhu et al.~\cite{Zhu2021}.
In this method, the acceleration of each particle is computed from the pressure force (see Eq.~\eqref{eq:pressure_force}), and the particle positions are then updated accordingly.
At equilibrium, the resultant pressure forces acting on each particle vanish, yielding a stable and isotropic particle distribution.
Consequently, when the particles are correctly arranged, the accuracy of the kernel gradient estimation is improved.
That is, for a constant pressure field, the gradient is zero as is shown in Eq.~\eqref{eq:zero_gradient}.

\subsection{Bounding method}
\label{sec:bounding_method}
Since packing moves particles towards areas of lower particle density, particles at the boundaries drift away when not constrained, as depicted in Fig.~\ref{fig:circle_sampled_b}.
To address this, a bounding method \cite{Zhu2021} is applied,
where after each time step the particle positions $\mathbf{r_i}$ are corrected as:
\begin{equation}
    \mathbf{\tilde{r}}_{i} =
    \begin{cases}
      \mathbf{r}_i  - \left( \phi_i + \frac{1}{2} \Delta x \right) \mathbf{n}_i & \text{if } \phi_i \geq - \frac{1}{2} \Delta x, \\
      \mathbf{r}_i & \text{otherwise.}
    \end{cases}
    \label{eq:bounding method}
\end{equation}
Here, $\phi_i$ and $\mathbf{n}_i$ are the signed distance and the normal direction to the surface at the position of particle~$i$, respectively.
For methods where particles need to be placed on the boundary instead of half a particle spacing away, notably the Total Lagrangian SPH (TLSPH) method, as presented in \cite{O_Connor_2021}, the update condition in Eq.~\eqref{eq:bounding method} changes accordingly to $\phi_i \geq 0$.
To determine $\phi$ for each particle, Zhu et al.~\cite{Zhu2021} use trilinear interpolation to interpolate the level-set field at the particle's position from a background mesh.
In this work we use Shepard interpolation~\cite{Shepard1968}.
That is, $\phi_i$ and $\mathbf{n}_i$ are interpolated from the point cloud generated in Sec.~\ref{sec:determine_signed_distances} according to Eq.~\eqref{eq:interpolate phi} and Eq.~\eqref{eq:interpolate n}.

\subsection{Dynamic boundary packing}
\label{sec:dynamic_boundary_packing}
Another challenge arises at the boundaries, where particles lack full kernel support.
In this case, the gradient of the kernel is estimated inaccurately, resulting in a distorted pressure force.
Consequently, the particles experience a stronger repulsive force from the inner region.
This leads to undesirable particle distribution near the surface, as can be seen in Fig.~\ref{fig:circle_packed_a} (with detailed discussion of results in the next section),
which is particularly noticeable in regions of high curvature.
To address this, Yu et al.~\cite{Yu2023a} introduced the static confinement boundary condition.
In this method, the missing kernel support is compensated by incorporating a volume contribution from outside the body, which is calculated from a level-set mesh.
Inspired by their method, we use a similar approach, where we sample the boundary with particles as described in Sec.~\ref{sec:boundary_particle_generation}.
These boundary particles are packed together with the inner particles during relaxation.
This means that the position of the boundary particles is also updated.
The advantage is that a natural full kernel support is ensured, with the boundary particles aligning uniformly.
Additionally, the boundary particles together with the inner particles form an isotropic distribution.
To ensure a proper boundary-interior interface representing the geometry's surface,
the bounding method applied to the boundary particles is adjusted as follows:
\begin{equation}
  \mathbf{\tilde{r}}_{b, i} =
  \begin{cases}
    \mathbf{r}_{b, i}  - \left( \phi_{b, i} - \tau \right) \mathbf{n}_{b, i} & \text{if } \phi_{b, i} \geq \tau + \frac{1}{2} \Delta x_b, \\
    \mathbf{r}_{b, i}  - \left( \phi_{b, i} -  \frac{1}{2} \Delta x_b \right) \mathbf{n}_{b, i} & \text{if } \phi_{b, i} < \frac{1}{2} \Delta x_b, \\
    \mathbf{r}_{b, i} & \text{otherwise}.
  \end{cases}
  \label{eq:bounding method boundary}
\end{equation}
Here, $\tau$ denotes the boundary thickness, $\Delta x_b$ is the initial particle spacing of the boundary particles
and $\phi_{b, i}$ and $\mathbf{n}_{b, i}$ represent the signed distance and the normal direction at the position of boundary particle $i$, respectively.
The values of $\phi_{b, i}$ and $\mathbf{n}_{b, i}$ are computed using the same equations as those used for the inner particles (see Eq.\eqref{eq:interpolate phi} and Eq.~\eqref{eq:interpolate n}).
The first condition in Eq.~\eqref{eq:bounding method boundary} ensures that the boundary particles do not exit the boundary hull to the outside,
the second condition ensures that the boundary particles do not enter the inner domain of the geometry.
For methods where particles are placed directly on the surface, the second update condition in Eq.~\eqref{eq:bounding method boundary} changes accordingly to $\phi_{b,i} < \Delta x_b $
and the first condition changes to $\phi_{b,i} \geq \tau + \Delta x $.

% !TeX root = ../main.tex
\section{Results}
\label{sec:results}
In the following, we adopt a reference density $\rho_{0}$, from which we compute the particle mass by~${m_i = \rho_{0} V_{\Omega} / n_{\Omega}}$,
where $V_{\Omega}$ denotes the volume of the geometry and $n_{\Omega}$ represents the total number of interior particles.
This way, the total mass of the particles is independent of the particle resolution.

The quality of the particle packing is assessed both during and after the relaxation process.
To evaluate the quality during packing, we analyze the kinetic energy of the system, which indicates the magnitude of the zero-order consistency error, providing insight into the quality of the particle distribution during the packing process \cite{Zhu2021}.
The kinetic energy is defined as
\begin{equation}
 E_{\text{kin}} = \frac{1}{2} \sum_{i=1}^{n_{\Omega}} m_{i} \left( \Vert \mathbf{\tilde{v}}_{i} \Vert^2 \right),
\end{equation}
where $\mathbf{\tilde{v}}_{i}$ is the advection velocity.
Since the initial kinetic energy is much smaller for higher resolutions, to observe its evolution we define the
normalized kinetic energy by
\begin{equation}
 E_{\text{kin}, \text{n}} = \frac{E_{\text{kin}}}{E_{\text{kin}, \text{max}}},
\end{equation}
where $E_{\text{kin}, \text{max}}$ is the maximum kinetic energy during packing.

To evaluate the final particle distribution, the density $\rho_i$ at each particle position is computed using the kernel summation (via Eq.~\eqref{eq:summation_density}),
including boundary particles to ensure full kernel support.
We then quantify the deviation from the reference density using the $L_\infty$ and $L_2$ norms defined as:
\begin{equation}
  L_\infty = \max_{i} \left| \rho_i - \rho_0 \right|, \quad
  \label{eq:L_infty}
\end{equation}
\begin{equation}
  L_2 = \sqrt{\frac{1}{n_{\Omega}} \sum_{i=1}^{n_{\Omega}} \left( \rho_i - \rho_0 \right)^2},
  \label{eq:L_2}
\end{equation}
where $\rho_0 = 1$ denotes the reference density.
Note that for all visualizations of particle distributions presented in this paper, the marker size in the plots is set precisely to the corresponding particle spacing $\Delta x$.
This approach ensures that the quality of the packing can be directly and consistently assessed from the figures.

This section is structured as follows:
In the first part, we examine the qualitative result of the particle packing using 2D geometries.
Here, we focus on the influence of the smoothing length and the presence of boundary particles.
We also analyze the evolution of the kinetic energy during packing and demonstrate the robustness of the method.
In the second part, we shift the focus to a quantitative evaluation of the final particle configurations.
To this end, we present results for 3D geometries, investigating both the deviation from the reference density and the convergence behavior of the method.
Additionally, a multi-body scenario is presented to highlight the versatility of the proposed approach.
Finally, to validate the practical applicability of our method,
we perform fluid-structure interaction (FSI) simulations of an airfoil,
comparing lattice-based and packed particle configurations across different resolutions.
For reproducibility, the numerical setups for all results can be found in our reproducibility repository \cite{Neher2025reproducibility}.

\subsection{Particle packing in 2D}
First, we present qualitative packing results for a simple circle,
where we investigate the influence of the smoothing length and the inclusion of boundary particles on the packing quality.
Subsequently, we show results for geometries exhibiting high curvature and sharp corners using airfoil profiles.
\subsubsection{Circle}
\begin{figure}[!h]
  \centering
  \begin{subfigure}[b]{0.24\textwidth}
    \centering
    \includegraphics[width=0.9\textwidth]{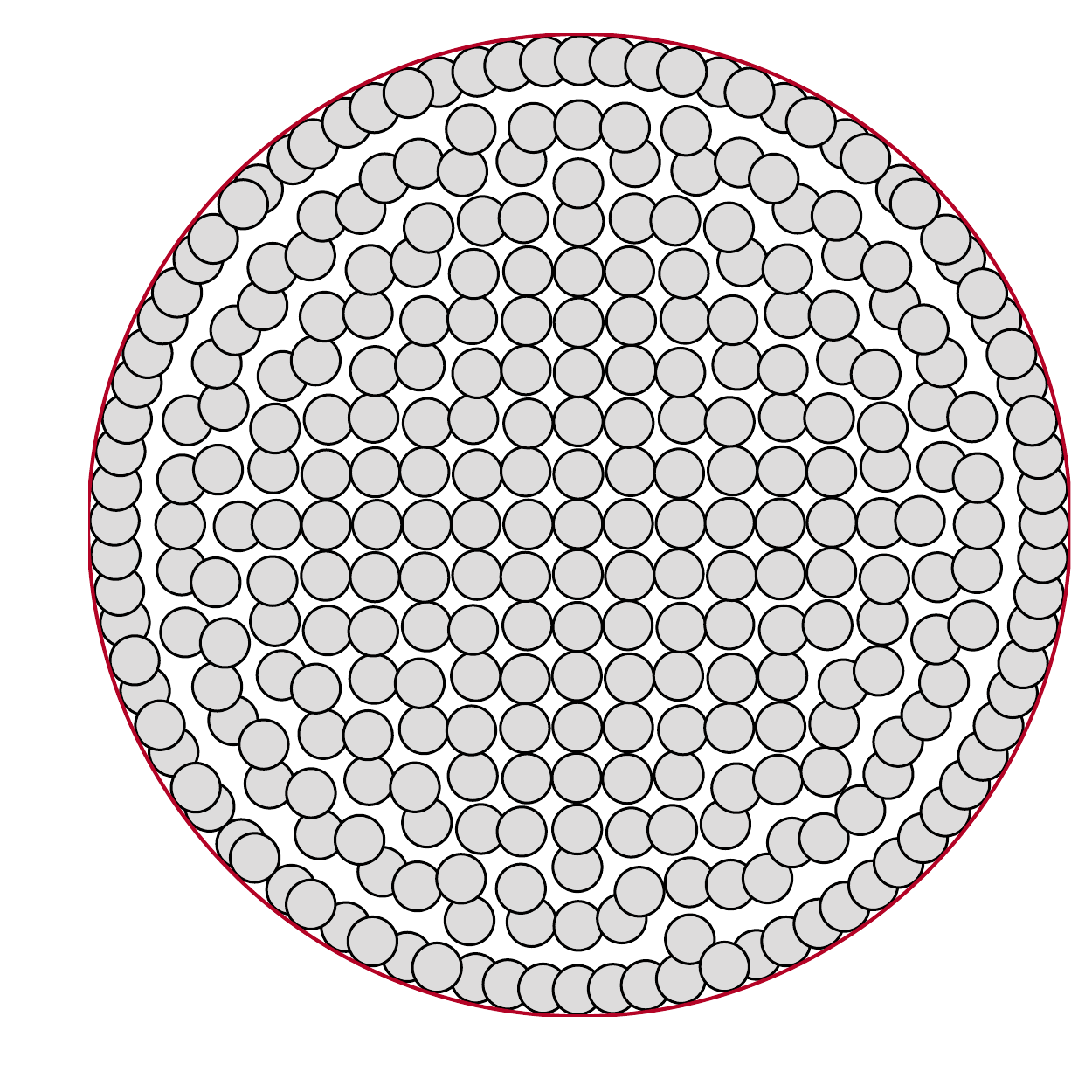}
    \caption{$h=1.2 \Delta x$,\\no boundary.}
    \label{fig:circle_packed_a}
  \end{subfigure}%
  \hfill
  \begin{subfigure}[b]{0.24\textwidth}
    \centering
    \includegraphics[width=0.9\textwidth]{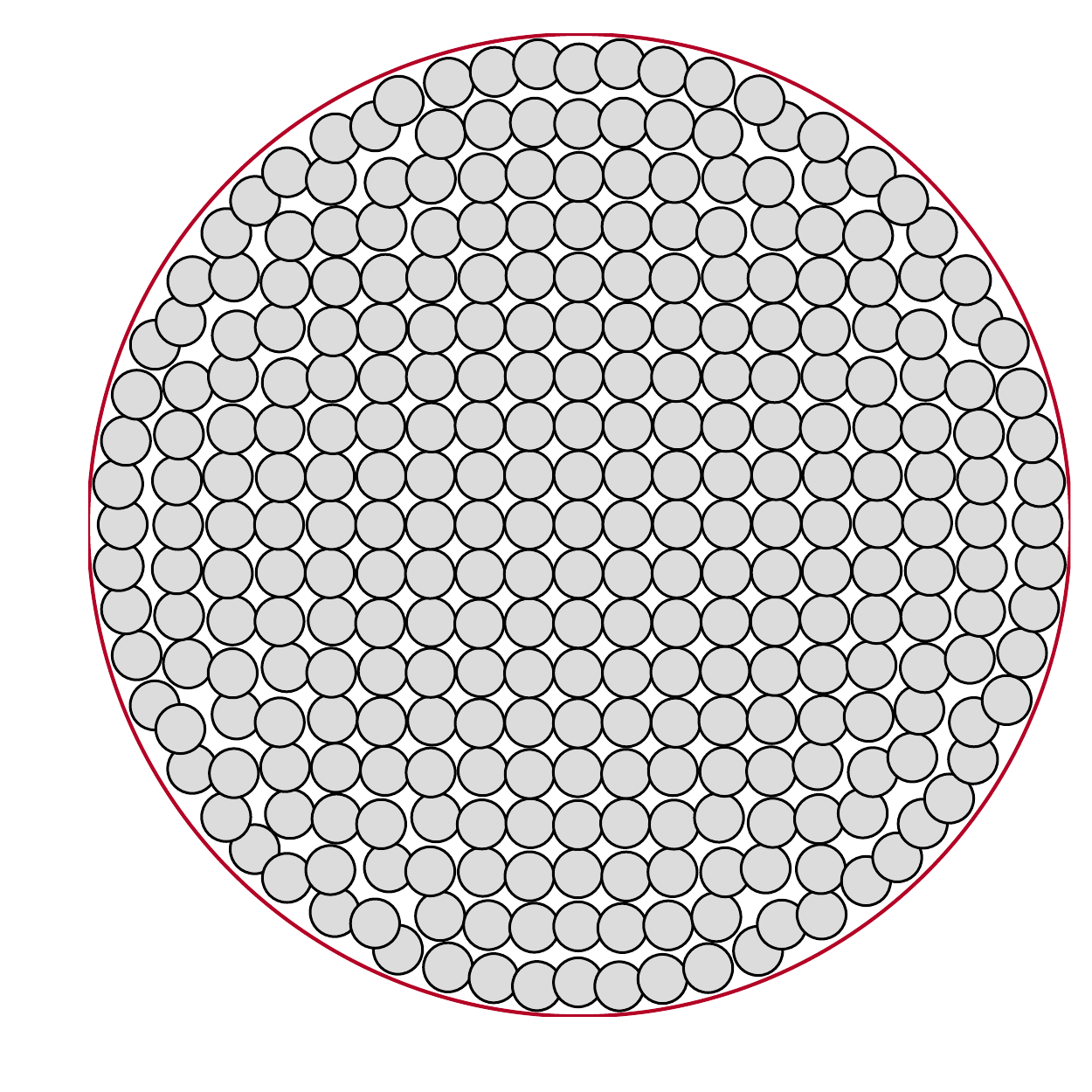}
    \caption{$h=1.2 \Delta x$,\\with boundary.}
    \label{fig:circle_packed_b}
  \end{subfigure}%
  \hfill
  \begin{subfigure}[b]{0.24\textwidth}
    \centering
    \includegraphics[width=0.9\textwidth]{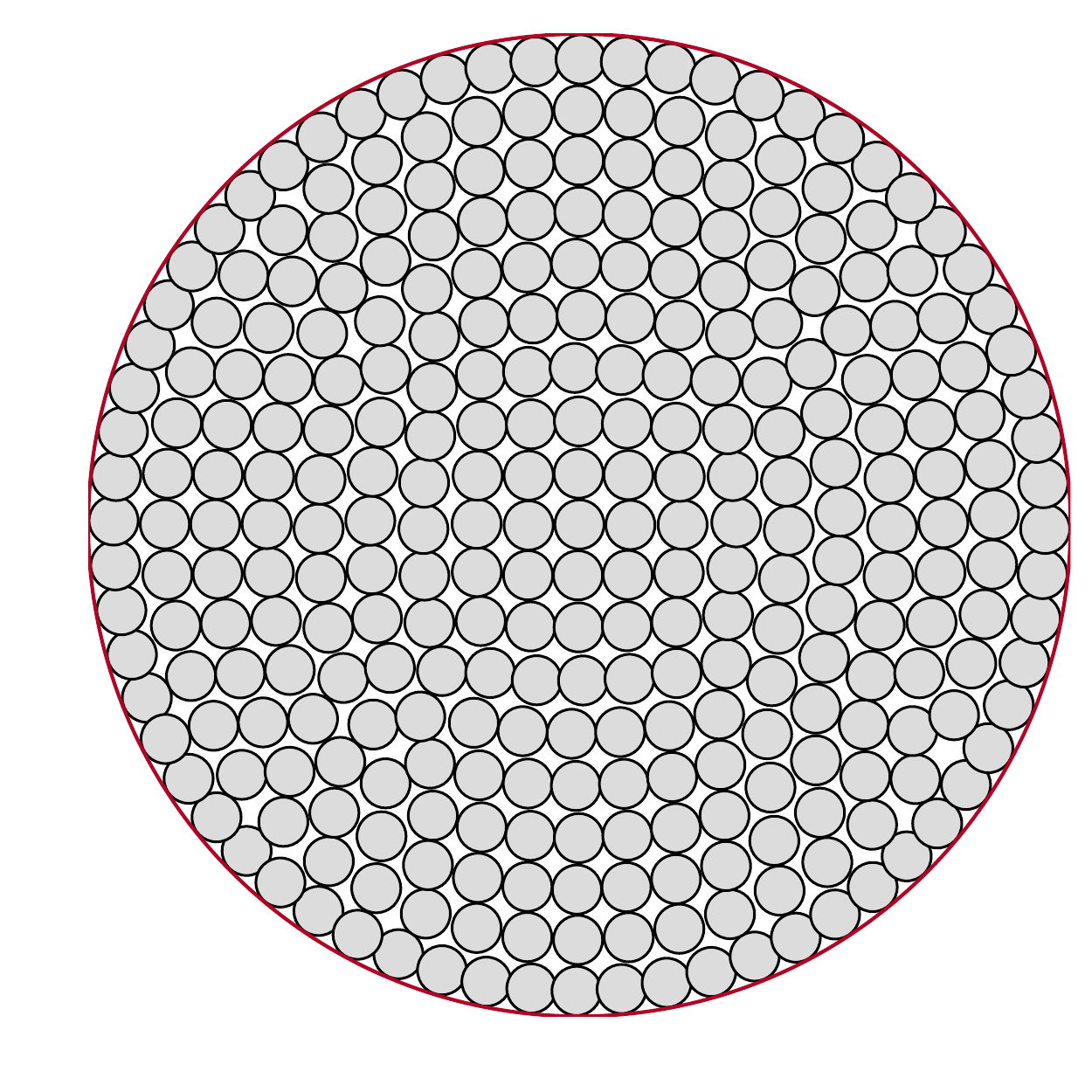}
    \caption{$h=0.8 \Delta x$,\\no boundary.}
    \label{fig:circle_packed_c}
  \end{subfigure}%
  \hfill
  \begin{subfigure}[b]{0.24\textwidth}
    \centering
    \includegraphics[width=0.9\textwidth]{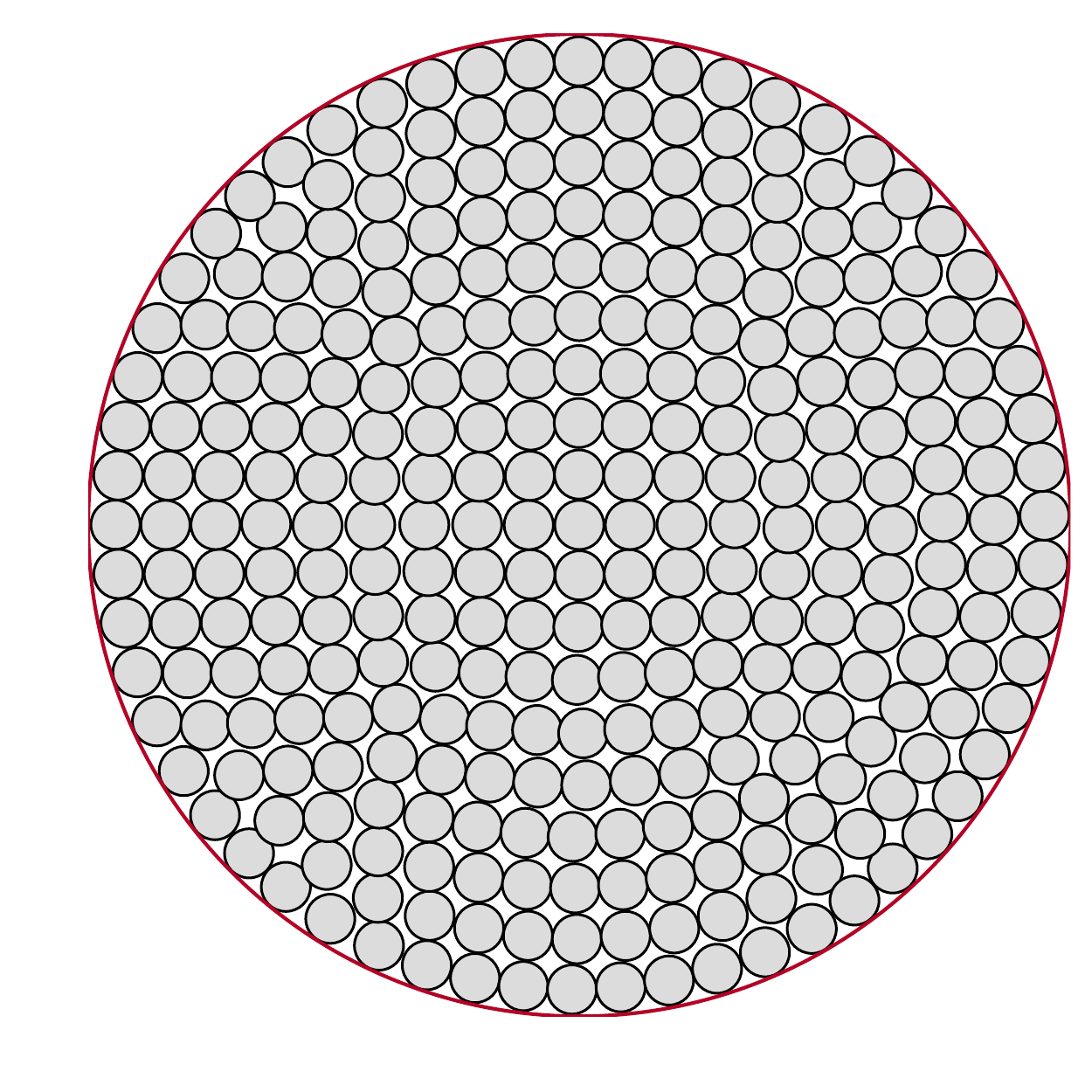}
    \caption{$h=0.8 \Delta x$,\\with boundary.}
    \label{fig:circle_packed_d}
  \end{subfigure}%
  \caption{Comparison of a packed circle geometry, represented with a red line, with and without boundary particles and with different smoothing lengths $h$. $\Delta x$ is the initial particle spacing.}
  \label{fig:circle_packed}
\end{figure}

\begin{figure}[!h]
  \center
  \includegraphics[width=0.7\textwidth]{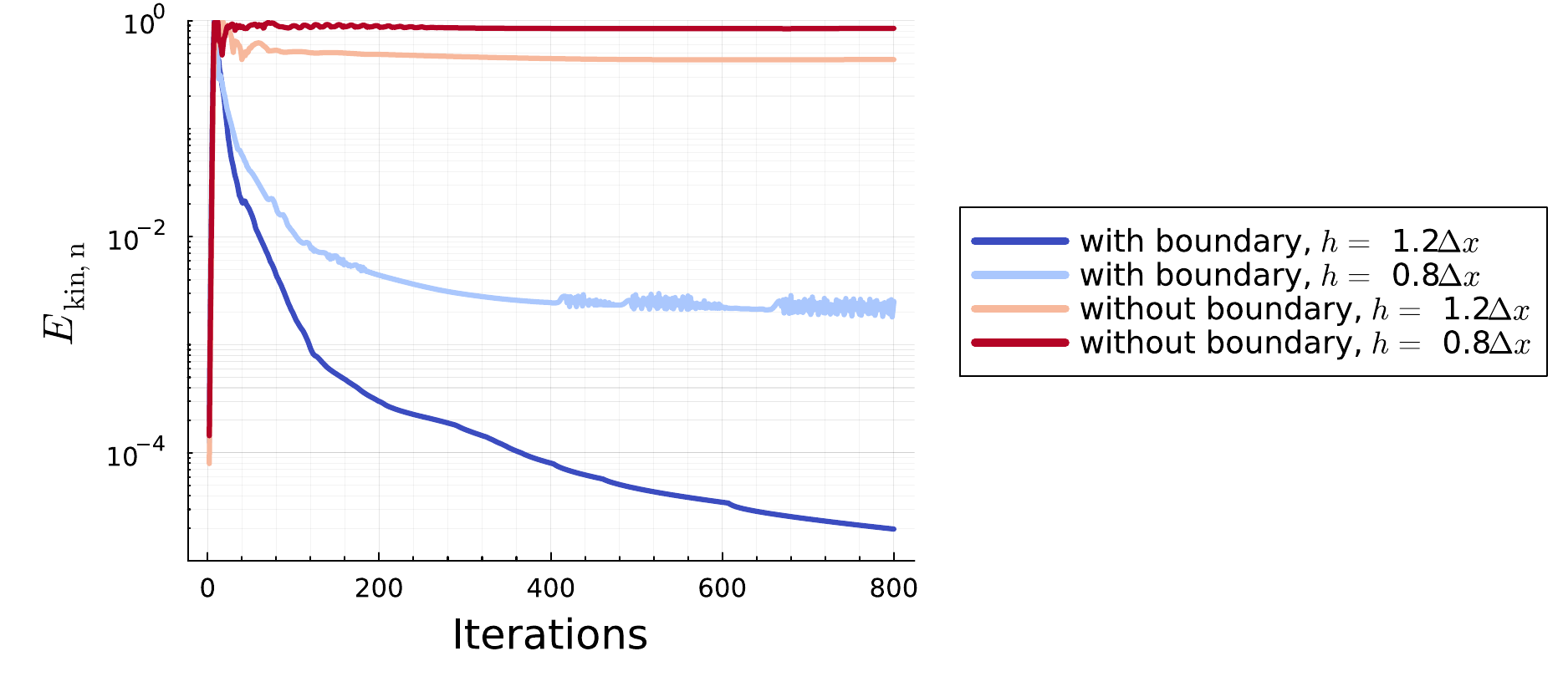}
  \caption{Normalized kinetic energy $E_{\mathrm{kin}, \mathrm{n}}$ during packing for the circle with and without boundary particles, and for different smoothing lengths $h$.
  The packing process is performed over 800 iterations.}%
  \label{fig:circle_kinetic_energy}
\end{figure}

\begin{figure}[!h]
  \centering
  \begin{subfigure}[b]{0.3\textwidth}
    \centering
    \includegraphics[width=0.9\textwidth]{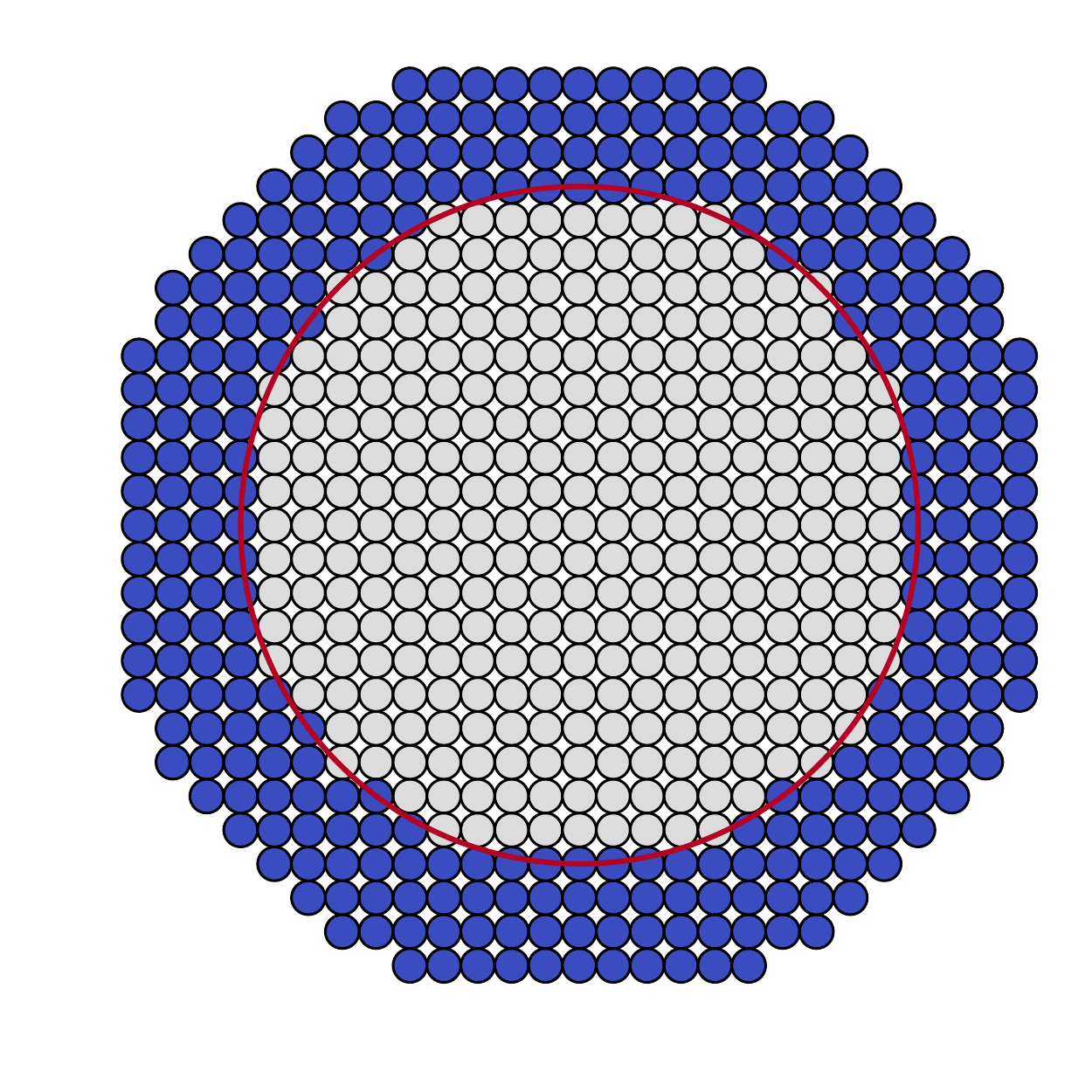}
    \caption{Initial configuration}
    \label{fig:sampled_boundary_a}
  \end{subfigure}%
  \hfill
  \begin{subfigure}[b]{0.3\textwidth}
    \centering
    \includegraphics[width=0.9\textwidth]{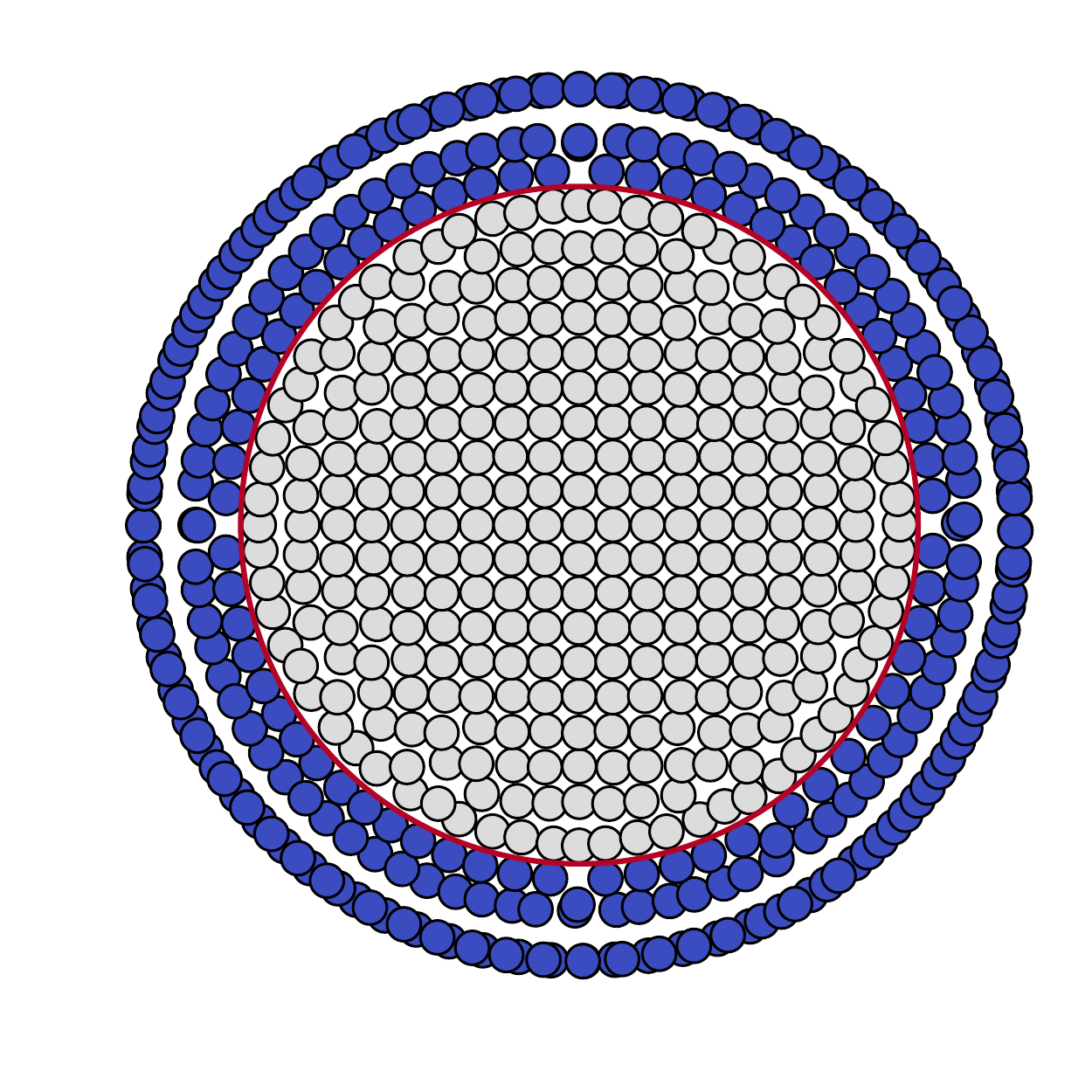}
    \caption{ $h=1.2 \Delta x$}
    \label{fig:sampled_boundary_b}
  \end{subfigure}%
  \hfill
  \begin{subfigure}[b]{0.3\textwidth}
    \centering
    \includegraphics[width=0.9\textwidth]{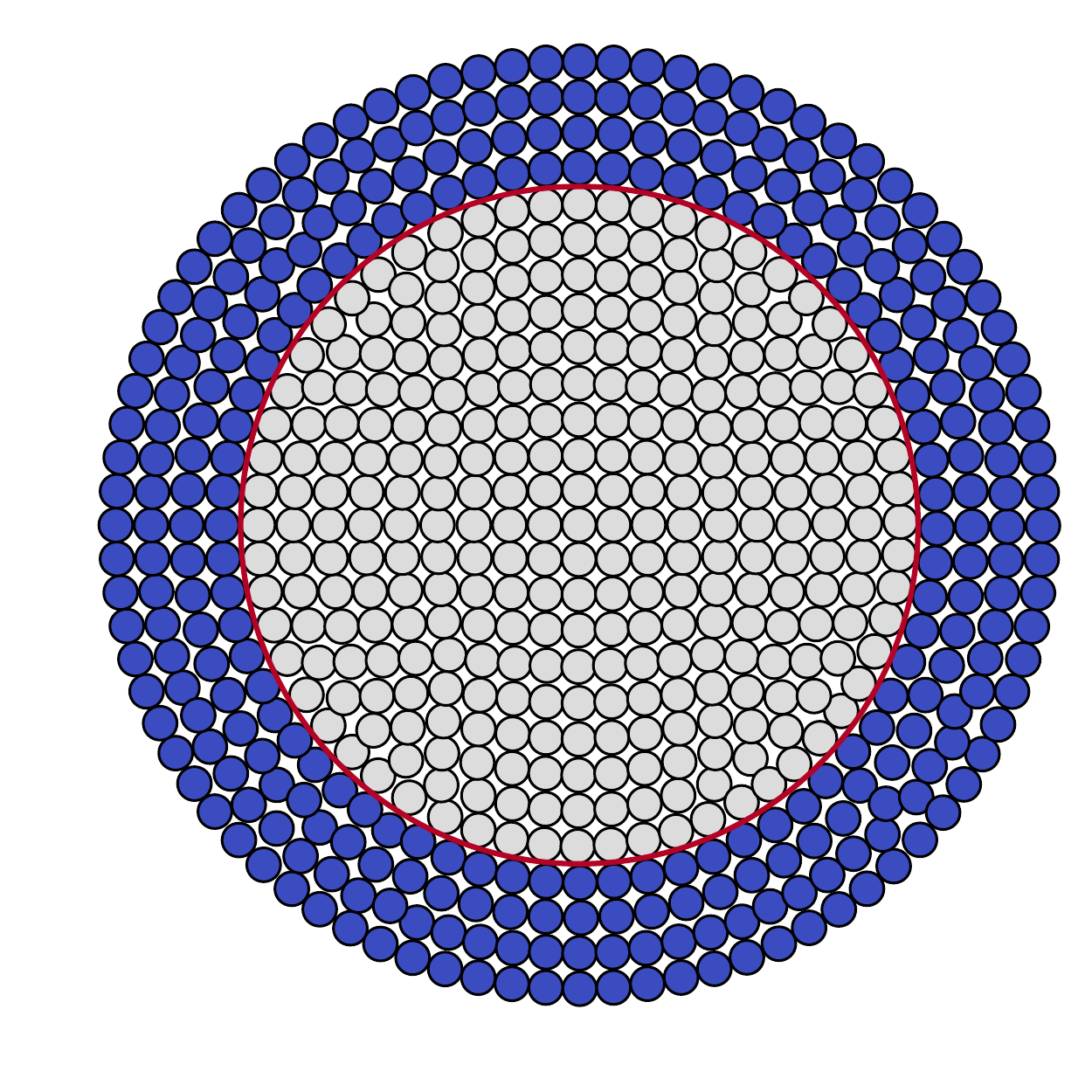}
    \caption{$h=0.8 \Delta x$}
    \label{fig:sampled_boundary_c}
  \end{subfigure}%
  \caption{Comparison of the initial boundary of the circle to the packed boundary with different smoothing lengths.}%
  \label{fig:sampled_boundary}
\end{figure}

\begin{figure}[!h]
  \center
  \includegraphics[width=0.4\textwidth]{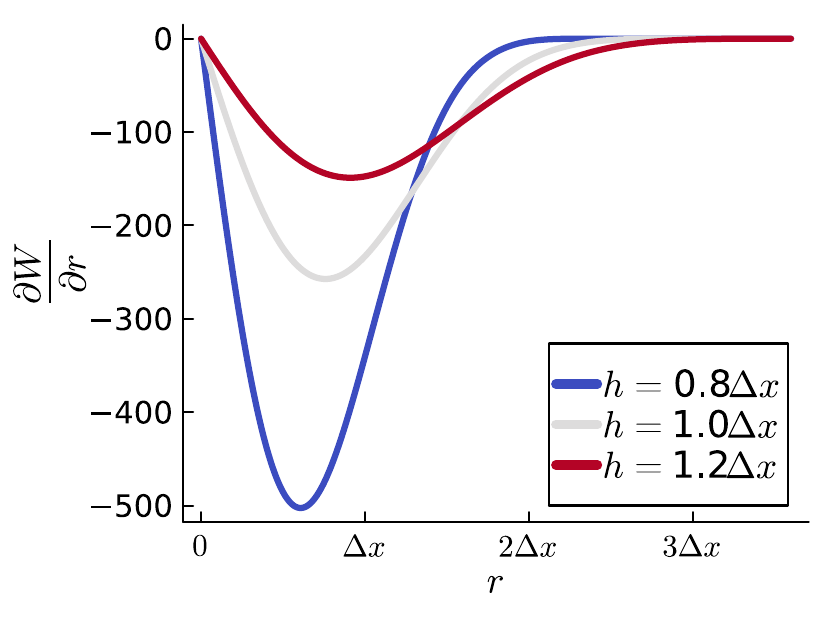}
  \caption{Derivative of the quintic spline kernel for different smoothing lengths $h$.}%
  \label{fig:kernel_derivatives}
\end{figure}

A circle of radius $r=1$ is considered and sampled with a particle spacing of $\Delta x = 0.1$.
The initial particle configuration after sampling is shown in Fig.~\ref{fig:circle_sampled_a}.
For the packing, a quintic kernel~\cite{Schoenberg1946} with $h=1.2 \Delta x$ is used in Fig.~\ref{fig:circle_packed_a} and Fig.~\ref{fig:circle_packed_b}.
In Fig.~\ref{fig:circle_packed_a}, the circle is packed without boundary particles
and in Fig.~\ref{fig:circle_packed_b} with boundary particles.
The qualitative improvement from Fig.~\ref{fig:circle_packed_a} to Fig.~\ref{fig:circle_packed_b} is clearly visible.
In Fig.~\ref{fig:circle_packed_c}, the circle is packed without boundary particles but the smoothing length is reduced to $h=0.8 \Delta x$,
which also improves the particle distribution.
In Fig.~\ref{fig:circle_packed_d}, both strategies are combined, yielding an improved result:
while a smaller smoothing length (Fig.~\ref{fig:circle_packed_c}) produces visually almost similar results to the packing with boundary (Fig.~\ref{fig:circle_packed_d}),
Fig.~\ref{fig:circle_kinetic_energy} demonstrates that packing with boundary particles lowers the kinetic energy by several orders of magnitude.
An explanation for this can be inferred by considering the boundary particles.
The initial configuration of the boundary particles is shown in Fig.~\ref{fig:sampled_boundary_a} in blue.
In Fig.~\ref{fig:sampled_boundary_b}, the packed boundary particles with $h=1.2 \Delta x$ is shown.
It is clear that clumping occurs, which affects the quality of the boundary.
By reducing the smoothing length, the clumping effects can be prevented, as seen in Fig.~\ref{fig:sampled_boundary_c}.
This is due to the properties of the kernel function.
The forces driving the particle packing are described by Eq.~\eqref{eq:pressure_force}, which depends on the gradient of the kernel function.
With smaller smoothing lengths, the minimum of the gradient moves closer to the particle center (see Fig.~\ref{fig:kernel_derivatives}).
Particles that are very close together stay clumped, since the repulsive forces between them are small.
When reducing the smoothing length, particles have to overcome a much larger repulsive force before reaching the clumped state.
While a smaller smoothing length improves the particle distribution,
it might negatively impact kernel estimation due to the reduced number of neighboring particles.
This issue can be circumvented by employing a dedicated, higher smoothing length for the interpolation process (see Eqs.~\eqref{eq:interpolate phi} and \eqref{eq:interpolate n}).
However, our investigations indicated that using two different smoothing lengths has a negligible effect, which is why we employ the same smoothing length for both processes.

Similar to Fig.~\ref{fig:circle_packed_d}, the quality of the outer layers of the boundary could be further improved by adding more particles outside of the boundary.
This would only be used for packing and not the simulation essentially applying the boundary packing to the boundary as well, but we do not expect any noticeable improvement from this extra step.

\subsubsection{Airfoil NACA0015}
\label{sec:airfoil_NACA0015}
To examine how the algorithm behaves under varying resolutions, high curvature, and sharp corners,
we consider a symmetric NACA0015 airfoil.
\begin{figure}[!h]
  \centering
  \begin{subfigure}[b]{0.48\textwidth}
    \centering
    \includegraphics[width=0.95\textwidth]{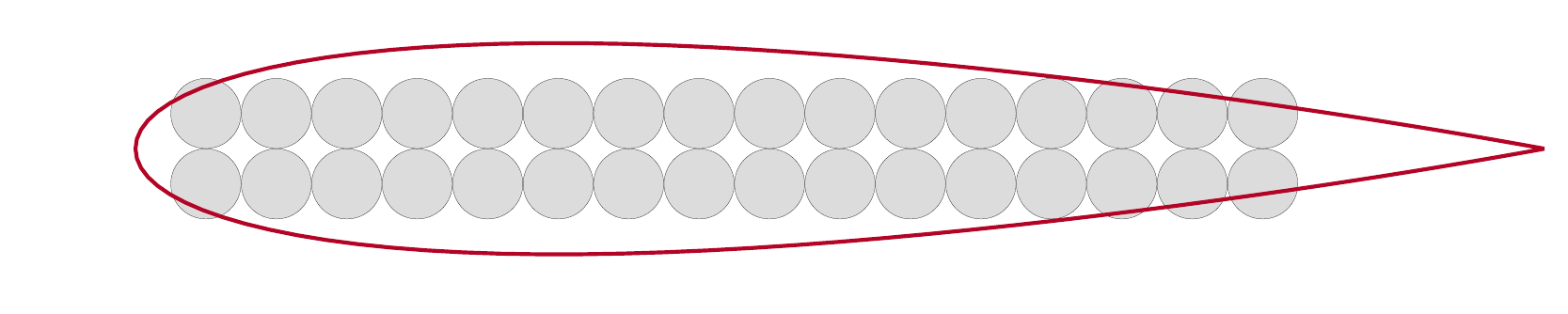}\\
    \includegraphics[width=0.95\textwidth]{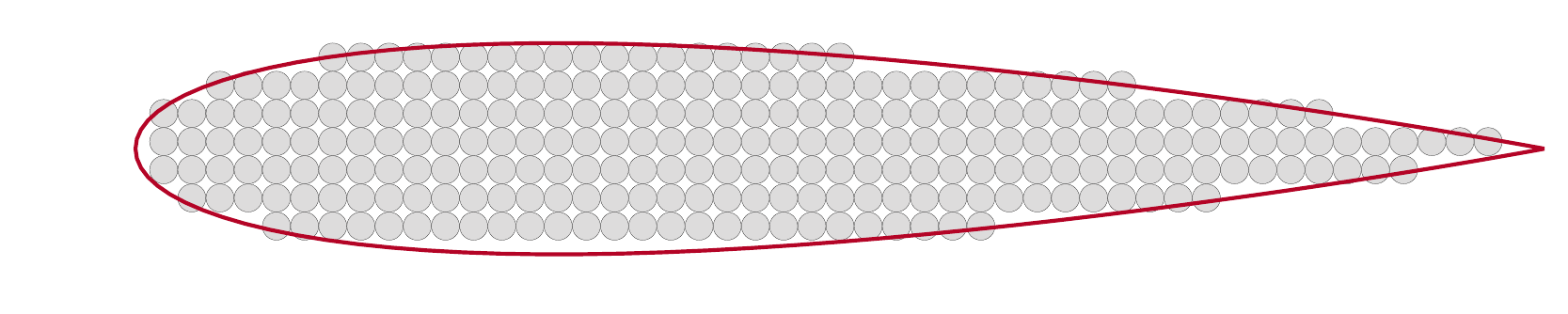}\\
    \includegraphics[width=0.95\textwidth]{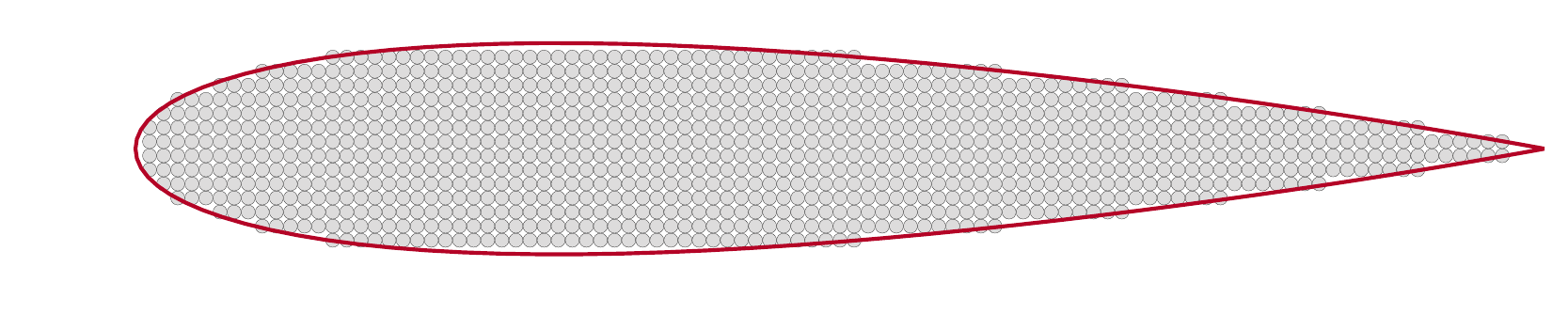}
    \caption{Initial configuration.}
  \end{subfigure}
  \hfill
  \begin{subfigure}[b]{0.48\textwidth}
    \centering
    \includegraphics[width=0.95\textwidth]{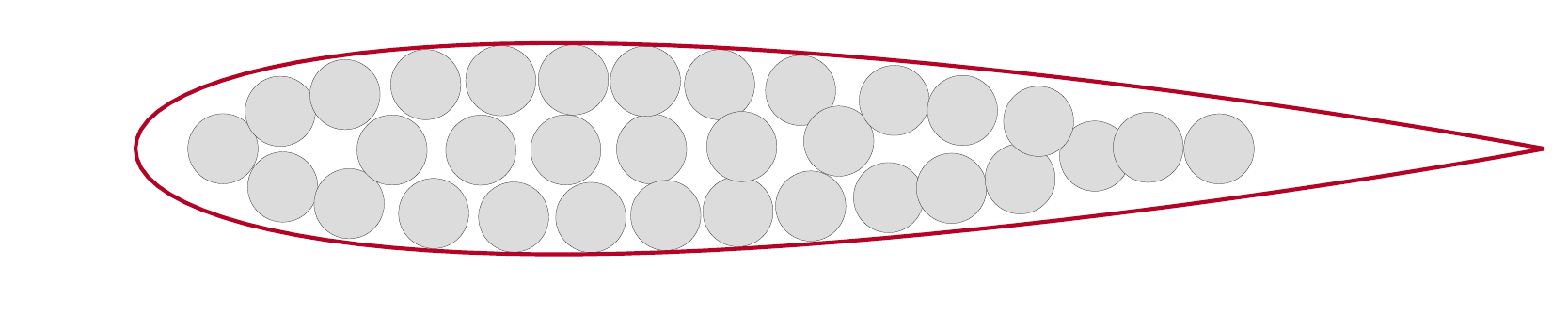}\\
    \includegraphics[width=0.95\textwidth]{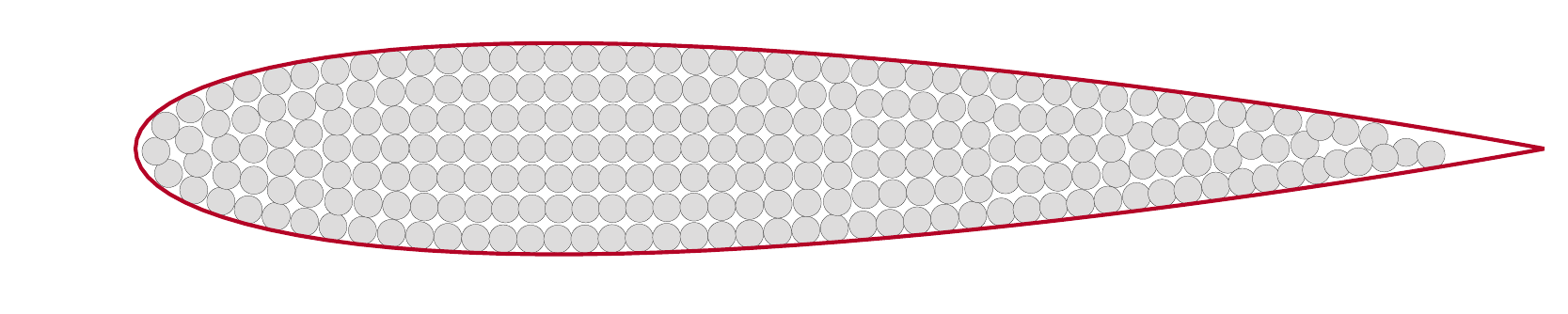}\\
    \includegraphics[width=0.95\textwidth]{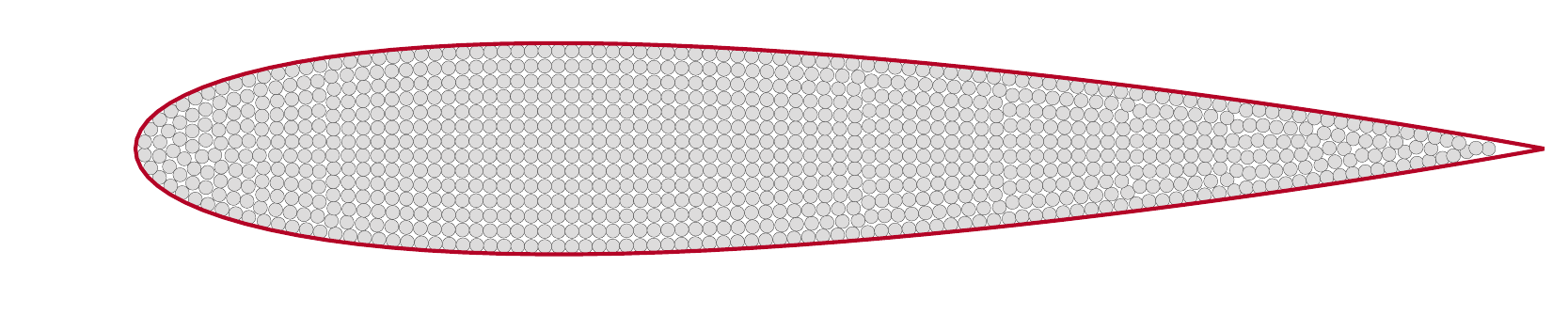}
    \caption{Packed with $h=0.8 \Delta x$.}
  \end{subfigure}
  \caption{Particle arrangement for the NACA0015 airfoil at different resolutions.
    From top to bottom: $\Delta x = 0.05$, $\Delta x = 0.02$, $\Delta x = 0.01$.}%
  \label{fig:airfoil}
\end{figure}
The particles are again packed using a quintic kernel with $h=0.8 \Delta x$.
In Fig.~\ref{fig:airfoil}, the initial configuration for different resolutions are shown on the left,
while the corresponding packed results are displayed on the right.
It can be observed that surface area with high curvature are increasingly better resolved with higher resolution.
This effect is most evident at the tip of the airfoil, where a sharp geometry causes the minimum distance, $\phi = \Delta x /2 $,
to shift further into the interior of the geometry at lower resolutions.
This behavior is not observed when the particles are placed directly on the surface, e.g. when the packing is used for TLSPH.
One potential remedy is to employ the scale separation method proposed by Luo et al.~\cite{Luo2016},
where the SDF is reconstructed according to the particle spacing,
meaning that the tip would be rounded off based on the resolution.
For example, Zhu et al.~\cite{Zhu2021} applied this method to reduce the kinetic energy during packing.
However, in the next section, we will demonstrate that our relaxation process remains stable even without geometry reconstruction.

\subsubsection{Airfoil 30P--30N}
\label{sec:airfoil_30P}
\begin{figure}[!h]
  \centering
  \begin{subfigure}[b]{0.4\textwidth}
    \centering
    \includegraphics[width=0.9\textwidth]{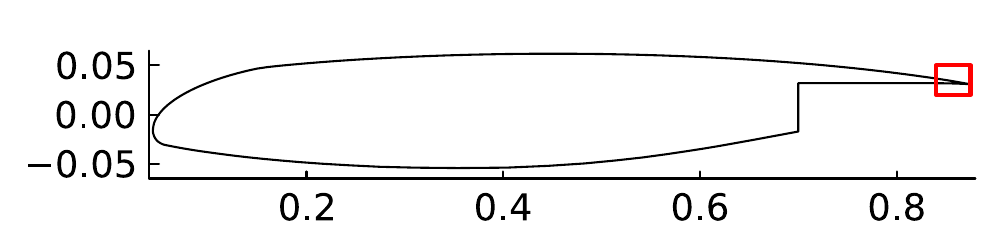}
    \caption{Airfoil 30P--30N with the area of interest highlighted in red.}%
    \label{fig:airfoil_30P_area_of_interest}
  \end{subfigure}
  \hfill
  \begin{subfigure}[b]{0.55\textwidth}
    \centering
    \includegraphics[width=\textwidth]{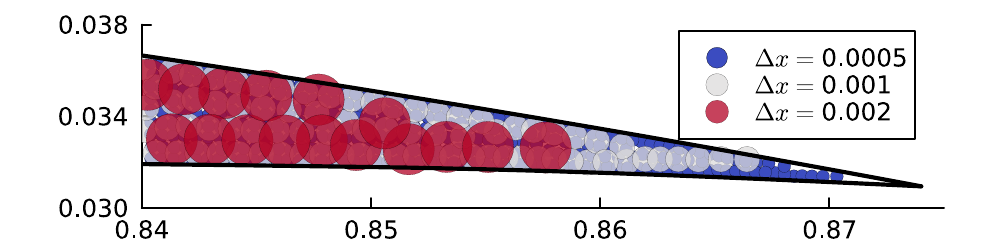}
    \caption{Results for different resolutions.}%
    \label{fig:airfoil_30P_different_resolutions}
  \end{subfigure}
  \caption{Left: Airfoil 30P-30N geometry.  Right: Zoomed area of interest.}%
\end{figure}

In the following, we highlight three key aspects of this work.
First, we demonstrate that increasing the resolution allows the geometry to be captured with better accuracy.
Second, we investigate the robustness of the final particle configuration with respect to different initial configurations.
Third, we show that a rapid convergence to a steady state is reached even for geometries with sharp corners.
To this end, we consider the Airfoil 30P–30N, which features a very sharp tip (Fig.~\ref{fig:airfoil_30P_area_of_interest}).
\begin{figure}[!h]
  \centering
  \begin{subfigure}[b]{0.48\textwidth}
    \centering
    \includegraphics[width=0.85\textwidth]{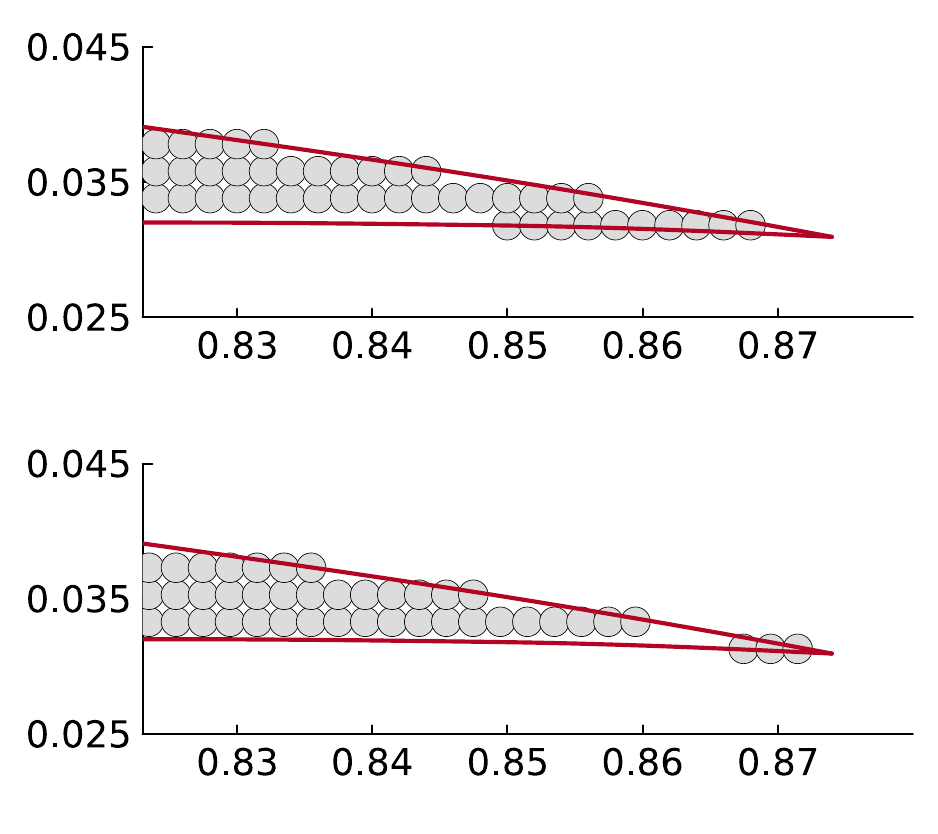}
    \caption{Initial configuration.}%
    \label{fig:airfoil_30P_a}
  \end{subfigure}
  \hfill
  \begin{subfigure}[b]{0.48\textwidth}
    \centering
    \includegraphics[width=0.85\textwidth]{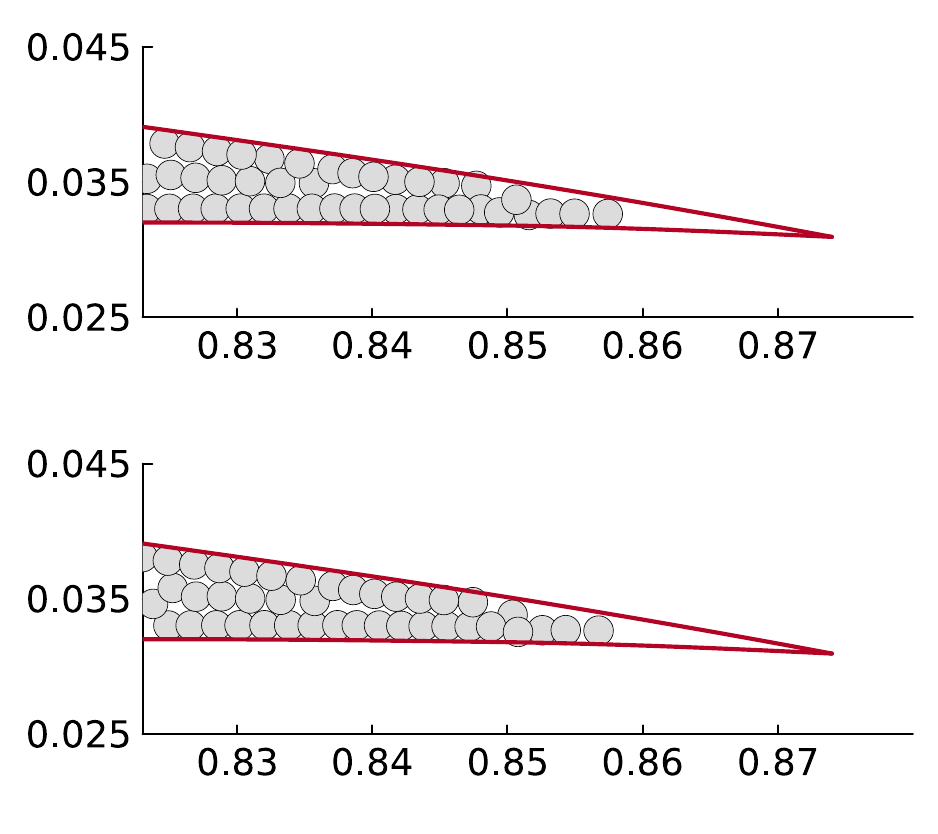}
    \caption{Packed configuration.}%
    \label{fig:airfoil_30P_b}
  \end{subfigure}\caption{
    Comparison of two different initial particle configurations at the tip of the Airfoil 30P-30N and their packed results.
    The initial configurations (left) differ by a slight spatial shift, yet the final relaxed states (right) converge to nearly identical particle distribution.
    }%
  \label{fig:airfoil_30P}
\end{figure}

Fig.~\ref{fig:airfoil_30P_different_resolutions} illustrates the tip of the 30P-30N airfoil for different resolutions.
It is clear that with increasing resolution, the geometry is represented more accurately.
This is enabled by two key factors:
first, the inside-outside segmentation based on the winding number resolves the geometry exactly,
allowing the initial particle configuration to converge towards the exact geometry as the resolution increases;
second, the SDF is constructed and evaluated according to the chosen resolution determined by the particle spacing
of the initial configuration.

Fig.~\ref{fig:airfoil_30P_a} shows two different initial particle configurations at the tip of the 30P--30N airfoil, where the edges of the geometry are highlighted in red.
Despite a slight spatial offset in the initial configuration, including detached particles,
both configurations evolve toward nearly identical final configurations, as shown in Fig.~\ref{fig:airfoil_30P_b}.
No particles detach from the main body.
This robustness is primarily attributed to the use of dynamic boundary particles.
They ensure that even initially detached particles remain influenced by an isotropic particle distribution,
stabilizing their motion and naturally drive them back toward the main body.
% In addition, the Shepard interpolation of the SDF further contributes to the robustness.
% The normal directions used for the bounding method (Sec.~\ref{sec:bounding_method}) are interpolated (Eq.~\ref{eq:interpolate n}) from the calculated normals of the discrete geometry surface,
% rather than being derived via numerical differentiation of a level-set gradient (Eq. ~\ref{eq:level_set_normal}).

\begin{figure}[!h]
  \centering
  \begin{subfigure}[b]{0.48\textwidth}
    \centering
    \includegraphics[width=0.95\textwidth]{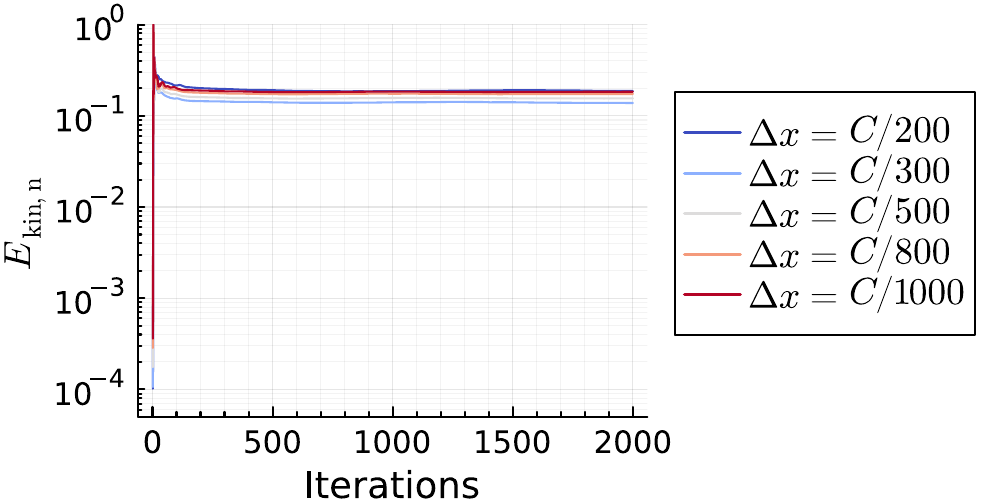}
    \caption{No boundary particles.}%
  \end{subfigure}
  \hfill
  \begin{subfigure}[b]{0.48\textwidth}
    \centering
    \includegraphics[width=0.95\textwidth]{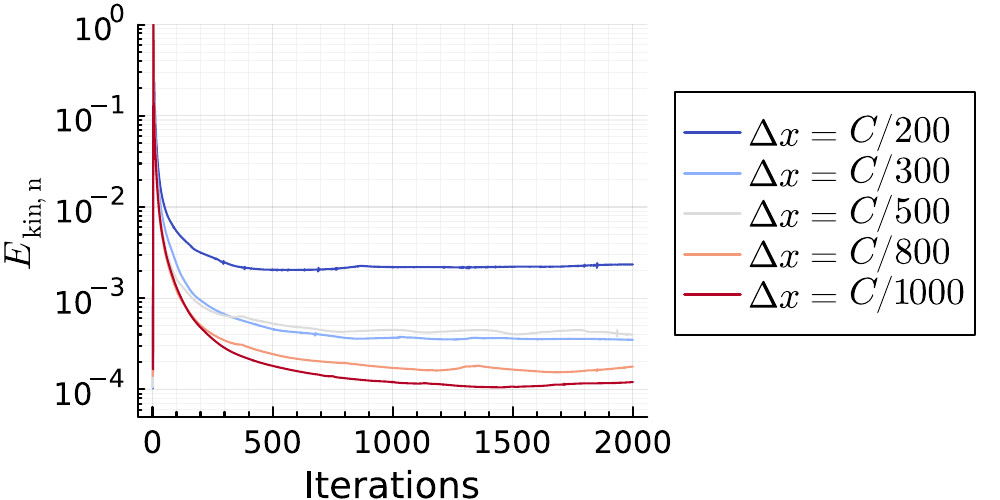}
    \caption{Dynamic boundary particles.}%
  \end{subfigure}%
  \caption{Normalized kinetic energy $E_{\mathrm{kin}, \mathrm{n}}$ during packing for the airfoil 30P--30N at different resolutions, where $C$ is the airfoil's chord length}%
  \label{fig:kinetic_energy_airfoil_30P_30N}
\end{figure}

To further demonstrate that the relaxation process remains stable,
Fig.~\ref{fig:kinetic_energy_airfoil_30P_30N} shows the normalized kinetic energy for the airfoil 30P--30N over each iteration.
It can be observed that, without boundary packing, the kinetic energy does not sufficiently decay,
although it reaches a steady state for each resolution.
With boundary packing, the kinetic energy rapidly drops by at least two orders of magnitude and stabilizes after 1000 iterations.
  Without adding algorithmic complexity, our dynamic boundary packing is able to achieve lower residual kinetic energy compared to the results in~\cite{Yu2023a}.

Previous studies have reported challenges with unresolved sharp features and fine details in the geometry~\cite{Yu2023a}.
During packing, unresolved singular points in the initial particle configuration lead to particles detaching from the main body and prevent the system from reaching a steady state.
To mitigate this, geometry reconstruction techniques, such as the scale separation method~\cite{Luo2016} cleaning non-resolved structures, have been employed,
resulting in improved packing behavior.
We demonstrated that without geometry reconstruction, we can achieve a rapid steady state, obtain robust results, and resolve the geometry with arbitrary accuracy.

\subsection{Particle packing in 3D}
\label{sec:particle_packing_in_3D}
Next, we focus on the deviation from the reference density and the convergence behavior of the proposed method
by considering two 3D geometries: a simple ellipsoid and a more complex biomechanical geometry, the aorta.
The deviation from the reference density is quantified using the $L_2$ and $L_\infty$ error norms (see Eq.~\eqref{eq:L_2} and Eq.~\eqref{eq:L_infty}).
The initial particle configuration is perfectly uniform and would yield only minimal error associated with the numerical density evaluation.
During the packing process, the particle distribution initially degrades due to the redistribution of particles toward the geometry surface,
leading to a temporary increase in density error.
As the process continues, the distribution gradually improves again, while preserving the geometry's surface representation.
The following investigations reflects the method's ability to balance uniform particle distribution with precise interface representation.
Before delving deeper into this, however, we first examine the impact of time integration on the packing process.

\subsubsection{Ellipsoid}
\label{sec:ellipsoid}
An adaptive time integration scheme based on the explicit Runge--Kutta method \texttt{RDPK3SpFSAL35}~\cite{Ranocha2021} is employed,
using a fixed relative tolerance of $10^{-3}$ while varying the absolute tolerance.
A higher absolute tolerance is expected to result in fewer iterations at the cost of a higher error,
whereas a lower absolute tolerance should yield improved accuracy but require more iterations.
\begin{figure}[!h]
  \centering
    \includegraphics[width=0.5\textwidth]{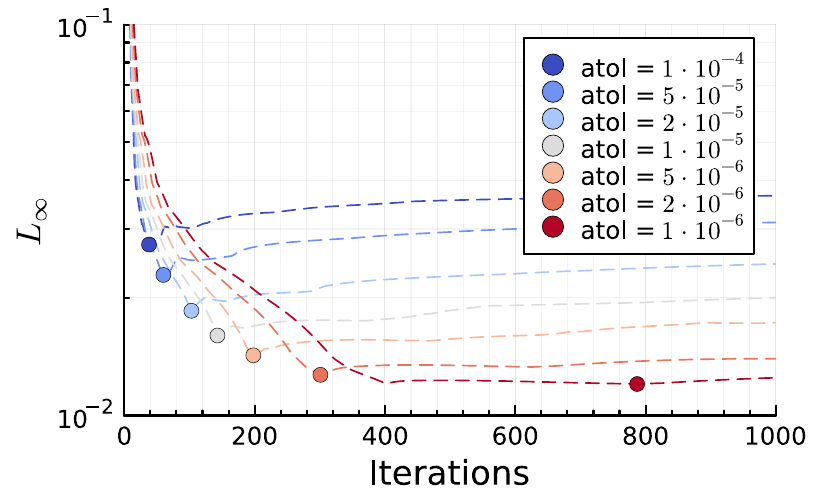}
    \caption{
      $L_\infty$ error in the density for the ellipsoid with different absolute tolerances (\texttt{atol}) for the adaptive time integration scheme \texttt{RDPK3SpFSAL35}~\cite{Ranocha2021}.
      Dashed lines represent the $L_\infty$ error per iteration, while circular markers indicate the minimum error achieved within each run.
      The initial particle spacing is $\Delta x = 0.1$.}
    \label{fig:L_inf_ellipsoid}
\end{figure}%
This implies that the overall error is bounded by the tolerances imposed by the time integration scheme.
Our objective is to evaluate whether a sufficiently low error can already be achieved with a relatively high absolute tolerance,
thus reducing computational cost while maintaining accuracy within an acceptable range.
For this purpose, we consider a simple ellipsoid geometry based on the dataset provided by~\cite{Negi2021}.
The packing process is performed for a total of $5 \cdot 10^{3}$ iterations across different values of the absolute tolerance.

Figure~\ref{fig:L_inf_ellipsoid} shows the evolution of the $L_\infty$ error (Eq.~\eqref{eq:L_infty}) over the number of iterations for different absolute tolerances (\texttt{atol}).
Dashed lines indicate the iteration-wise $L_\infty$ error progression for each tolerance level, while circular markers highlight the smallest $L_\infty$ error reached within each run.
It can be observed that, even with a higher absolute tolerance, an error on the same order of the lower tolerance is temporarily reached after fewer than 100 iterations.
The increase in error observed after reaching the minimum can be explained by the effect of the adaptive time integration scheme:
as the packing approaches a steady state, the time steps increase, which can temporarily lead to a slightly increased error.
In comparison to the results reported in~\cite{Negi2021}, our approach achieves a similar error level with significantly fewer iterations.

\subsubsection{Aorta geometry}
\label{sec:result_aorta}
Next, we assess an aorta as an example for a complex biomechanical geometry.
The geometry is obtained from the extracted data sets of an MRI scan reported by Mirzae et al.~\cite{Mirzaee2016}.
To provide context, we first present qualitative results that illustrate the overall setup.
Subsequently, we focus on the final particle configurations obtained after the packing process and demonstrate how increasing resolution improves the accuracy of the geometry representation.
In the following, the characteristic length $L_\text{A}$ is defined as the approximate length of the longest edge of the bounding box enclosing the aortic geometry.

\begin{figure}[!h]
  \centering
  \begin{subfigure}[b]{0.3\textwidth}
    \centering
    \includegraphics[width=\textwidth]{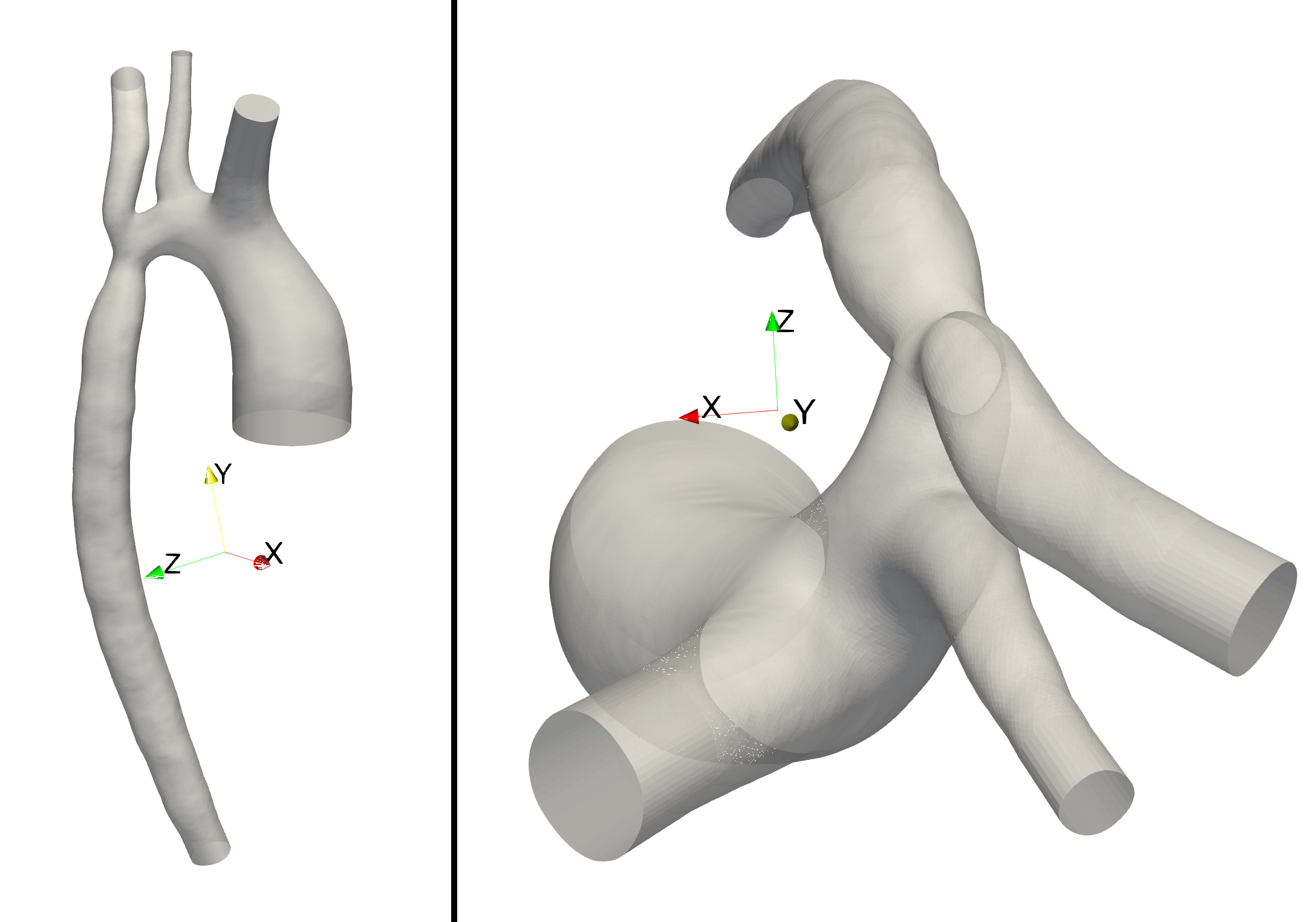}
    \caption{Aorta geometry.}%
    \label{fig:aorta_geometry}
  \end{subfigure}\\
  \begin{subfigure}[b]{0.45\textwidth}
    \centering
    \includegraphics[width=0.95\textwidth]{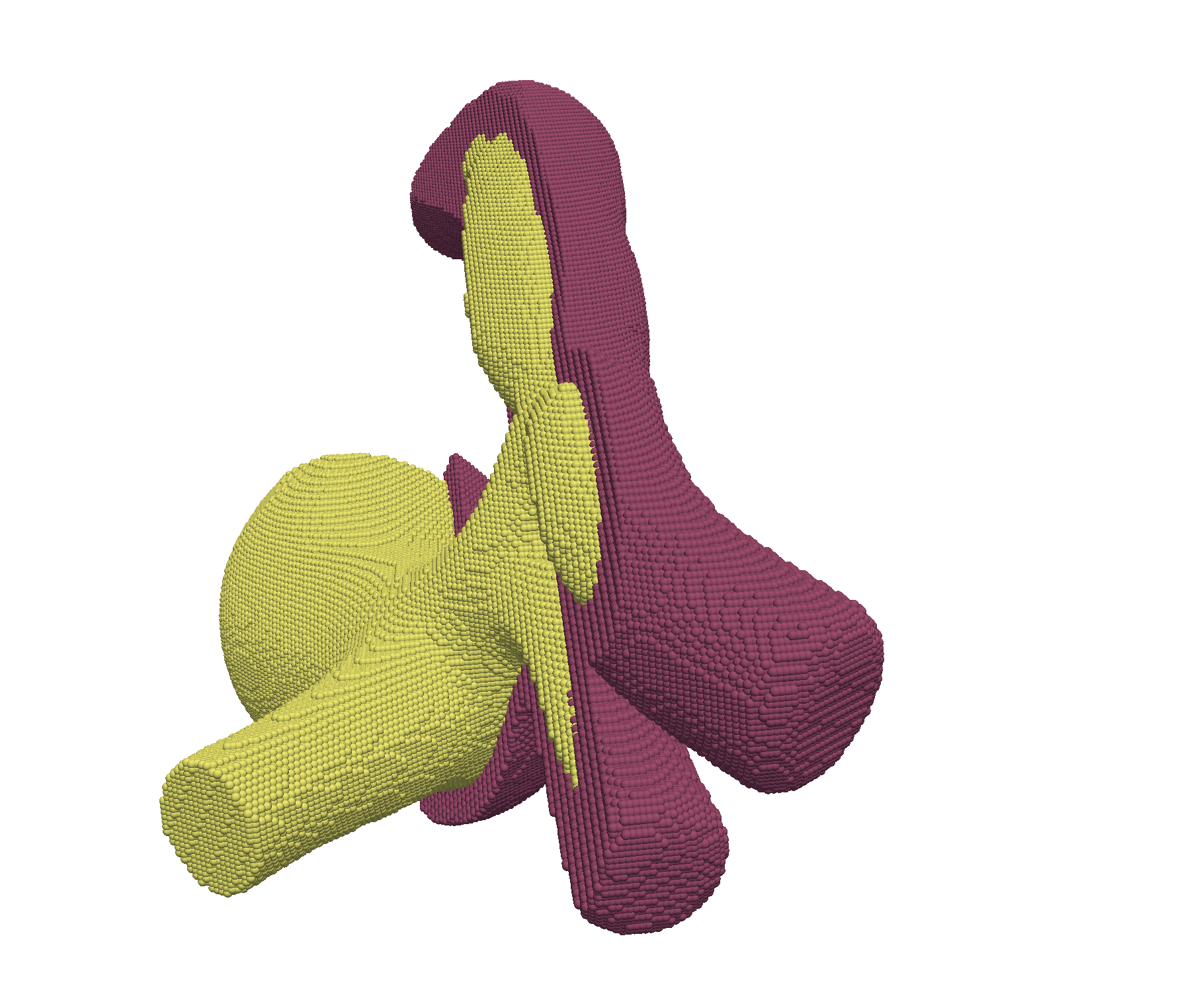}
    \caption{Initial configuration.}%
    \label{fig:aorta_sampled_a}
  \end{subfigure}
  \hfill
  \begin{subfigure}[b]{0.45\textwidth}
    \centering
    \includegraphics[width=0.95\textwidth]{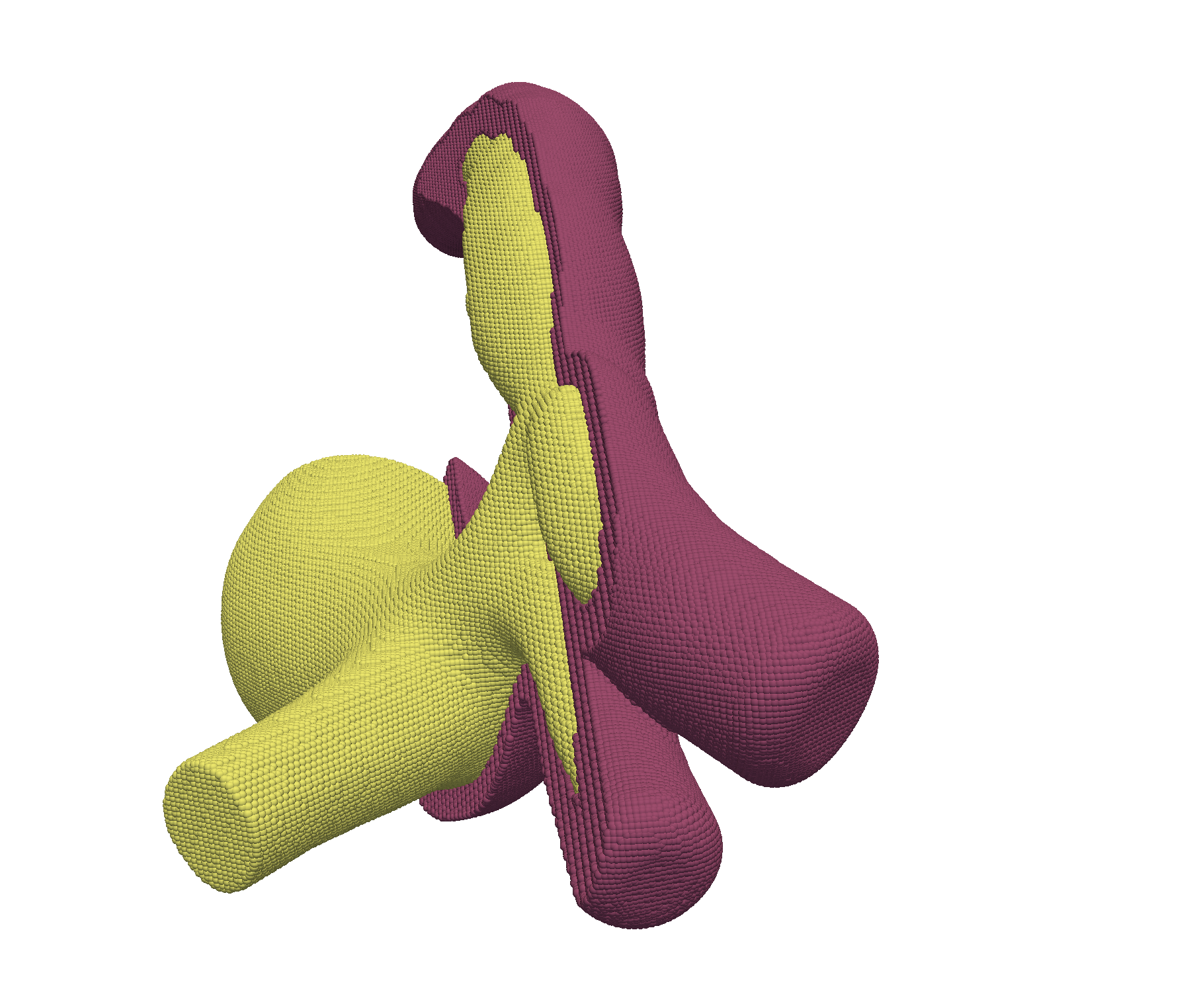}
    \caption{Packed configuration.}%
    \label{fig:aorta_sampled_b}
  \end{subfigure}%
  \caption{Visualization of the particle packing process for the aorta geometry.
    (a) The aorta geometry extracted from MRI data~\cite{Mirzaee2016}, rotated to provide a clearer view of the relevant features.
    (b) Initial particle configuration after segmentation, where yellow particles represent interior particles and purple particles represent boundary particles.
    (c) Final packed configuration of both the inner and boundary particles. The initial particle spacing  is $\Delta x = L_\text{A} / 400$ and the boundary thickness is $8\Delta x$.}%
  \label{fig:aorta_sampled}
\end{figure}%
Fig.~\ref{fig:aorta_geometry} shows the aorta geometry extracted from MRI data~\cite{Mirzaee2016}, rotated to provide a clearer view of the relevant features.
In Fig.~\ref{fig:aorta_sampled_a}, the initial particle configuration is visualized:
yellow particles represent the initial configuration of the interior particles,
while purple particles represent the initial configuration of the boundary particles.
Here, the boundary particles are generated with a thickness of $8\Delta x$.
For visualization purposes, the boundary particles have been clipped, making the interior particles partially visible.
Fig.~\ref{fig:aorta_sampled_b} illustrates the configuration after packing, revealing a smooth, body-fitted particle arrangement.
Since the boundary particles are dynamically packed, they smoothly adapt to the geometry.

\begin{figure}[!h]
  \centering
  \begin{subfigure}[b]{0.48\textwidth}
    \centering
    \includegraphics[width=0.85\textwidth]{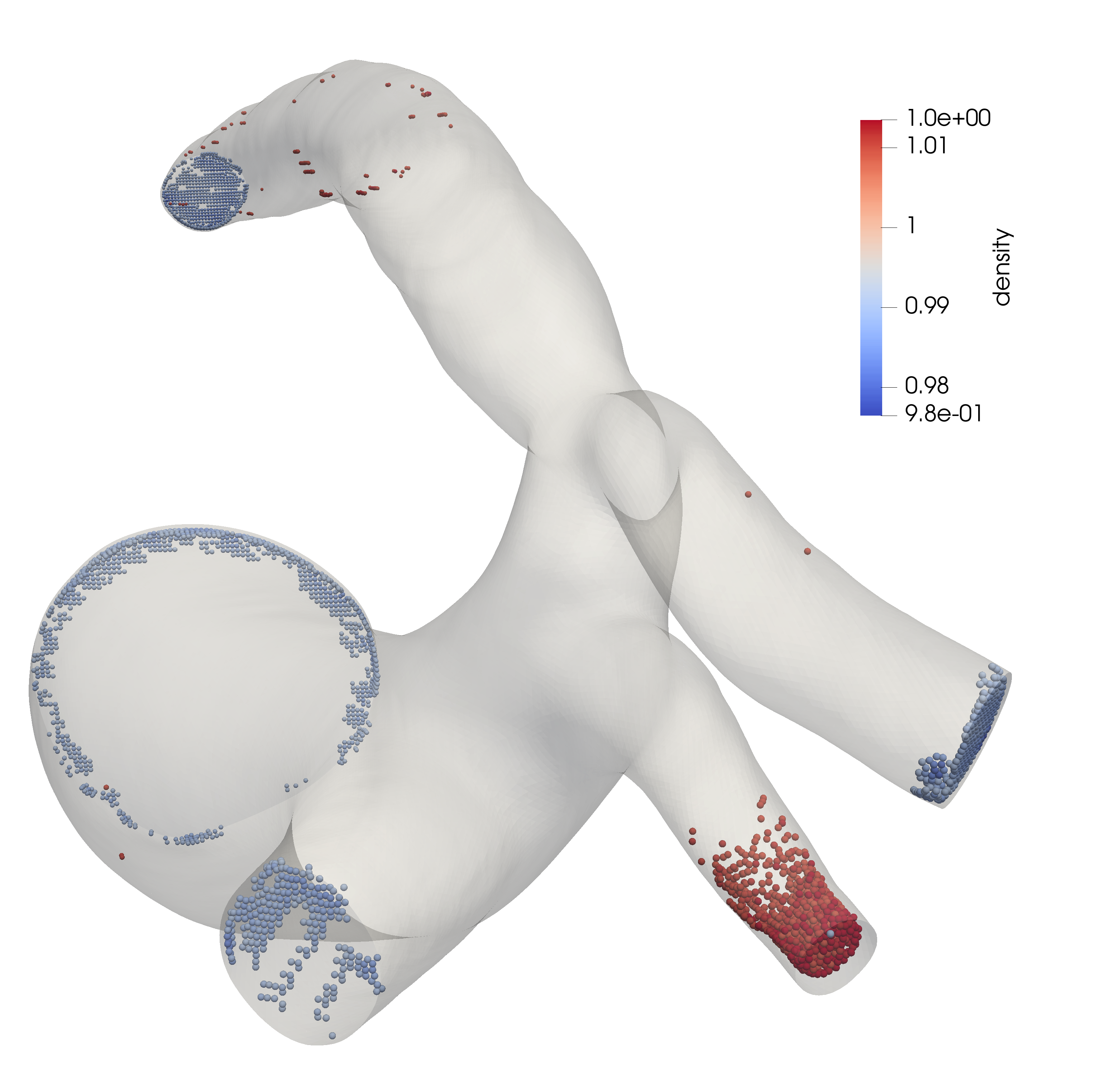}
    \caption{$\Delta x / L_\text{A} = 1/ 400$}%
    \label{fig:density_lower_a}
  \end{subfigure}
  \hfill
  \begin{subfigure}[b]{0.48\textwidth}
    \centering
    \includegraphics[width=0.85\textwidth]{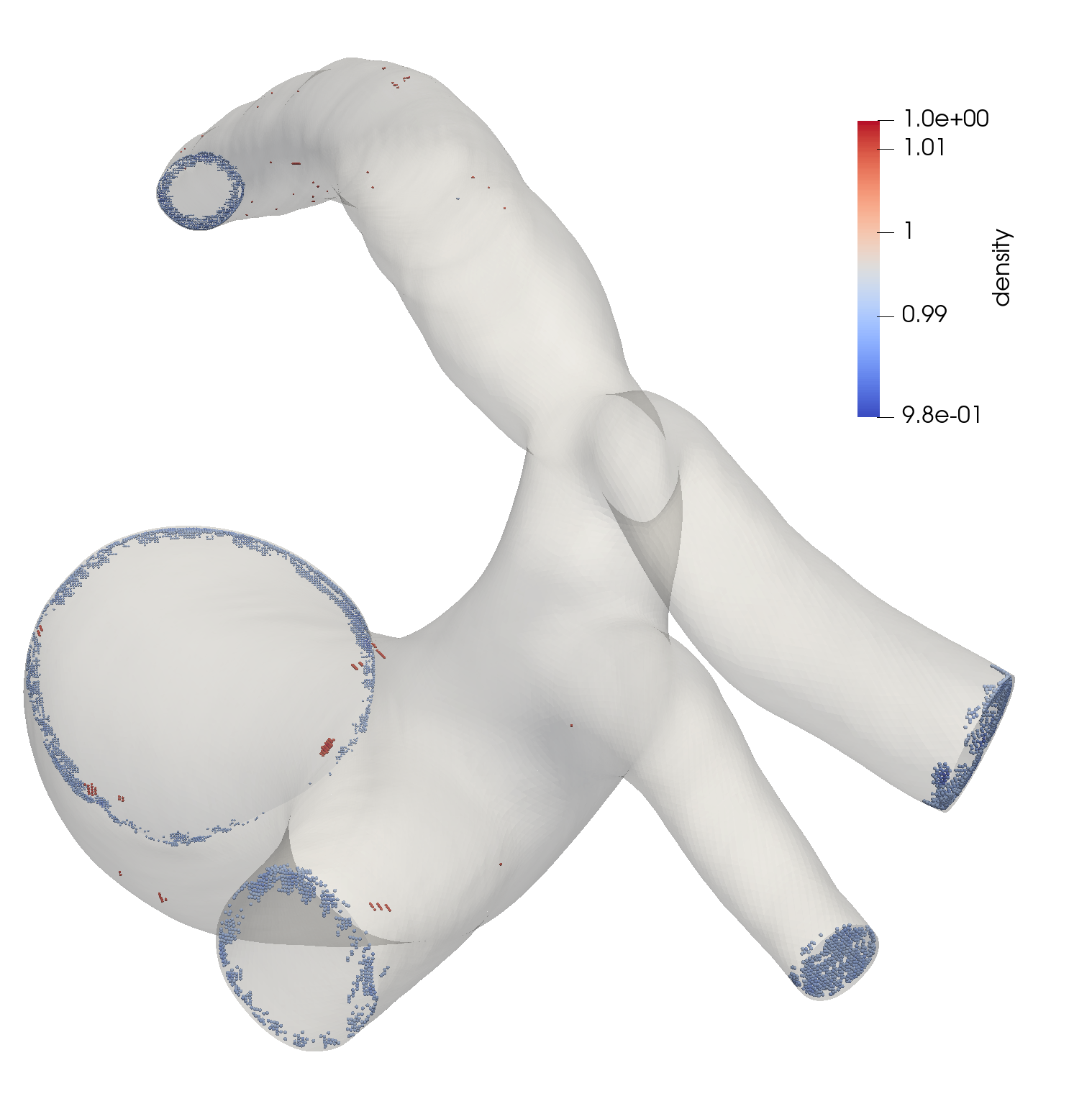}
    \caption{$\Delta x / L_\text{A} = 1 / 800$}%
    \label{fig:density_lower_b}
  \end{subfigure}%
  \caption{Particle packing for an aorta geometry at different resolutions. Only particles with a density deviation more than 1\% from the reference density are shown.
   The fraction of deviating particles decreases from $0.44\,\%$ in (a) to $0.096\,\%$ in (b).}%
  \label{fig:density_lower}
\end{figure}

Next, we investigate the density deviation in the final configuration at different resolutions.
Fig.~\ref{fig:density_lower_a} presents results for the aorta geometry at $\Delta x/ L_\text{A} = 1/ 400$ and Fig.~\ref{fig:density_lower_b} at $\Delta x/ L_\text{A} = 1/ 800$.
Only particles whose density values deviate from the reference density by more than 1\% are displayed.
In Fig.~\ref{fig:density_lower_a}, approximately $0.44\,\%$ of the particles deviate from the reference density, whereas in Fig.~\ref{fig:density_lower_b}, this fraction decreases to $0.096\,\%$.
The observed deviations are mainly concentrated near regions with sharp edges,
where a coarse resolution results in non-resolved geometric features and thus larger distances to the boundary particles.
As the resolution increases, these deviations are reduced, indicating an improved representation of the geometry interface and, thus, a more uniform particle distribution.
In this case, such deviations are less critical, since inflow and outflow buffer zones are created at these surfaces.

To quantify the convergence behavior of the particle distribution with increasing resolution,
we compute the experimental order of convergence (EOC), which provides a measure of how the error decreases when the particle spacing $\Delta x$ is refined, and is defined as
\begin{equation}
  \text{EOC} = \frac{\log\left(\frac{L_{2, i}}{L_{2, i+1}}\right)}{\log\left(\frac{\Delta x_i}{\Delta x_{i+1}}\right)},
  \label{eq:EOC}
\end{equation}
where $L_{2,i}$ and $L_{2,i+1}$ denote the $L_2$ errors corresponding to two successive resolutions $\Delta x_i$ and $\Delta x_{i+1}$, respectively.
Table~\ref{tab:l2_error_particle_spacings} summarizes the $L_2$ errors and the corresponding EOC values for various resolutions.
As the resolution increases, the $L_2$ error decreases, indicating improved accuracy of the particle distribution and better representation of the geometry interface.
The EOC values further confirm this trend, with convergence rates that vary depending on the resolution level, but generally indicate a consistent reduction in error with refinement.
\begin{table}[!t]
  \centering
   \caption{Convergence behavior for the aorta geometry.}
    \begin{tabular}{ccccc}
      \toprule
      $\Delta x/ L_\text{A}$ & $L_2$ error & EOC($L_2$) & Fraction of particles   \\
                             &             &            & with $>1\%$ density deviation  \\
      \midrule
        $ 1/50$ & $7.681 \cdot 10^{-2}$  & --- &  $99.92\,\%$ \\
        $ 1/100$ & $1.161 \cdot 10^{-2}$& 2.726  &  $31.11\,\%$ \\
        $ 1/200$ & $4.797 \cdot 10^{-3}$ & 1.275& $5.62\,\%$  \\
        $ 1/400$ & $1.857 \cdot 10^{-3}$& 1.369 & $0.44\,\%$  \\
        $ 1/800$ & $9.670 \cdot 10^{-4}$ & 0.942& $0.096\,\%$  \\
        % $1.25 \cdot 10^{-3}$ & $5.820 \cdot 10^{-4}$ & 0.732 \\
    \bottomrule
  \end{tabular}
  \label{tab:l2_error_particle_spacings}
\end{table}

\subsubsection{Multi-body packing}
Packing multiple geometries simultaneously is especially useful for applications such as airfoils in channel flows.
Several approaches are available:
pack all geometries simultaneously,
sequentially pack already processed particle configurations with additional geometries,
or combine both approaches and then extract the particles corresponding to a specific geometry using segmentation techniques (e.g. the winding number method).
This last approach is particularly beneficial for creating buffer zones at open boundaries,
ensuring a smooth transition from the buffer region to the interior domain.

\begin{figure}[!h]
  \centering
  \begin{subfigure}[b]{0.48\textwidth}
    \centering
    \includegraphics[width=0.75\textwidth]{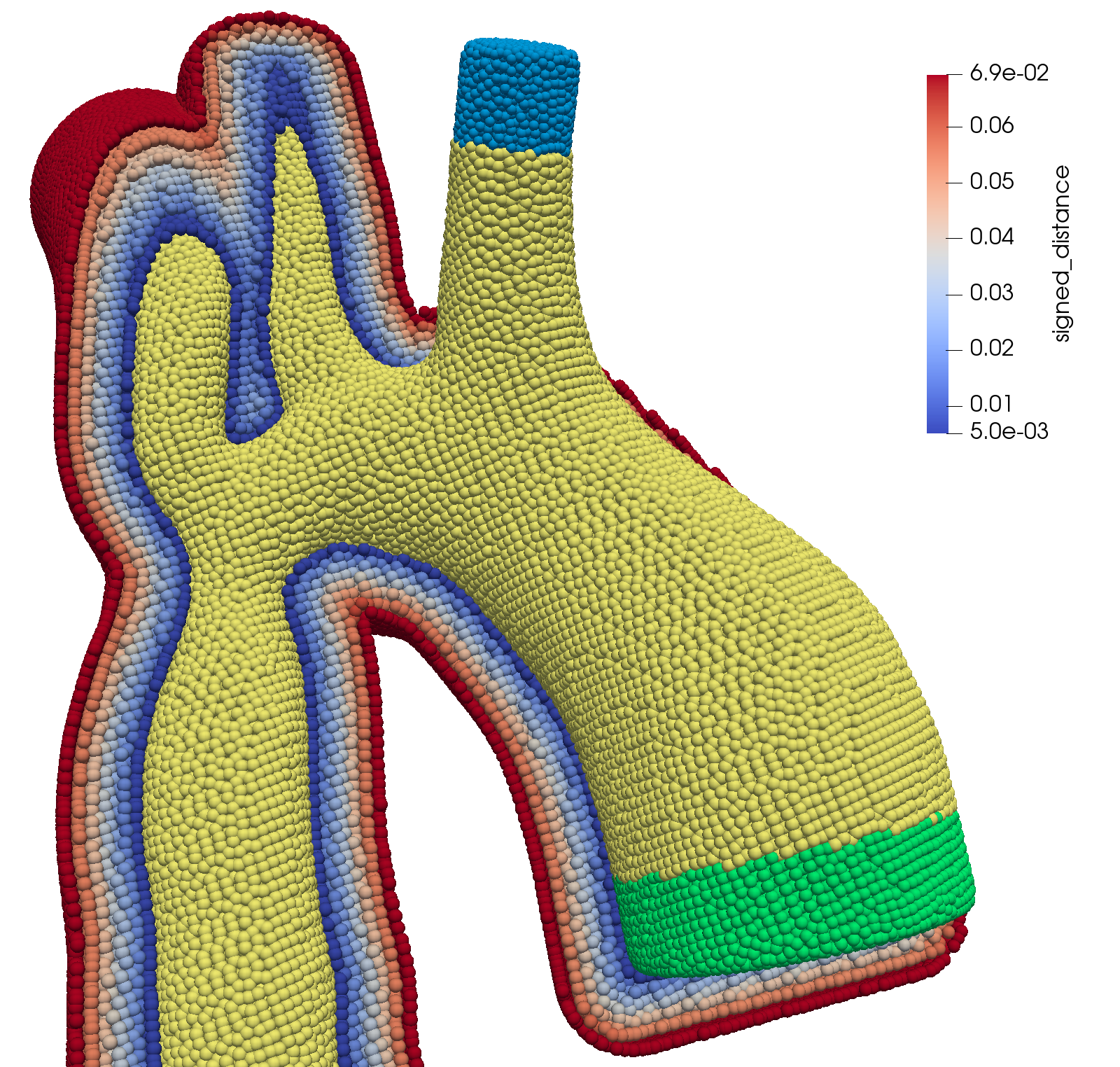}
    \caption{$\Delta x / L_\text{A} = 1 / 200$}%
  \end{subfigure}
  \hfill
  \begin{subfigure}[b]{0.48\textwidth}
    \centering
    \includegraphics[width=0.75\textwidth]{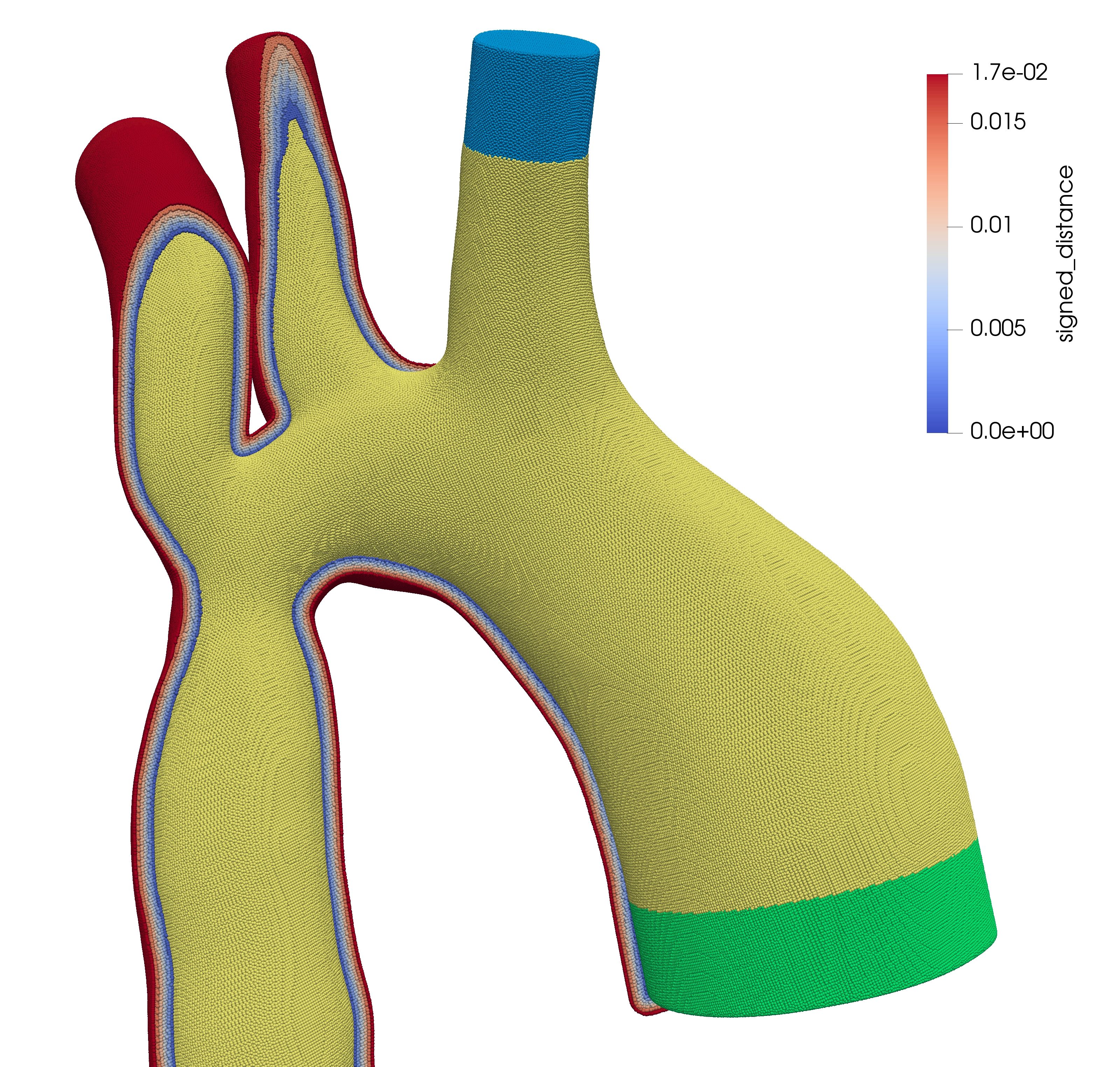}
    \caption{$\Delta x / L_\text{A} = 1 / 800$}%
  \end{subfigure}%
  \caption{Packed particle distribution of the aorta geometry~\cite{Mirzaee2016} with different resolutions. The boundary particles are colored with the signed distance to the surface of the geometry.
  The yellow particles are the interior fluid particles and the blue and green particles are the buffer particles for the out- and inflow, respectively.}%
  \label{fig:multi_body}
\end{figure}%
Fig.~\ref{fig:multi_body} illustrates the packed configurations of the aorta geometry at two different resolutions.
The boundary particles are colored according to their signed distance to the surface of the geometry,
demonstrating that the boundary thickness automatically adjusts with resolution.
In this figure, the green particles represent the inflow zone and the blue particles represent the outflow zone.
The inflow and outflow geometries were previously extruded using Blender~\cite{blender}.
After packing, the particles in these extruded regions were separated using the winding number method
to facilitate the boundary conditions for simulating blood flow through the aorta.

\subsection{Simulation of an Airfoil}
\label{sec:simulation_airfoil}
\begin{figure}[!h]
    \centering
    \includegraphics[width=0.9\textwidth]{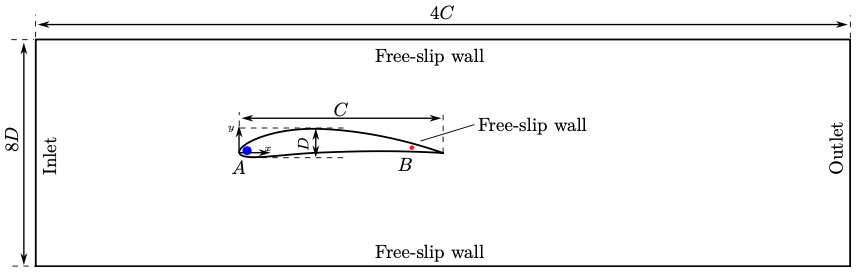}
  \caption{Simulation setup for the FSI study with a NACA6412 airfoil.
  The computational domain has a height of $8D$ and a length of $4C$, with free-slip conditions applied at the top and bottom walls as well as the airfoil surface.
  The airfoil is fixed at the circular area $A$, and the vertical displacement is measured at point $B$.}
    \label{fig:setup_simulation}
\end{figure}

To demonstrate the practical applicability of our preprocessing approach,
we conduct a fluid-structure interaction (FSI) simulation, in which a NACA6412 airfoil serves as an elastic structure immersed within a fluid domain.
This study aims to compare lattice-based and packed particle configurations.
We investigate the vertical displacement of the trailing edge as well as the von Mises strain distribution within the airfoil.
All walls, including the airfoil surface, are simulated with free-slip boundary conditions to investigate the surface smoothness resulting from the particle packing.
Fig.~\ref{fig:setup_simulation} presents a detailed overview of the simulation setup, closely following the configuration described in \cite{Yu2023a}.
All parameters are given in dimensionless form.
The NACA6412 airfoil is positioned such that its leading edge coincides with $(0.0,0.0)$.
Fixation is applied at a circular area $A$ centered at $(0.04,0.01)$ with a radius of $r=0.025$.
Displacement measurements are taken at the trailing edge, located at point $B=(0.85,0.03)$.
At the inlet, a constant inflow velocity of $v_{0}=2.0$ is imposed, while a pressure boundary condition $p_{\text{out}}=0.0$ is maintained at the outlet.
Characteristic geometric and material properties include a chord length of $C=1.0$ and a thickness of $D=0.13895$;
fluid and solid densities of $\rho^f=1.0$ and $\rho^s=10 \rho^f$, respectively;
a dimensionless Young's modulus defined as $E^* = E / (\rho^f v_0^2) = 1.4 \times 10^3$, with Poisson's ratio $\nu = 0.4$;
and flow conditions corresponding to a Reynolds number $Re = \rho^f v_0 C / \eta = 200$.

\begin{figure}[!h]
  \centering
    \begin{subfigure}[b]{0.4\textwidth}
    \centering
    \includegraphics[width=0.9\textwidth]{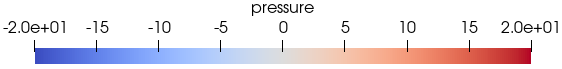}
  \end{subfigure}\\
  \begin{subfigure}[b]{0.48\textwidth}
    \centering
    \includegraphics[width=0.75\textwidth]{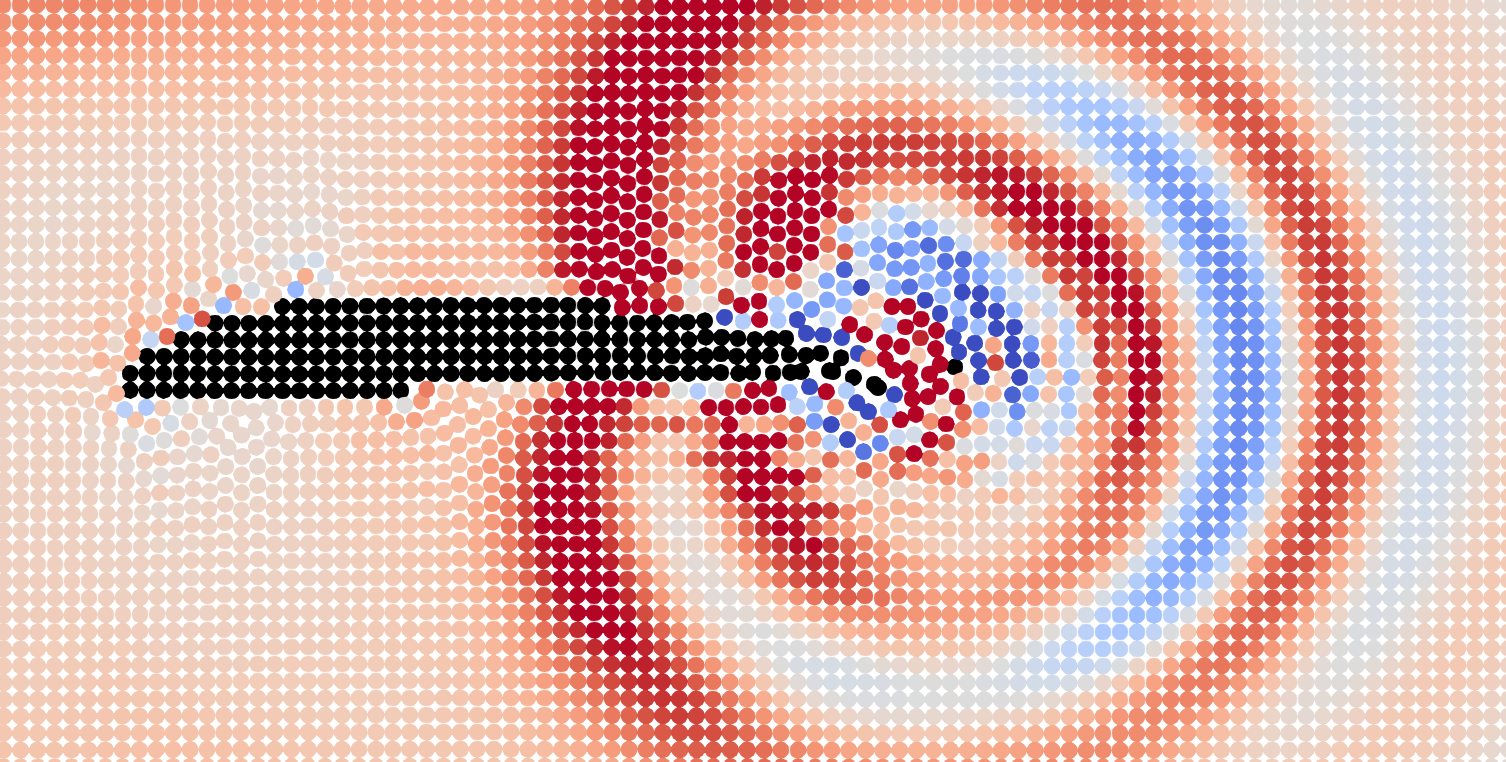}
    \caption{Lattice configuration $\Delta x = C /50$.}%
  \end{subfigure}
  \hfill
  \begin{subfigure}[b]{0.48\textwidth}
    \centering
    \includegraphics[width=0.75\textwidth]{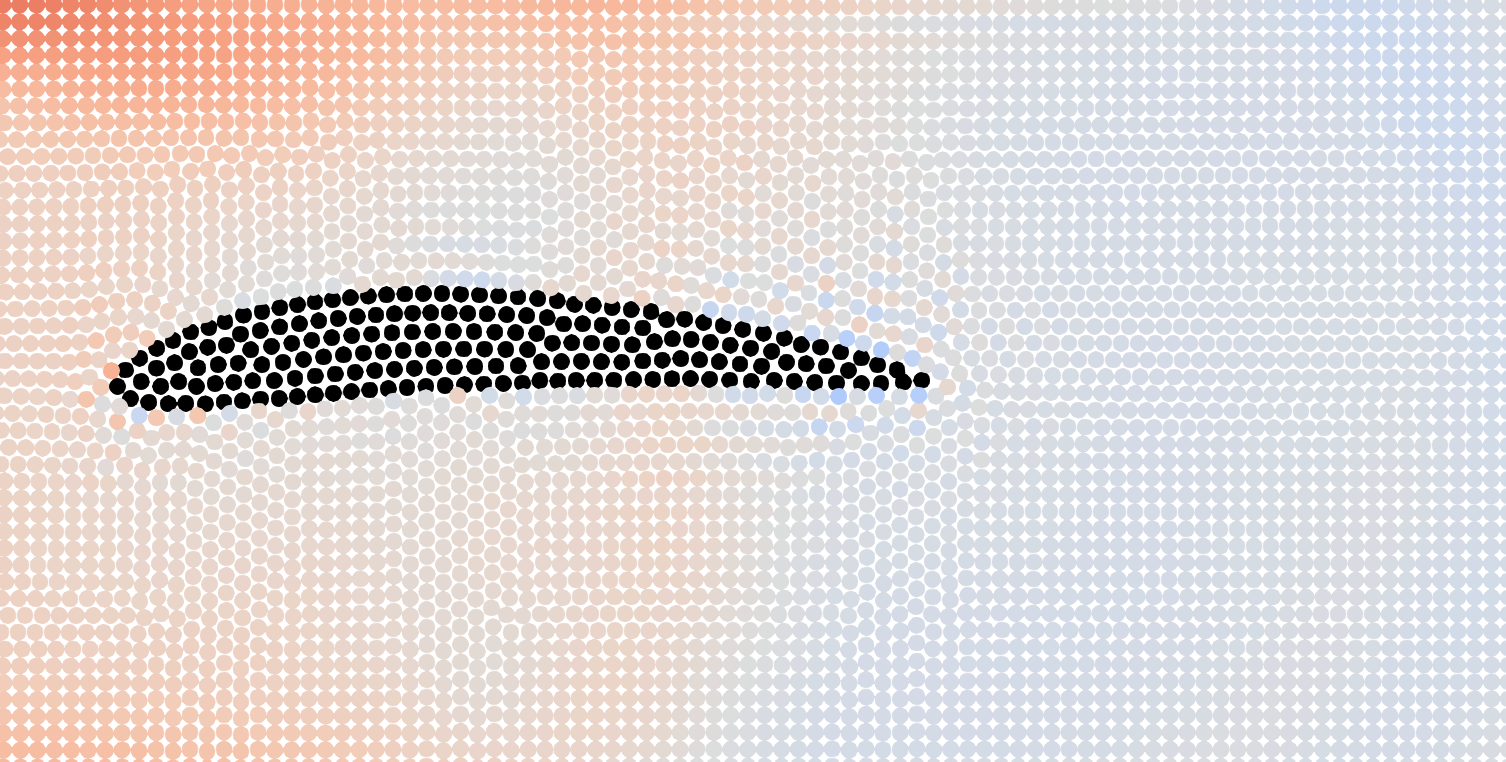}
    \caption{Packed configuration $\Delta x = C /50$.}%
  \end{subfigure}\\
  \caption{Comparison of pressure fields around a NACA6412 airfoil for lattice-based (left) and packed (right) particle configurations at $t = 0.02$.}%
  \label{fig:solid_fluid_1}
\end{figure}

\begin{figure}[!h]
  \centering
    \includegraphics[width=0.65\textwidth]{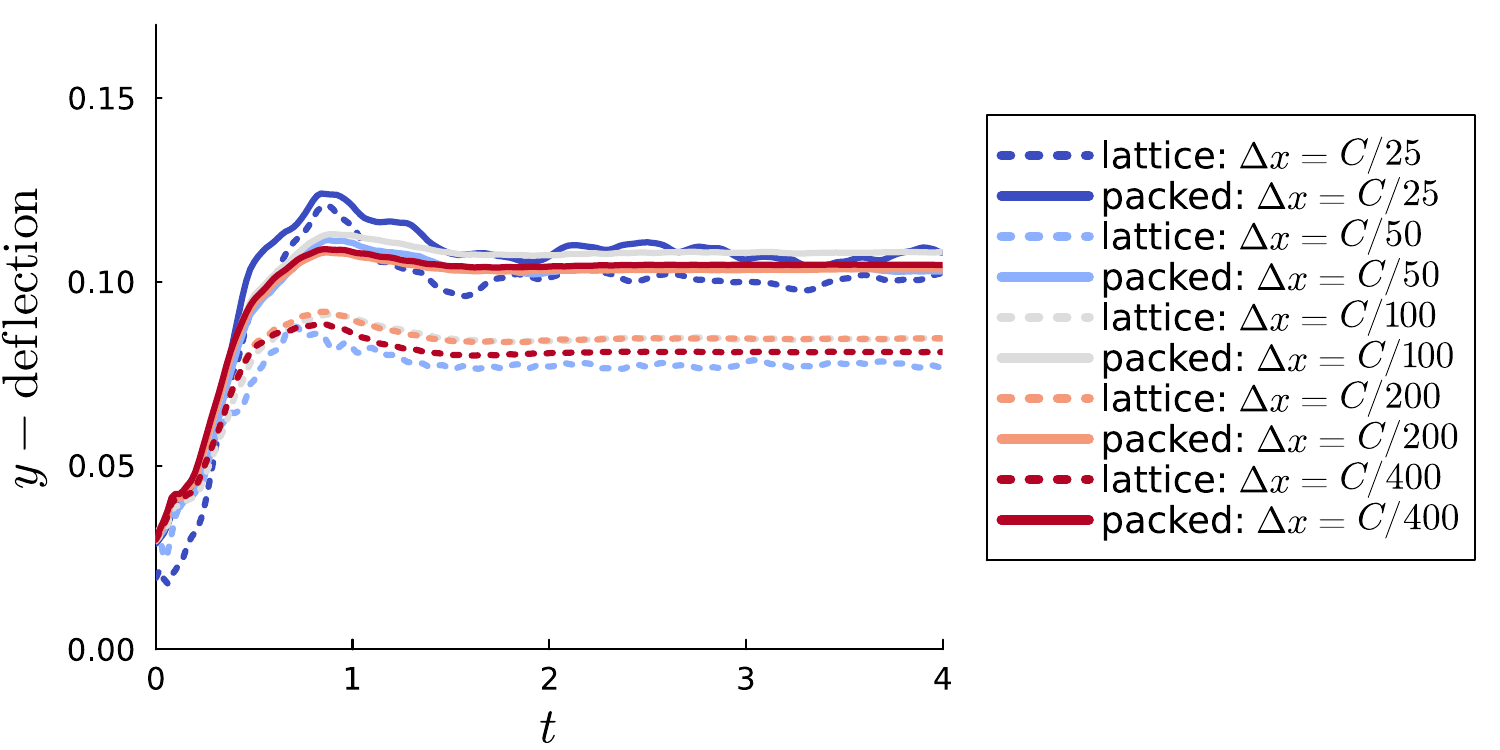}
    \caption{Vertical deflection of the trailing edge of the NACA6412 airfoil over time for different configurations.}
    \label{fig:y_deflection}
\end{figure}

\begin{figure}[!h]
  \centering
    \begin{subfigure}[b]{0.6\textwidth}
    \centering
    \includegraphics[width=0.9\textwidth]{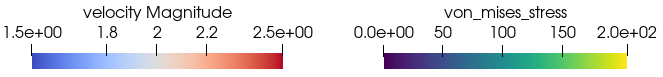}
  \end{subfigure}\\
  \begin{subfigure}[b]{0.48\textwidth}
    \centering
    \includegraphics[width=0.75\textwidth]{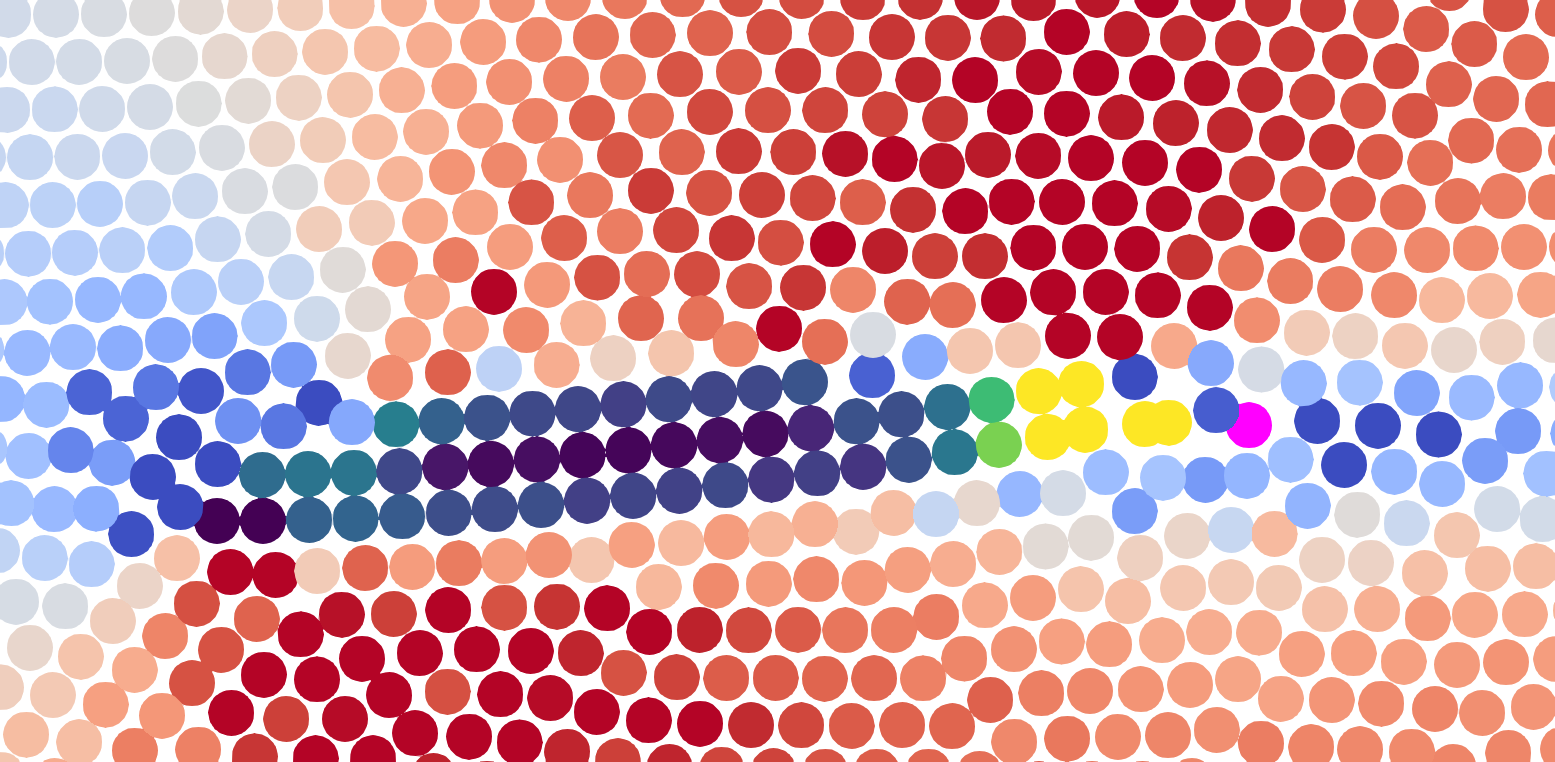}
    \caption{Lattice configuration $\Delta x = C /25$.}%
  \end{subfigure}
  \hfill
  \begin{subfigure}[b]{0.48\textwidth}
    \centering
    \includegraphics[width=0.75\textwidth]{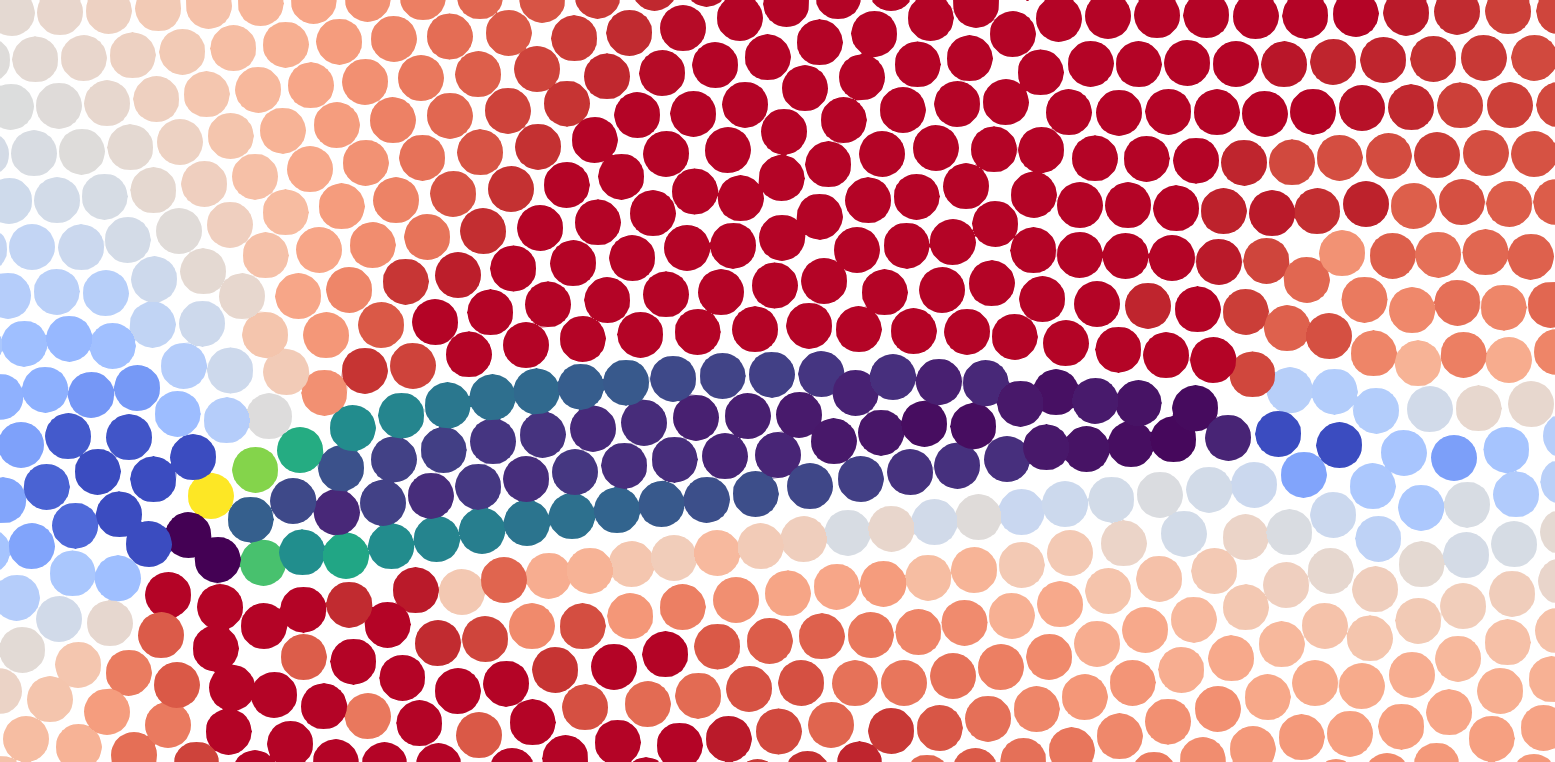}
    \caption{Packed configuration $\Delta x = C /25$.}%
  \end{subfigure}\\
  \begin{subfigure}[b]{0.48\textwidth}
    \centering
    \includegraphics[width=0.75\textwidth]{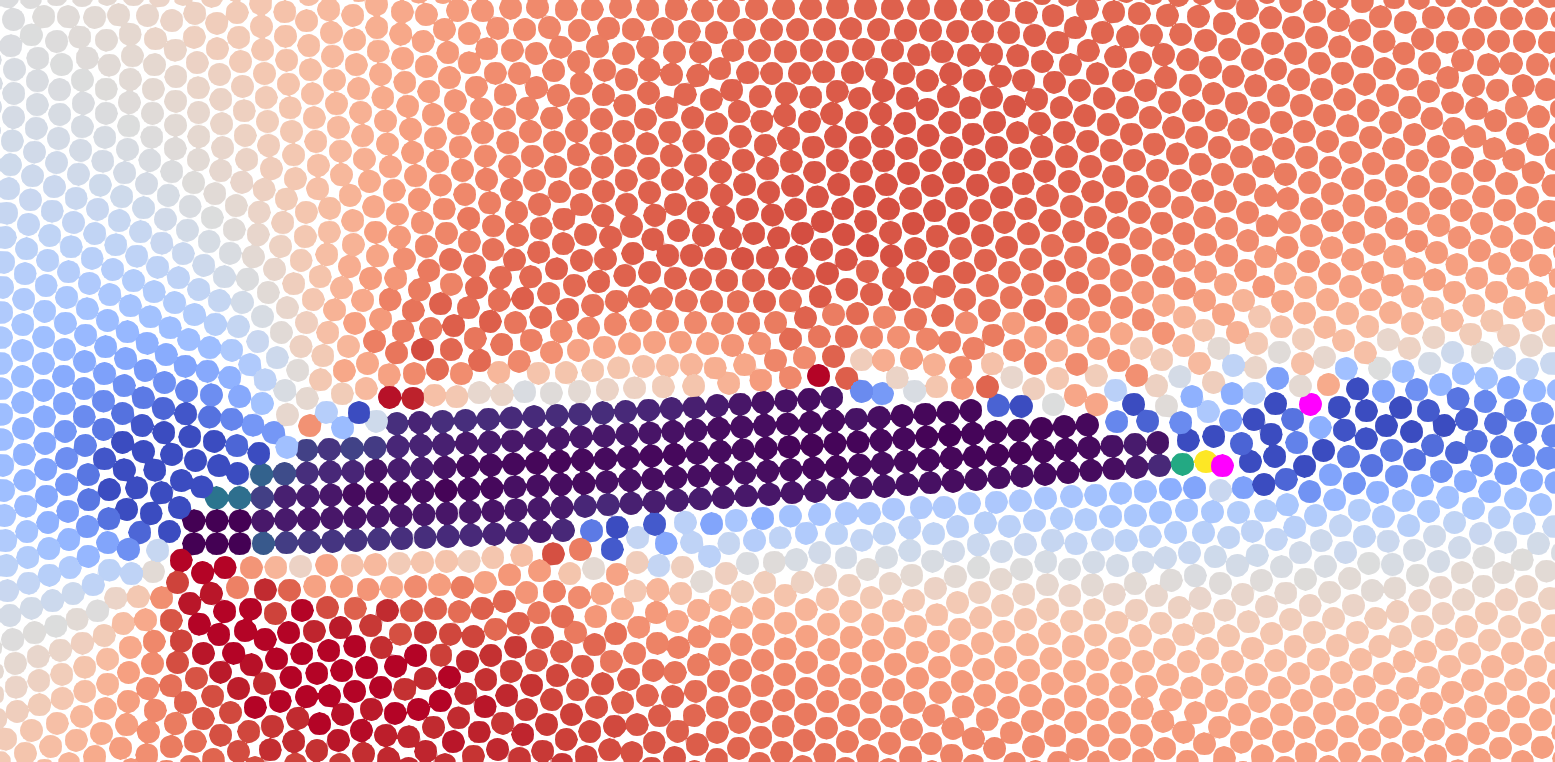}
    \caption{Lattice configuration $\Delta x = C /50$.}%
  \end{subfigure}
  \hfill
  \begin{subfigure}[b]{0.48\textwidth}
    \centering
    \includegraphics[width=0.75\textwidth]{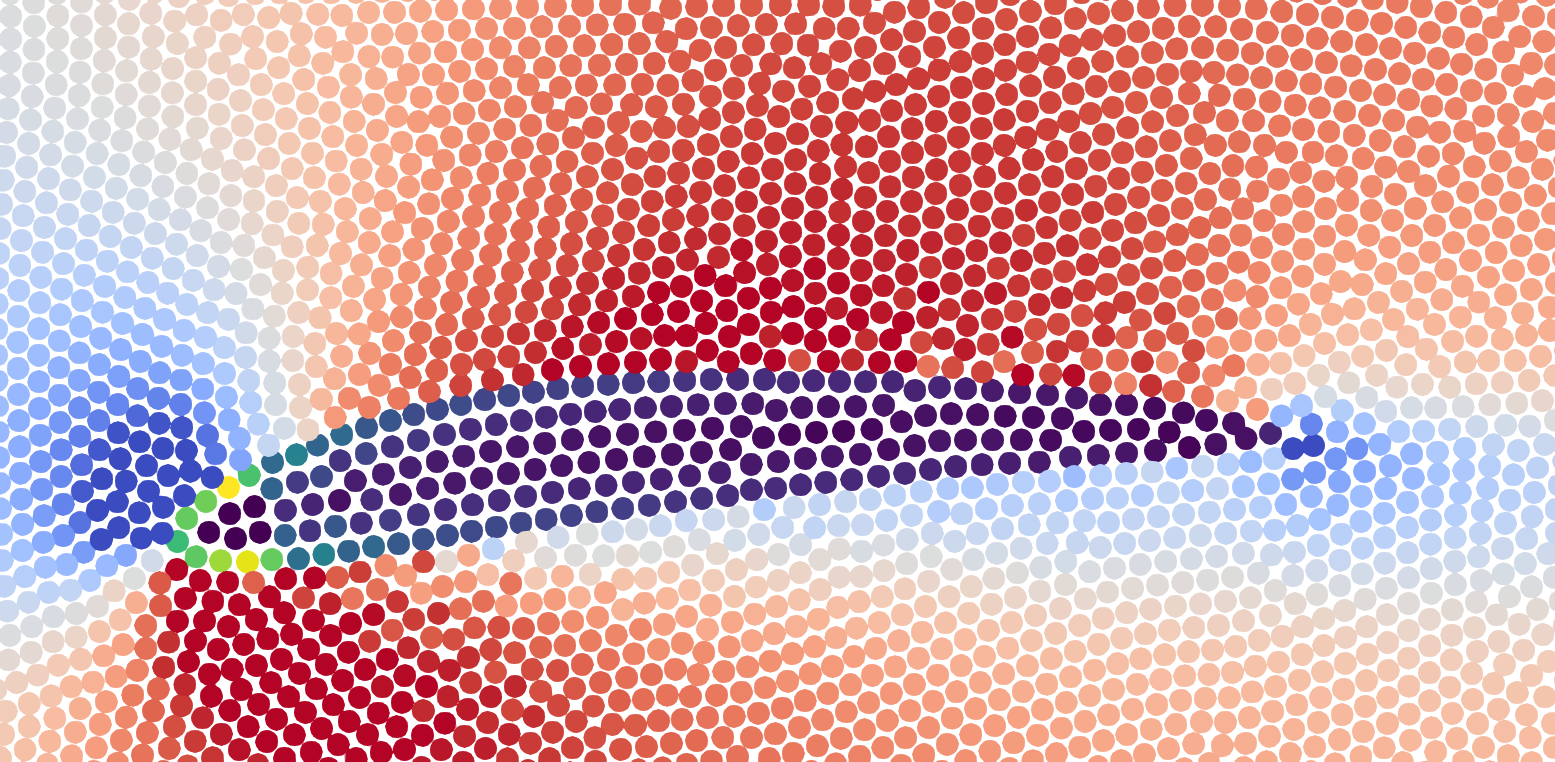}
    \caption{Packed configuration $\Delta x = C /50$.}%
  \end{subfigure}\\
  \begin{subfigure}[b]{0.48\textwidth}
    \centering
    \includegraphics[width=0.75\textwidth]{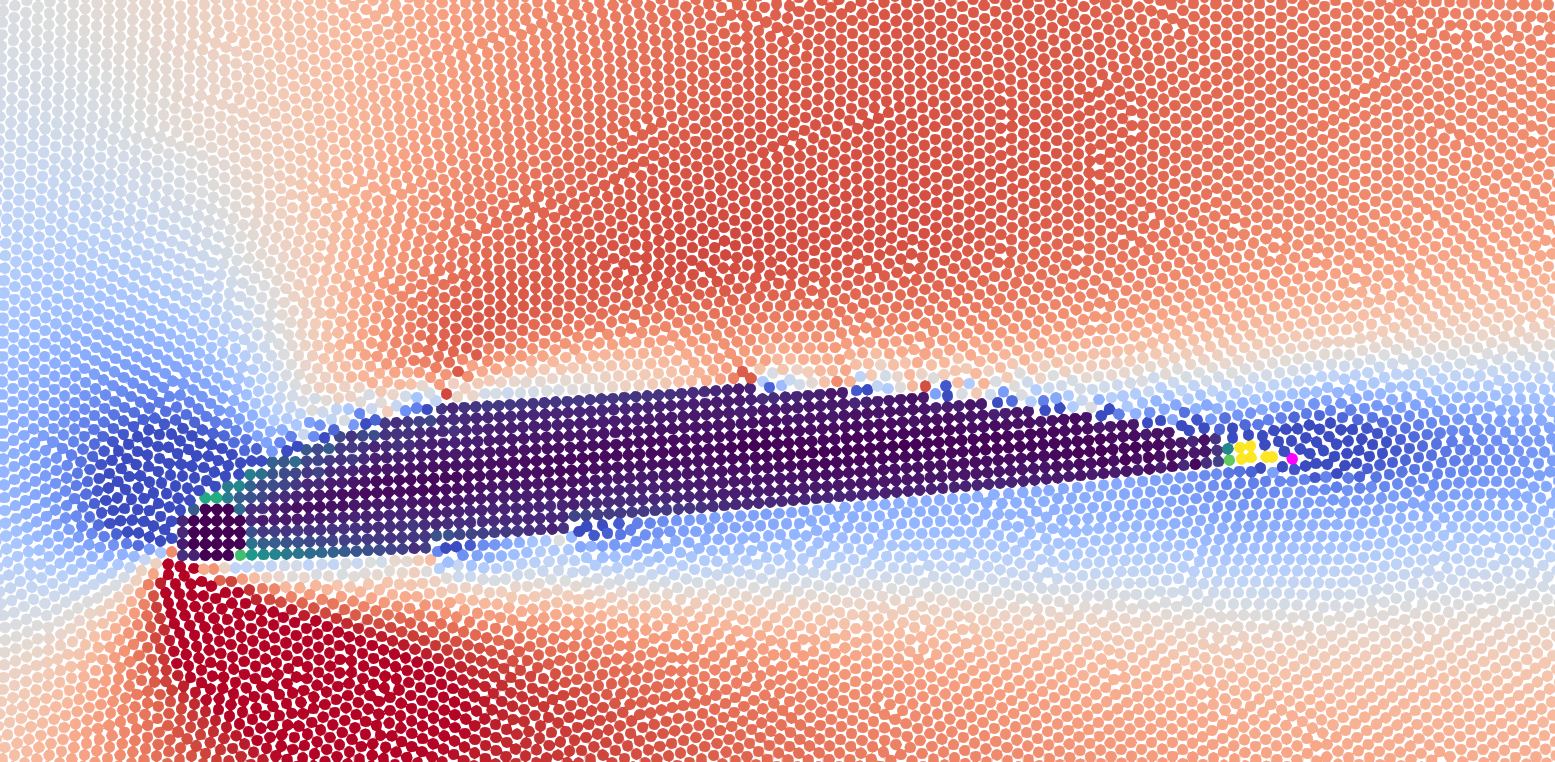}
    \caption{Lattice configuration $\Delta x = C /100$.}%
  \end{subfigure}
  \hfill
  \begin{subfigure}[b]{0.48\textwidth}
    \centering
    \includegraphics[width=0.75\textwidth]{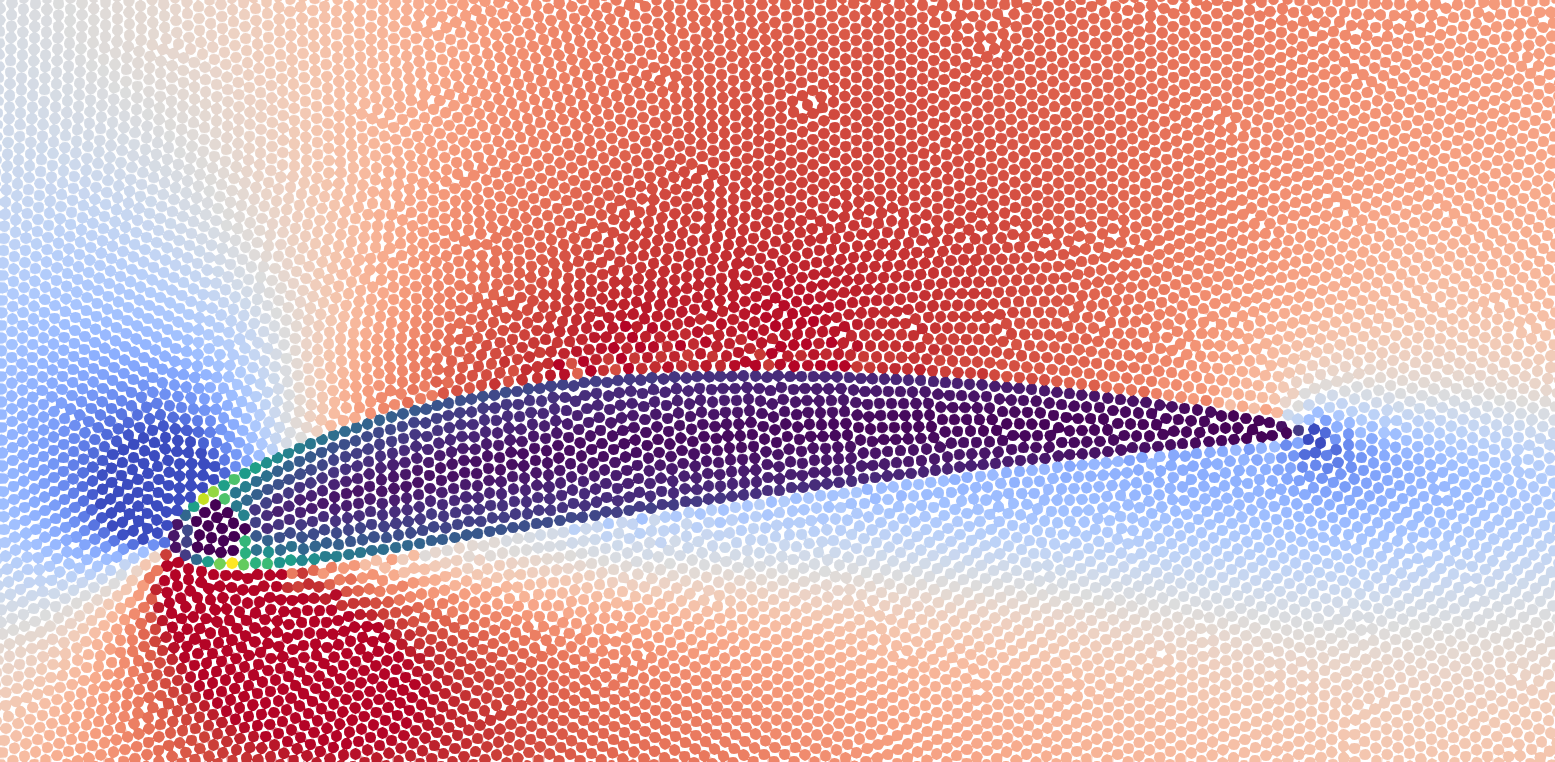}
    \caption{Packed configuration $\Delta x = C /100$.}%
  \end{subfigure}\\
  \begin{subfigure}[b]{0.48\textwidth}
    \centering
    \includegraphics[width=0.75\textwidth]{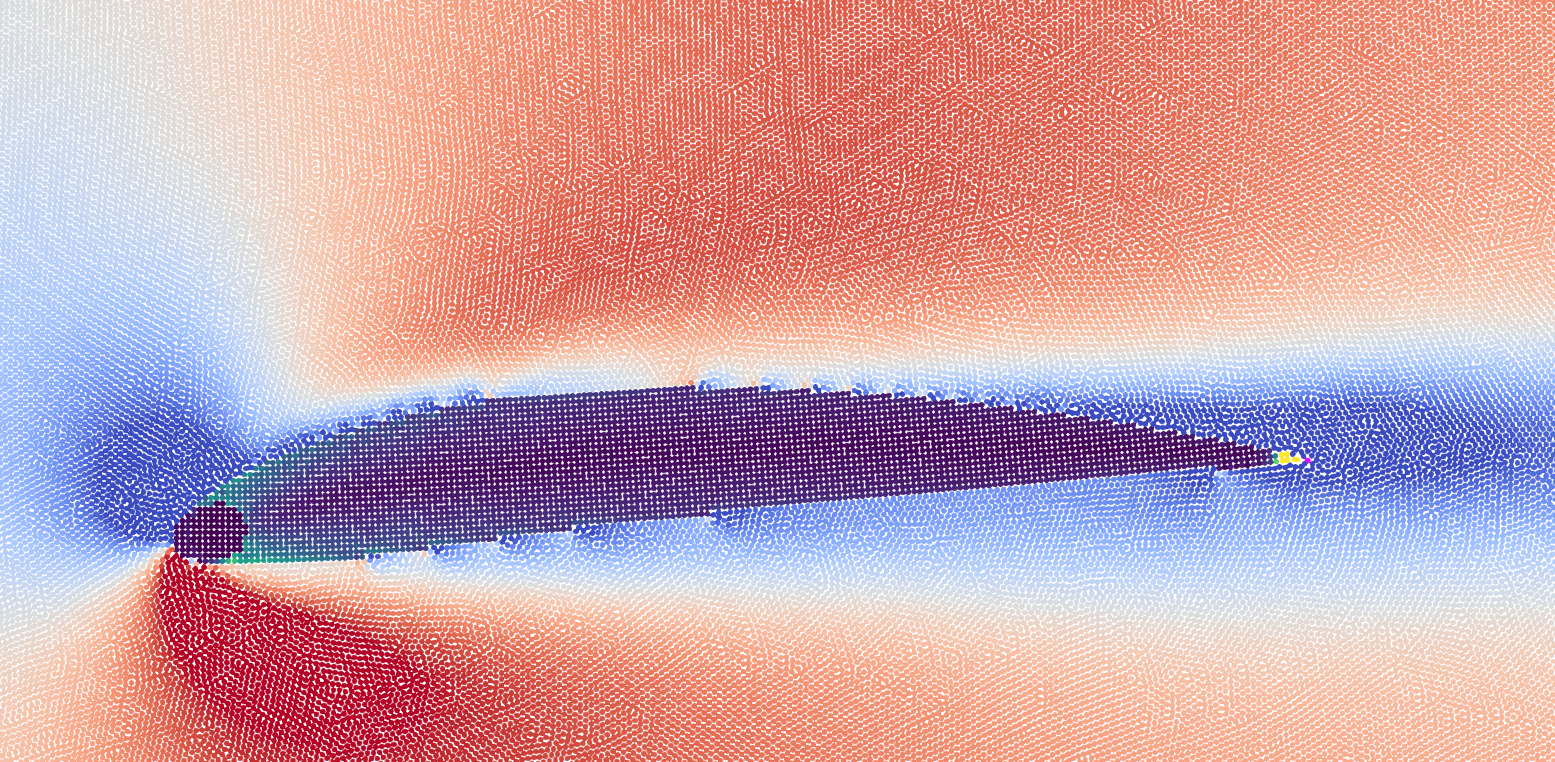}
    \caption{Lattice configuration $\Delta x = C /200$.}%
  \end{subfigure}
  \hfill
  \begin{subfigure}[b]{0.48\textwidth}
    \centering
    \includegraphics[width=0.75\textwidth]{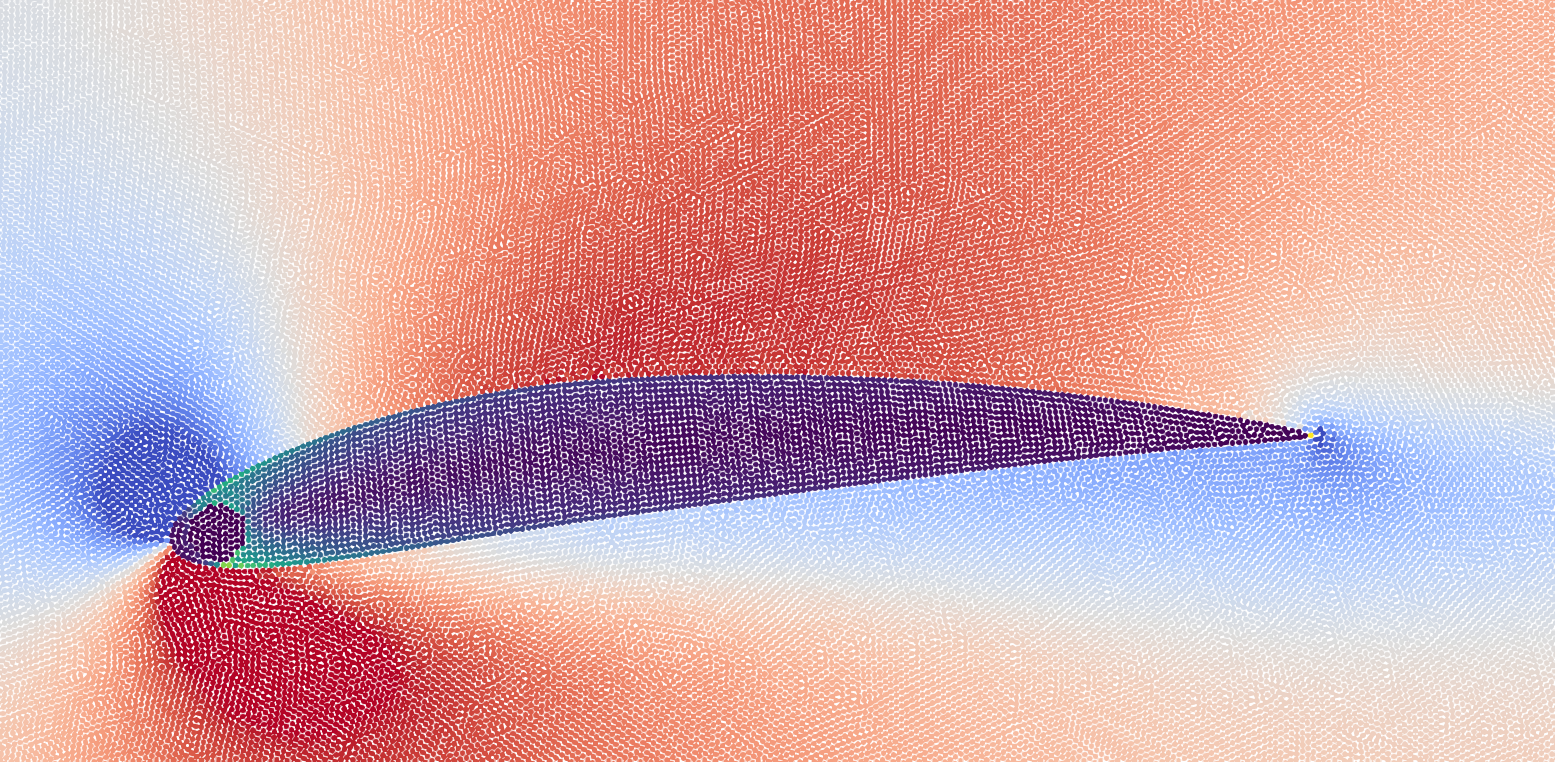}
    \caption{Packed configuration $\Delta x = C /200$.}%
  \end{subfigure}\\
  \begin{subfigure}[b]{0.48\textwidth}
    \centering
    \includegraphics[width=0.75\textwidth]{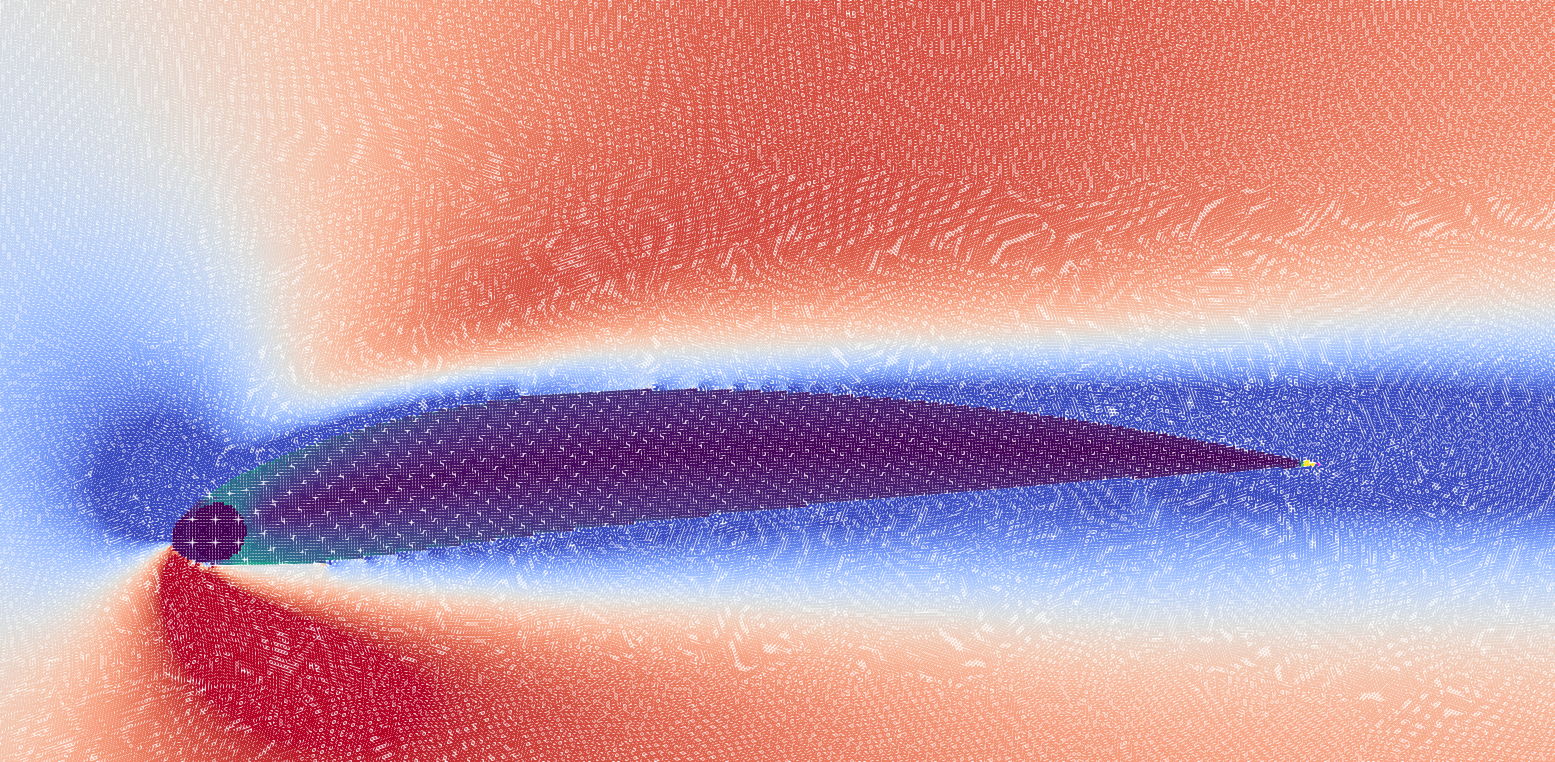}
    \caption{Lattice configuration $\Delta x = C /400$.}%
  \end{subfigure}
  \hfill
  \begin{subfigure}[b]{0.48\textwidth}
    \centering
    \includegraphics[width=0.75\textwidth]{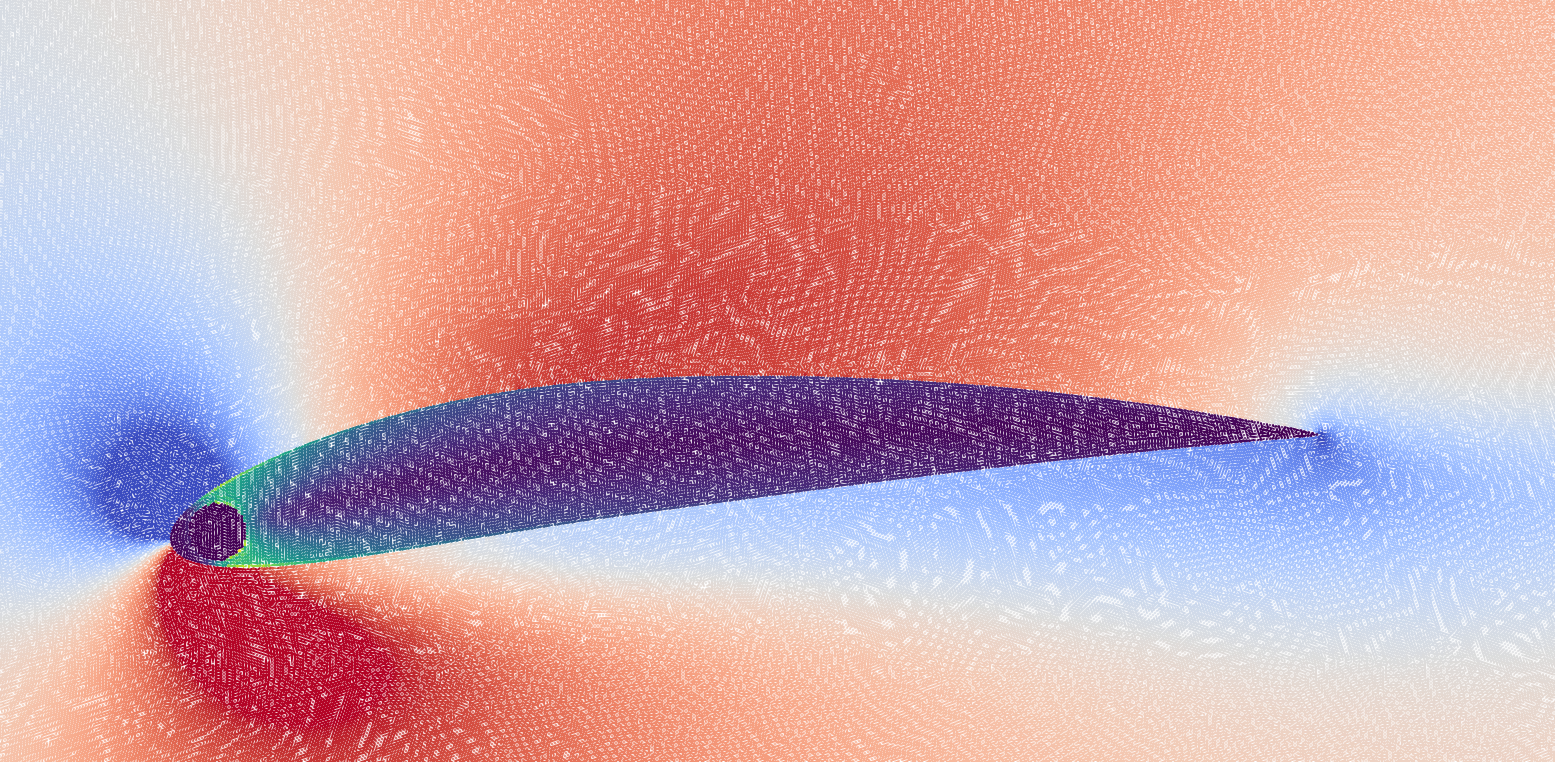}
    \caption{Packed configuration $\Delta x = C /400$.}%
  \end{subfigure}%
  \caption{Comparison of velocity magnitude in the flow around a NACA6412 airfoil as well as the von~Mises strain distribution within the airfoil structure for lattice-based (left) and packed (right) particle configurations at various resolutions at $t = 4.0$.}%
  \label{fig:solid_fluid_2}
\end{figure}

Fig.~\ref{fig:solid_fluid_1} provides a direct visual comparison of the pressure fields around the airfoil for both lattice-based (left) and packed (right) particle configurations at an early stage ($t=0.02$).
The lattice configuration exhibits pressure waves, which are especially pronounced at the trailing edge of the airfoil.
In contrast, the packed configuration establishes a smooth and physically consistent pressure profile, with no spurious oscillations present.
Fig.~\ref{fig:y_deflection} shows the time evolution of the vertical deflection of the trailing edge for different particle resolutions and configurations.
Both lattice-based (dashed line) and packed (solid line) particle configurations are compared to assess the influence of particle distribution on the dynamic response of the airfoil.
It is evident that the packed configuration results in a larger deflection. Even at the lowest resolution ($\Delta x = C/25$),
the stationary deflection is of similar magnitude as for higher resolutions.
In contrast, the lattice configuration produces a smaller deflection.
The reason for this discrepancy can be explained by examining the velocity profile around the airfoil:
Fig.~\ref{fig:solid_fluid_2} shows the velocity field at $t = 4.0$ for different resolutions and configurations (lattice and packed).
It can be observed that, in the lattice configuration, the fluid velocity near the airfoil surface is reduced due to the staircase-like structure of the particle arrangement.
This effect alters the dynamic response of the airfoil compared to the packed configuration,
where the fluid velocity at the surface remains largely unaffected, indicating a smoother interface.

To further evaluate the mechanical response of the airfoil, Fig.~\ref{fig:solid_fluid_2} also presents the von~Mises strain distribution
within the structure for different resolutions and for both lattice and packed configurations at $t=4.0$.
At low resolution, the lattice configuration exhibits non-physical strain concentrations at the trailing edge, where detached particles without any neighboring particles may occur (indicated by magenta).
In contrast, the packed configuration shows stress primarily localized at the clamped region, which is physically consistent.
Moreover, the stresses at the clamped region are higher in the packed case than in the lattice case, reflecting the larger deflection of the structure.

The results demonstrate that the packed configuration yields a dynamic response of the airfoil that accurately reflects the underlying geometry,
while---unlike the lattice-based configuration---no artificial stress concentrations at the trailing edge are observed.
Note that the geometry does not require any preprocessing or resolution-dependent modifications.

% !TeX root = ../main.tex
\section{Performance}
\label{sec:performance}
In this section, we analyze the runtime contributions of each preprocessing step in our workflow to understand their impact on overall performance.
The performance analysis was conducted on a single node with two AMD EPYC 7742 CPUs, each of which has 64 cores, and with a clock speed of 2.25~GHz.
For the performance analysis, we examine the following steps.
\begin{itemize}
  \item Load geometry: this involves reading the geometry's input data plus performing preprocessing.
  Specifically, ensuring unique vertex and edge IDs to facilitate the efficient identification of exterior faces, which is required for closing the bounding box in the subsequent step.
  \item Construct bounding box hierarchy: as detailed in Sec.~\ref{sec:hierarchical_winding}, a hierarchical structure of axis-aligned bounding boxes is created for fast winding number evaluation.
  \item Sample the geometry: perform inside-outside segmentation (Sec.~\ref{sec:segmentation}) to generate the initial particle configuration.
  \item Generate SDF: this includes initializing the face-based NHS and computing the SDF using Alg.~\ref{alg:face-nhs} and Alg.~\ref{alg:sdf}.
\end{itemize}
Note that the actual particle packing depends on the specific software framework, thus its performance numbers likely vary for different implementations.
For the performance analysis, we utilize the Stanford Bunny geometry \cite{stanfordbunny} due to its complexity and well-documented characteristics.
The software tool Blender~\cite{blender} is used to generate variants of the geometry with different face counts,
allowing us to evaluate performance at various levels of geometric complexity.
The number of particles in each of the following results refers to the number of generated query points in the initial point grid of the bounding box of the geometry.

\begin{figure}[!h]
  \centering
  \begin{subfigure}[b]{0.65\textwidth}
    \centering
    \includegraphics[width=0.95\textwidth]{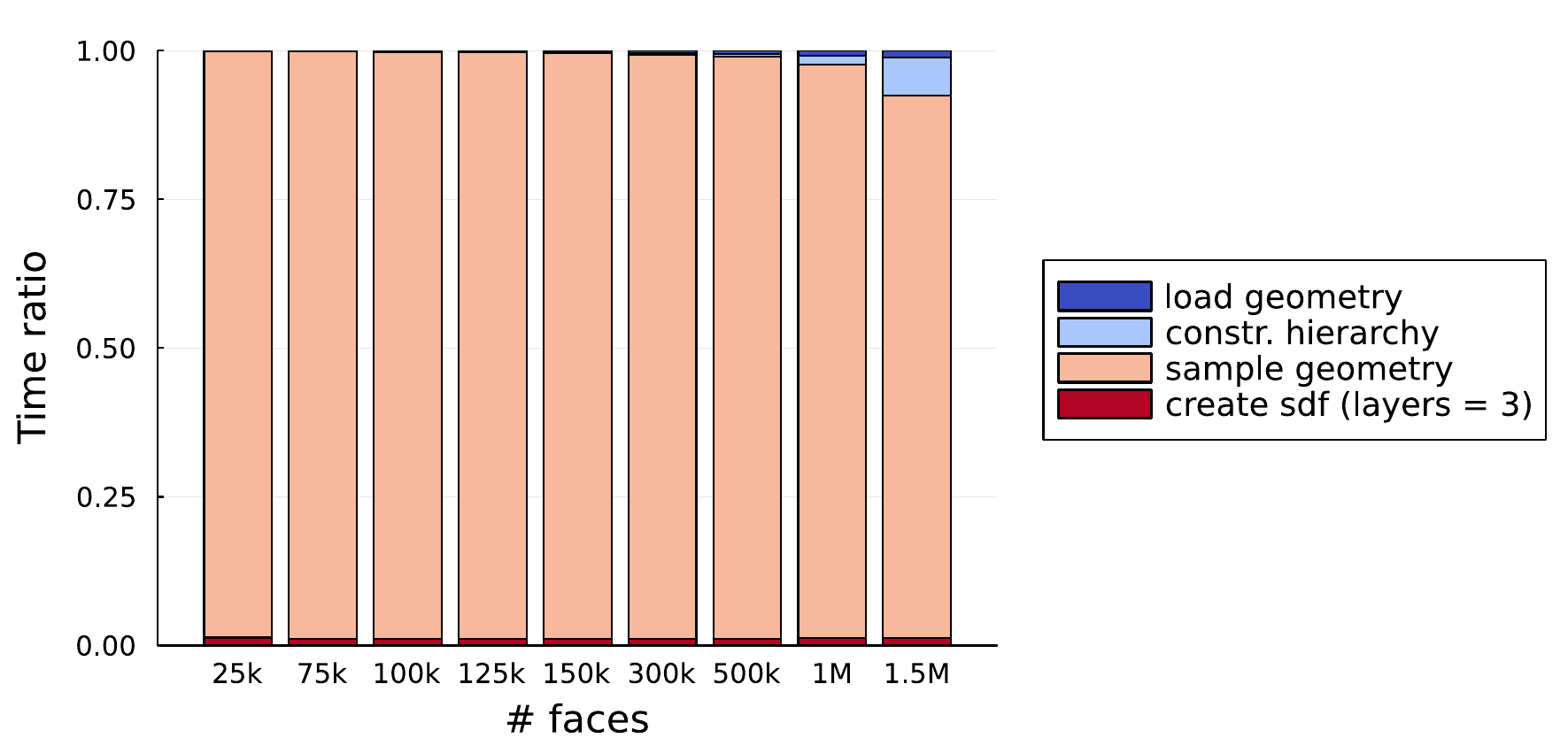}
    \caption{Fixed number of particles: $15.6\text{M}$.}%
    \label{fig:time_ratios_a}
  \end{subfigure}\\
  \begin{subfigure}[b]{0.65\textwidth}
    \centering
    \includegraphics[width=0.95\textwidth]{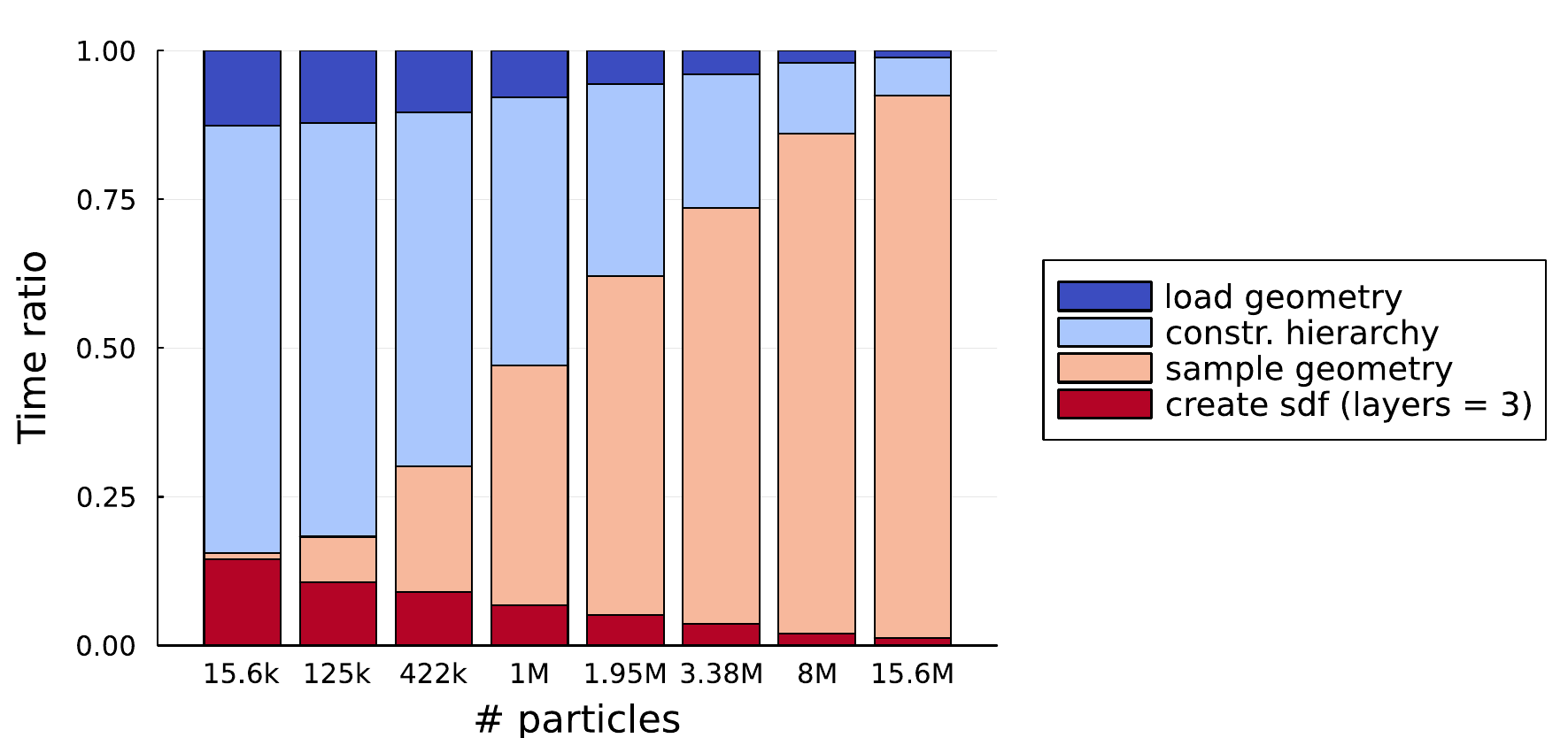}
    \caption{Fixed number of faces: $1.5\text{M}$.}%
    \label{fig:time_ratios_b}
  \end{subfigure}%
  \caption{Runtime contribution of the preprocessing steps executed in parallel on eight threads and with either a fixed number of particles (a) or faces (b).}
  \label{fig:time_ratios}
\end{figure}
To gain an overview of the runtime contributions of each individual step,
we analyze the proportions for scenarios with a relatively high fixed number of particles compared to a relatively high fixed number of faces.
Fig.~\ref{fig:time_ratios} illustrates the runtime contributions of each step.
Fig.~\ref{fig:time_ratios_a} shows the performance for a fixed particle count of approximately 15.6~million particles,
whereas Fig.~\ref{fig:time_ratios_b} varies the particle count for a fixed face count of approximately 1.5~million faces.
It is evident that for smaller numbers of particles, the time required for hierarchy construction and SDF creation dominates.
However, as the particle count increases, sampling becomes the dominant factor.
Note that the construction of the bounding box hierarchy depends solely on the geometry,
specifically the number of faces, and is independent of the particle count.
That is, in Fig.~\ref{fig:time_ratios_b} the time for creating the bounding box hierarchy remains constant regardless of the number of particles.

\subsection{Sampling the geometry}
\begin{figure}[!h]
  \centering
  \begin{subfigure}[b]{0.48\textwidth}
    \centering
    \includegraphics[width=0.95\textwidth]{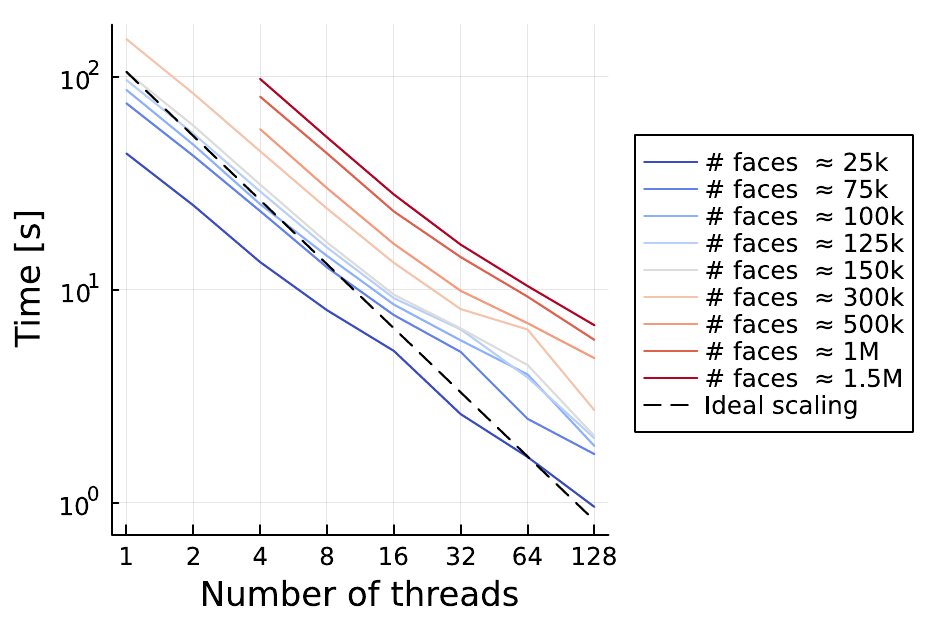}
    \caption{Fixed number of particles: $1\text{M}$.}%
    \label{fig:performance_sample_geometry_a}
  \end{subfigure}
  \hfill
  \begin{subfigure}[b]{0.48\textwidth}
    \centering
    \includegraphics[width=0.95\textwidth]{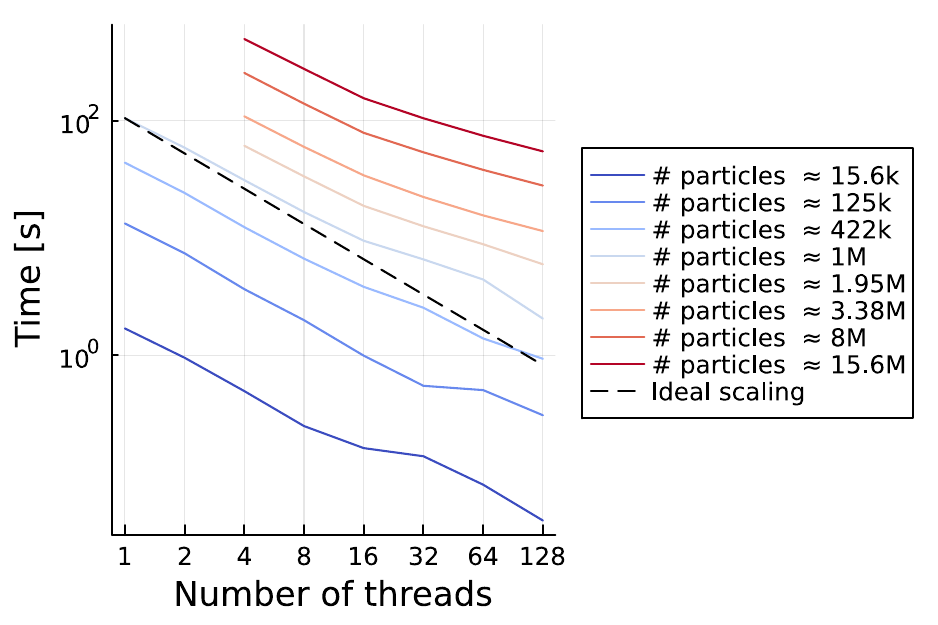}
    \caption{Fixed number of faces: $150\text{k}$.}%
    \label{fig:performance_sample_geometry_b}
  \end{subfigure}%
  \caption{Runtime for creating the initial particle contribution of the bunny geometry with different parameter variations plotted over the number of threads. The ideal scaling is shown as a dashed line.
  (a) shows the runtimes with a fixed number of particles and the ideal scaling for $150\text{k}$ faces.
  (b) shows the runtimes with a fixed number of faces and the ideal scaling for $1\text{M}$ particles.}
  \label{fig:performance_sample_geometry}
\end{figure}

Fig.~\ref{fig:performance_sample_geometry} shows the runtimes for sampling the geometry with different parameter variations as a function of the number of threads.
For one and two threads, no performance measurements were conducted for larger problem sizes, as the expected runtime would be disproportionately high and not representative of practical usage scenarios.
In Fig.~\ref{fig:performance_sample_geometry_a}, the number of particles is held constant, while in Fig.~\ref{fig:performance_sample_geometry_b}, the number of faces is fixed.
The ideal scaling (dashed line) is determined based on the fixed particle count or the fixed face count, respectively.
The observed scaling behavior nearly follows the trend of the corresponding ideal scaling.
However, deviations can be attributed to the fact that the bounding box hierarchy is not perfectly balanced, as the geometry is not symmetric.%

\subsection{Creating the SDF}
\label{sec:performance_sdf}
In the following, we analyze the runtime required for generating the SDF.
Key performance-related parameters in this step include the number of faces and the particle spacing $\Delta x$,
which is defined relative to the length of the cuboid bounding box of the geometry, denoted by $L$.
The particle spacing directly affects the search radius used in the face-based NHS,
and thus also influences the number of faces per cell.
To better understand this dependency, we first compare the runtime for SDF construction both with and without the use of the NHS.

\begin{figure}[!h]
  \centering
  \begin{subfigure}[b]{0.48\textwidth}
    \centering
    \includegraphics[width=0.65\textwidth]{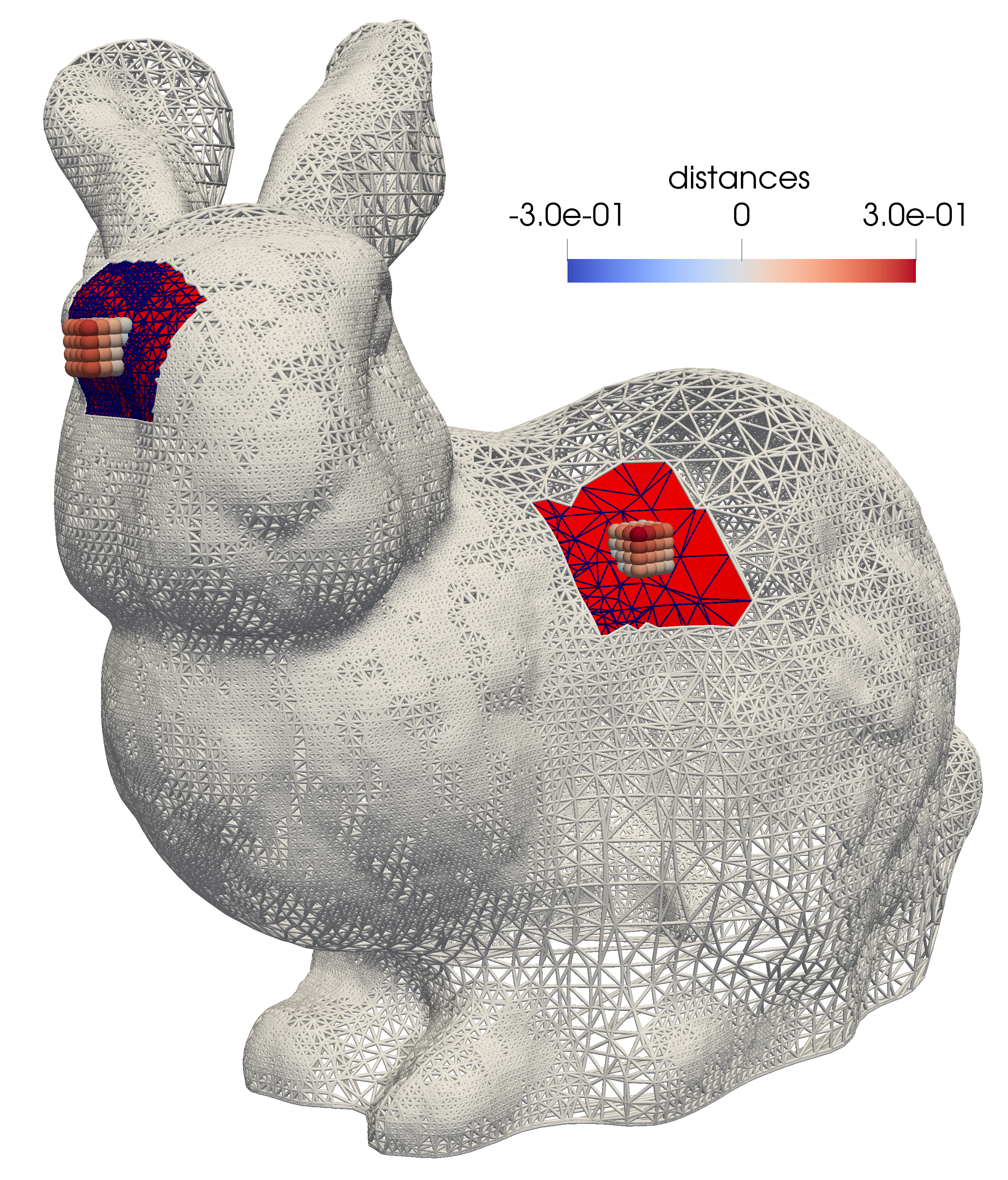}
    \caption{Max signed distances: $4 \Delta x$, where $\Delta x = L/100$}
    \label{fig:nhs_faces_bunny_a}
  \end{subfigure}
  \hfill
  \begin{subfigure}[b]{0.48\textwidth}
    \centering
    \includegraphics[width=\textwidth]{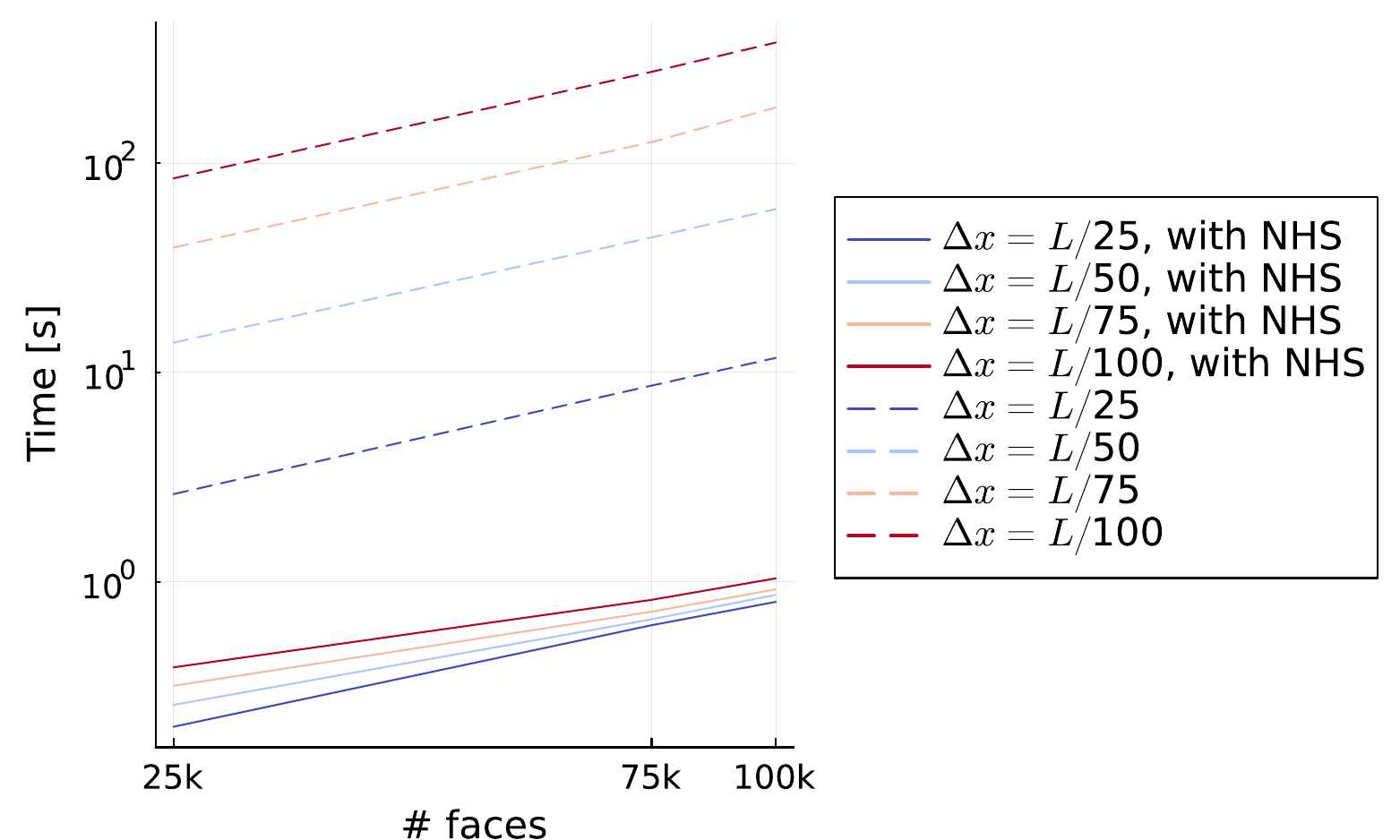}
    \caption{Runtimes for creating the SDF}
    \label{fig:nhs_faces_bunny_b}
  \end{subfigure}
  \caption{Right: Runtimes creating the SDF with 8 threads on an AMD EPYC 7742 node. The dashed lines represent the duration without NHS and the solid with NHS.
  Left: The geometry shows the SDF points for two exemplary NHS cells and the corresponding faces in the cell represented in red.
  $L$ corresponds to the bounding box size of the geometry.}%
  \label{fig:nhs_faces_bunny}
\end{figure}%
\begin{figure}[!h]
  \center
    \includegraphics[width=0.55\textwidth]{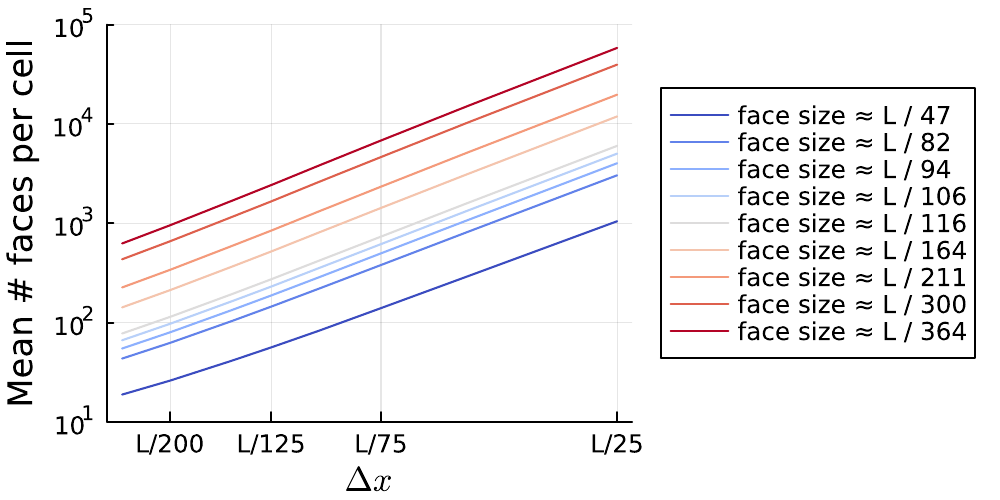}
  \caption{Average number of faces per cell in the face-NHS for a search radius of $3\Delta x$. The legend shows the average size of the faces in the domain.
  $L$ is the length of the cuboid bounding box of the geometry. The examined geometry is the Stanford Bunny.}
  \label{fig:face_nhs_mean_nfaces}
  \end{figure}%
Fig.~\ref{fig:nhs_faces_bunny} illustrates the performance of creating SDF with and without the use of the face-based NHS.
On the left (Fig.~\ref{fig:nhs_faces_bunny_a}), the geometry of the Stanford Bunny is shown, highlighting two exemplary NHS cells with their corresponding faces (red) and SDF points.
The figure on the right (Fig.~\ref{fig:nhs_faces_bunny_b}) presents the runtime for creating the SDF as a function of the number of faces for different resolutions.
Dashed lines depict the runtime without NHS, whereas solid lines represent the runtime with NHS, including its initialization.
The diagram shows a significant speedup when using the NHS.
It can also be seen that with an increasing number of faces, the runtime seems to converge.
The dominant factor influencing runtime in the face-based NHS is the ratio between the size of the mesh faces, which in this case correlates with the number of faces, and the resolution $\Delta x$.
Coarser resolutions result in a larger search radius, which leads to more faces per cell. Conversely, if the faces are large and the resolution is small,
the  NHS becomes more finely resolved, resulting in a higher number of cells.
This behavior is illustrated in Fig.~\ref{fig:face_nhs_mean_nfaces},
which shows the average number of faces per cell for a fixed search radius of $3\Delta x$ as a function of the resolution.
In this case, the variation in face size for each line is minimal, since the Stanford bunny geometry was remeshed in Blender.
It is clear that the number of faces per cell increases with coarser resolution (increasing search radius).
In summary, when the resolution is relatively fine compared to the face size, the runtime becomes more sensitive to $\Delta x$.
Conversely, if the resolution is coarse relative to the face size, the impact of $\Delta x$ on runtime diminishes.
However, if this ratio becomes too large, i.e., very coarse resolution combined with very small faces, the number of faces per cell increases significantly, leading to higher runtimes.
This effect can also be observed in the following performance diagrams.
Since the resolution $\Delta x$ is usually smaller than the size of the faces, in practice no further optimizations are necessary.

\begin{figure}[!h]
  \centering
  \begin{subfigure}[b]{0.48\textwidth}
    \centering
    \includegraphics[width=0.95\textwidth]{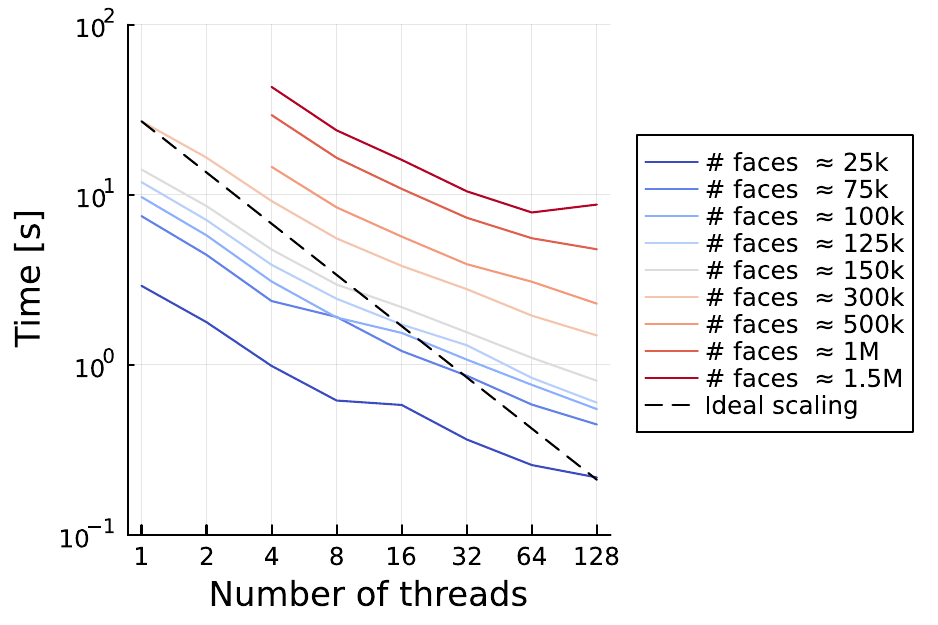}
    \caption{$\Delta x = L/100 $}%
    \label{fig:performance_create_sdf_a}
  \end{subfigure}
  \hfill
  \begin{subfigure}[b]{0.48\textwidth}
    \centering
    \includegraphics[width=0.95\textwidth]{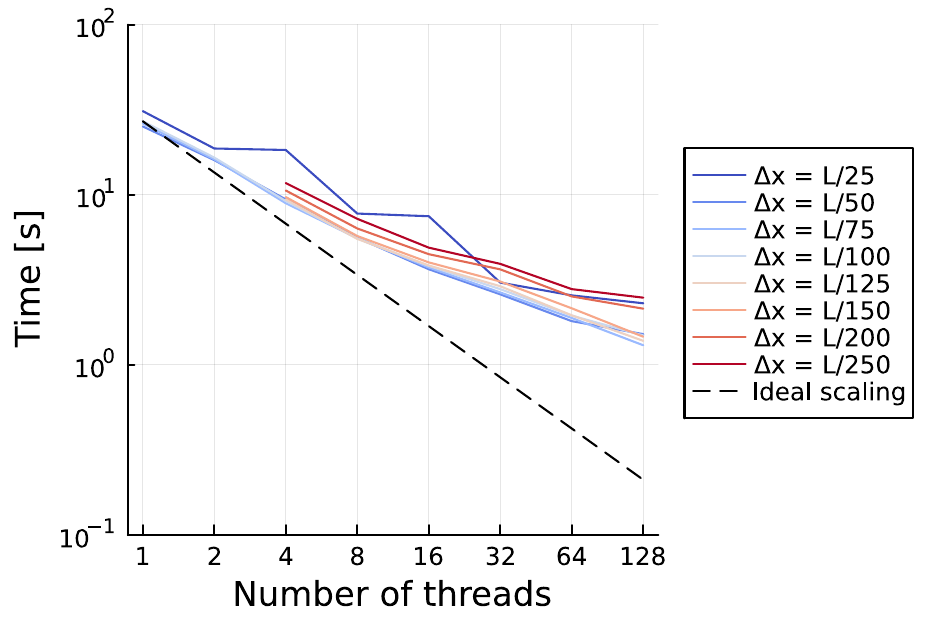}
    \caption{Number of faces: $\approx 300\text{k}$}%
    \label{fig:performance_create_sdf_b}
  \end{subfigure}\\
  \begin{subfigure}[b]{0.48\textwidth}
    \centering
    \includegraphics[width=0.95\textwidth]{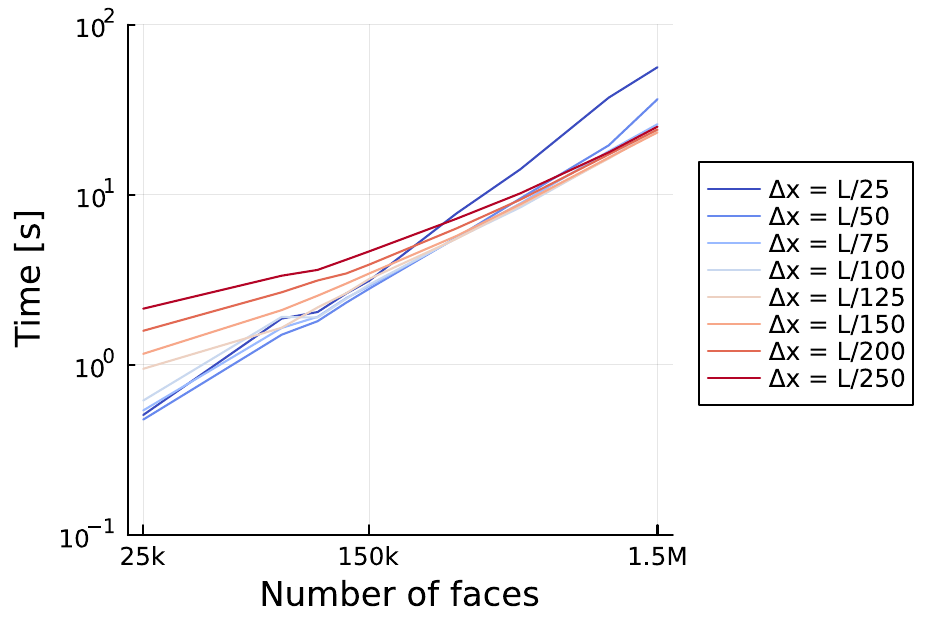}
    \caption{Number of threads: $8$}%
    \label{fig:performance_create_sdf_c}
  \end{subfigure}%
  \hfill
  \begin{subfigure}[b]{0.48\textwidth}
    \centering
    \includegraphics[width=0.95\textwidth]{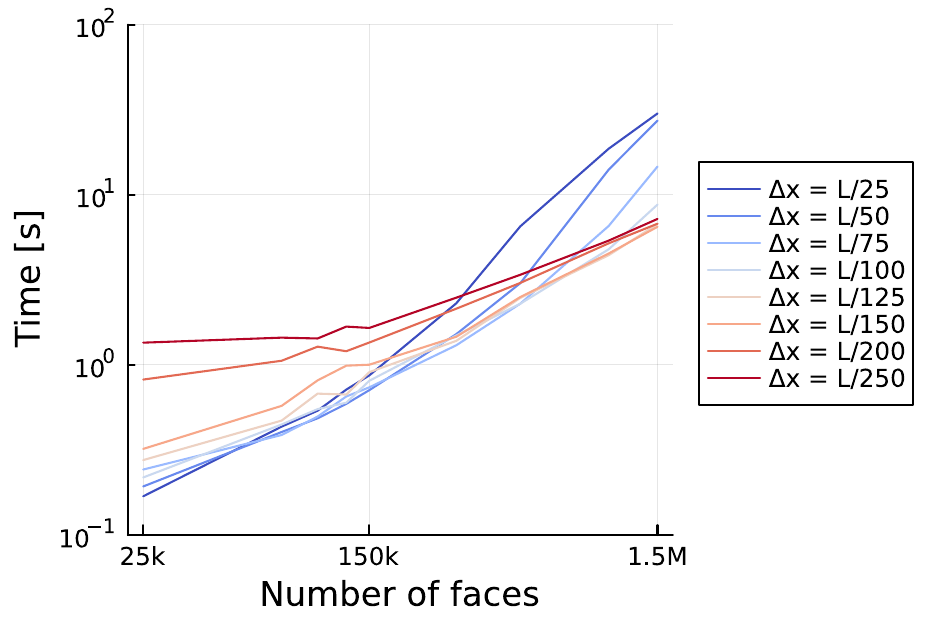}
    \caption{Number of threads: $128$}%
    \label{fig:performance_create_sdf_d}
  \end{subfigure}%
  \caption{Performance plots for creating the signed distance field with different parameter variations. (a) and (b) show the runtimes plotted over the number of threads.
            The dashed lines represent the ideal scaling for $\Delta x = L/100 $ particles and $3 \cdot 10^5$~faces.
          (c) and (d) show the runtimes plotted over the number of faces. $L$ is the length of the cuboid bounding box of the geometry.}
  \label{fig:performance_create_sdf}
\end{figure}%
Fig.~\ref{fig:performance_create_sdf} shows the runtimes for creating the SDF with different parameter variations.
The following runtime measurements were conducted with a search radius of $3\Delta x$.
The first row shows the runtimes as a function of the number of threads. The dashed lines represent the ideal scaling for approximately $3 \cdot 10^5$ faces and a resolution of $\Delta x = L / 100$.
In Fig.~\ref{fig:performance_create_sdf_a} we see that for constant $\Delta x$ the number of faces influences the runtime,
while Fig.~\ref{fig:performance_create_sdf_b} indicates that $\Delta x$ has no significant impact on the runtime.
This behavior is also evident when looking at the second row in Fig.~\ref{fig:performance_create_sdf},
where the runtime is over the number of faces and with a fixed number of threads.
We observe that as the number of faces increases (i.e., as the individual faces become smaller), the influence of $\Delta x$ on the runtime diminishes,
which is related to the face-based NHS, as described earlier.
The scaling behavior is not perfectly ideal because the current implementation only parallelizes the loop for SDF computation (line 4 in Alg.~\ref{alg:sdf}) and the removal of points in empty cells (line~3 in Alg.~\ref{alg:sdf}).
Additionally, this could also be related to the imperfect balance of the number of faces per cell.%

The runtime analysis shows that the majority of computational cost is associated with sampling the geometry, which strongly depends on the chosen resolution.
However, this step has been shown to scale efficiently with the number of threads.
In contrast, the runtime for computing the SDF is primarily influenced by the search radius.
For a practical search radius, typically in the order of the geometry's characteristic face size,
the use of the face-based NHS leads to a significant speedup.
While the parallel scalability of the SDF generation is not ideal,
its runtime contribution remains relatively small compared to the geometry sampling step and is thus acceptable for practical applications.

% !TeX root = ../main.tex
\section{Conclusions}
We presented an automated sampling and packing technique for SPH and other particle-based methods that is robust and efficient,
yielding high-quality particle distributions.
A key contribution is the construction of a memory-efficient signed distance field (SDF) through a face-based neighborhood search,
which significantly accelerates the SDF point cloud generation.
This approach also localizes computations to surface regions, enhancing performance during the packing process.
The combined use of a resolution-dependent SDF and the winding number method ensures that the results converge towards the exact geometry.
Additionally, boundary particles can be generated and packed dynamically, improving the particle distribution and ensuring a high-quality boundary representation through particles.
Including dynamic boundary particles also enhances the robustness of our method, facilitating rapid convergence to a steady state and is further supported by adaptive time-step control.
The effectiveness of our approach is further demonstrated by a simulation with fluid-structure interaction, which underlines the practical applicability and robustness of the method.
Since our method is entirely particle-based and does not rely on a background mesh,
it can be seamlessly integrated into existing SPH frameworks, making it highly versatile.
A comprehensive performance analysis of the preprocessing reveals potential for effective multi-threading scalability,
making our approach well-suited for large-scale problems.
Looking ahead, a detailed analysis of GPU performance is planned.
Furthermore, we aim to develop a multiresolution packing technique and introduce a scale separation method for the SDF point cloud to resolve sharp edges and corners.

% =========================================================================================
% ==== Data availability
\section*{Data availability}\label{sec:data_availability}

All code used in this work and the data generated during this study are available via our
reproducibility repository \cite{Neher2025reproducibility}.

% =========================================================================================
% ==== CRediT
\section*{CRediT authorship contribution statement}\label{sec:credit_authorship}

\textbf{Niklas S. Neher}:
Conceptualization, Methodology, Investigation, Data curation, Visualization, Software, Validation,
Writing --- original draft, Writing --- review \& editing.
\textbf{Erik Faulhaber}:
Software, Validation, Writing --- review \& editing.
\textbf{Sven Berger}:
Software, Validation, Writing --- review \& editing.
\textbf{Christian Weißenfels}:
Writing --- review \& editing.
\textbf{Gregor J. Gassner}:
Conceptualization, Writing --- review \& editing.
\textbf{Michael Schlottke-Lakemper}:
Conceptualization, Writing --- review \& editing.

% =========================================================================================
\section*{Declaration of competing interest}\label{sec:competing_interests}
The authors declare that they have no known competing financial interests or personal
relationships that could have appeared to influence the work reported in this paper.

% =========================================================================================
\section*{Declaration of generative AI and AI-assisted technologies in the writing
process}\label{sec:declaration-ai}

During the preparation of this manuscript the authors used ChatGPT to enhance
readability. They further used the GitHub Copilot plugin for Visual Studio Code for
assistance with mathematical notation and formatting in LaTeX. After using these
tools/services, the authors reviewed and edited the content as needed and take full
responsibility for the content of the published article.

% =========================================================================================
\section*{Acknowledgments}\label{sec:acknowledgments}
We gratefully acknowledge the computing time on Hawk provided by the High-Performance Computing Center
Stuttgart (HLRS), University of Stuttgart.
Gregor J. Gassner and Michael Schlottke-Lakemper acknowledge funding through the DFG research unit FOR 5409 ``Structure-Preserving Numerical Methods for Bulk- and Interface Coupling of Heterogeneous Models (SNuBIC)''
and through the BMBF project ``ADAPTEX''.
Michael Schlottke-Lakemper further received funding through a DFG individual grant (project number 528753982).
We thank Alec Jacobson for the helpful discussions on closing surfaces in hierarchical winding number computations.

\bibliography{references}

\end{document}